\pdfoutput=1
\documentclass[a4paper,11pt]{amsart}
\usepackage[english]{babel}
\usepackage{amsfonts}
\usepackage{amsthm}
\usepackage{amsmath}
\usepackage{amssymb}
\usepackage{enumerate}
\usepackage{mathtools}
\usepackage{tikz}
\usepackage{tikz-cd}
\usepackage[all]{xy}
\usepackage{graphicx}
\usepackage{amsmath}
\usepackage{adjustbox}
\usepackage{graphicx}
\usepackage{extpfeil}
\usepackage{comment}
\usepackage{float}
\usepackage[labelformat=empty]{caption}
\usepackage[left=3.2cm,top=3cm,right=3.2cm,bottom=3cm]{geometry}
\usepackage[toc]{appendix}
\usepackage{hyperref}
\hypersetup{
colorlinks=true,
linkcolor=[RGB]{0,0,204},
filecolor=magenta,      
urlcolor=cyan,
citecolor=[RGB]{0,204,0},
}
\usepackage{resizegather}
\usepackage[activate={true, nocompatibility}, final, tracking=true, kerning=true, spacing=true, factor=1100, stretch=10, shrink=10]{microtype}
\usepackage{tikz}
\usepackage{dsfont}
\usetikzlibrary{matrix}

\usepackage{color}
\usepackage{resizegather}
\usepackage[nameinlink]{cleveref}
\usepackage{ulem}
\usepackage{enotez}
\setenotez{backref=true}

\usepackage{setspace}
\usepackage{circuitikz}

\usepackage{pbox}

\theoremstyle{plain}
\newtheorem{thm}{Theorem}[section]
\newtheorem{coroll}[thm]{Corollary}
\newtheorem{lemma}[thm]{Lemma}
\newtheorem{prop}[thm]{Proposition}

\newtheorem*{thm*}{Main Theorem}
\newtheorem{introthm}{Theorem}

\theoremstyle{definition}
\newtheorem{defn}[thm]{Definition}

\newtheorem{notn}[thm]{Notation}
\newtheorem{context}[thm]{Context}

\newtheorem{remark}[thm]{Remark}

\newcommand{\dash}{{\operatorname{-}}}

\tikzset{
  symbol/.style={
    draw=none,
    every to/.append style={
      edge node={node [sloped, allow upside down, auto=false]{$#1$}}}
  }
}

\microtypecontext{spacing=nonfrench}

\makeatletter 
\def\makeCal#1{%
\expandafter\newcommand\csname c#1\endcsname{\mathcal{#1}}}

\def\makeBB#1{%
\expandafter\newcommand\csname b#1\endcsname{\mathbb{#1}}}

\def\makeFrak#1{%
\expandafter\newcommand\csname f#1\endcsname{\mathfrak{#1}}}

\def\makeScr#1{%
\expandafter\newcommand\csname s#1\endcsname{\mathscr{#1}}}

\count@=0
\loop
\advance\count@ 1
\edef\y{\@Alph\count@} 
\expandafter\makeCal\y
\expandafter\makeBB\y
\expandafter\makeFrak\y
\expandafter\makeScr\y
\ifnum\count@<26
\repeat

\def\makelowercaseFrak#1{%
\expandafter\newcommand\csname mf#1\endcsname{\mathfrak{#1}}}
\count@=0
\loop
\advance\count@ 1
\edef\y{\@alph\count@} 
\expandafter\makelowercaseFrak\y
\ifnum\count@<26
\repeat

\DeclareMathOperator{\GL}{GL}
\DeclareMathOperator{\Bun}{Bun}
\DeclareMathOperator{\ev}{ev}

\DeclareMathOperator{\wt}{wt}

\DeclareMathOperator{\Cor}{Cor}

\DeclareMathOperator{\Dqc}{D_{qc}}

\DeclareMathOperator{\Perf}{Perf}

\DeclareMathOperator{\Grad}{Grad}

\DeclareMathOperator{\gr}{gr}
\DeclareMathOperator{\Map}{\operatorname{Map}}

\DeclareMathOperator{\Gr}{Gr}
\DeclareMathOperator{\transfer}{\operatorname{\Omega}}

\newcommand{\B}{\mathcal{B}}

\newcommand{\SVect}{\mathcal{V}ect}

\newcommand{\on}{\operatorname}

\newcommand{\Coh}{\on{Coh}}

\newcommand{\pol}{\on{pol}}

\newcommand{\Aut}{ \on{Aut} } 

\newcommand{\Rep}{\on{Rep}}
\newcommand{\Hom}{ \on{Hom}}
\newcommand{\Ext}{ \on{Ext}}

\renewcommand{\ker}{ \on{ker}}

\newcommand{\im}{ \on{Im}}

\newcommand{\Spec}{\on{Spec}}

\newcommand{\Pic}{\on{Pic}}

\newcommand{\BGL}{B \!\GL}

\newcommand{\rk}{\on{rk}}

\newcommand\dirac{/\kern-1.2ex\partial} 
\newcommand\qu{/\kern-.7ex/} 
\newcommand\lqu{\backslash \kern-.7ex \backslash} 

\newcommand\dr{r_+ \kern-.7ex - \kern-.7ex r_-}
\newcommand{\stabi}{c}

\tikzset{
  symbol/.style={
    draw=none,
    every to/.append style={
      edge node={node [sloped, allow upside down, auto=false]{$#1$}}}
  }
}

    \newtheoremstyle{TheoremNum}
        {\topsep}{\topsep}              
        {\itshape}                      
        {}                              
        {\bfseries}                     
        {.}                             
        { }                             
        {\thmname{#1}\thmnote{ \bfseries #3}}
    \theoremstyle{TheoremNum}

\DeclareMathOperator{\StMap}{\overline{\mathcal{M}}}
\DeclareMathOperator{\GBun}{GBun}
\DeclareMathOperator{\CBun}{\overline{GBun}}
\DeclareMathOperator{\Pur}{Pur}
\DeclareMathOperator{\ggw}{I}
\DeclareMathOperator{\tot}{tot}
\DeclareMathOperator{\GMap}{GMap}
\newcommand{\lev}{\mathrm{lev}}

\newcommand{\cal}[1]{{\mathcal #1}} \newcommand{\prestable}{\cal{P}}
\newcommand{\precont}{\cal{P}^{\mathrm{ cont}}}
\newcommand{\prestablepol}{\cal{P}^{\mathrm{\pol}}}

\author{Daniel Halpern-Leistner}
\address{Mathematics Department, Cornell University, 310 Malott Hall, Cornell
University, Ithaca, New York 14853, U.S.A.}
\email{daniel.hl@cornell.edu}

\author{Andres Fernandez Herrero}
\address{Department of Mathematics, University of Pennsylvania, 209 South 33rd
Street, Philadelphia, PA 19104, USA} 
\email{andresfh@sas.upenn.edu}

\begin{document}
\title{Quantum operations on the ring of symmetric functions}

\begin{abstract}
    We define a version of stable maps into the classifying stack
    $B\mathrm{GL}_N$, and develop a corresponding notion of $K$-theoretic
    Gromov-Witten invariants. In this setting, the evaluation morphisms are not
    of finite type; the definition of the $K$-theoretic invariants proceeds by
    constructing a $\Theta$-stratification of the moduli stack. In the absence
    of markings, the semistable locus of the stratification recovers moduli
    spaces of bundles on nodal curves considered by Gieseker, Nagaraj-Seshadri,
    Schmitt and Kausz.

    We also define versions of stable maps into quotient stacks of the form
    $Z/\mathrm{GL}_N$, where $Z$ is a projective $\mathrm{GL}_N$-scheme. We
    construct corresponding stability stratifications, whose semistable loci
    provide new proper moduli spaces of gauged maps from a varying nodal curve
    into $Z/\mathrm{GL}_N$.
\end{abstract}
\maketitle
\tableofcontents

\addtocontents{toc}{\protect\setcounter{tocdepth}{1}}

\section{Introduction}

Modern enumerative geometry has benefitted tremendously from the interaction of
two flavors of mathematics. On one side there is moduli theory, in which
researchers have developed sophisticated ways to count geometric
objects -- Gromov-Witten invariants, Donaldson invariants, PT invariants, DT invariants, etc. On the
other side, there are beautiful algebraic structures that give order to these
counting problems, often inspired by mathematical physics -- cohomological field
theories, BV algebras, Hall algebras, vertex algebras, chiral algebras, etc.

Inspiration has
flowed in both directions between these two kinds of mathematics. For instance,
the definition of a cohomological field theory inspired the rigorous
construction of Gromov-Witten invariants \cite{kontsevich-manin}. In many other instances, such as Donaldson-Thomas theory \cite{donaldson_thomas_gauged, thomas_casson_invariants},
the development of new counting problems enabled the development of rich algebraic structures -- in this case wall-crossing formulas and Hall algebras \cite{kontsevich_soibelman_2, joyce_song,  kontsevich_soibelman}.

This project sits on the moduli-theoretic side of enumerative geometry. We
identify a smooth, unbounded stack $\StMap_{g,n}(\BGL_N,d)$ that enlarges the stack
of rank $N$ and degree $d$ vector bundles on smooth genus $g$ curves with $n$
marked points. There is a flat stabilization morphism $\stabi: \StMap_{g,n}(\BGL_N,d) \to
\overline{\cM}_{g,n}$ to the stack of Deligne-Mumford stable curves, as well as an
evaluation-at-marked-points morphism $\ev : \StMap_{g,n}(\BGL_N,d) \to (\BGL_N)^n$. We were led to this moduli problem by the technical tools developed in the
beyond Geometric Invariant Theory paradigm for moduli theory \cite{torsion-freepaper, gauged_theta_stratifications}. The tools apply
very cleanly to this stack, and they
provide a semistable locus that admits a proper good moduli space, which
actually recovers a classical compactification due to Nagaraj-Seshadri \cite{nagaraj-seshadri-2} and Schmitt \cite{schmitt-hilbert-compactification}.

 The stack $\StMap_{g,n}(\BGL_N,d)$ also has a
$\Theta$-stratification with cohomologically proper centers. This allows us to
prove that for every $E \in \Perf((\BGL_N)^n$) and $V \in \Perf(\StMap_{g,n})$ the derived global sections $R\Gamma(\StMap_{g,n}(\BGL_N,d), \cL_{\lev} \otimes \stabi^\ast(V) \otimes \ev^\ast(E))$ are finite dimensional. A precise definition of this stack, and the ``level'' line bundle $\cL_{\lev}$ will be defined precisely later in this introduction. Taking virtual dimension gives a well-defined group homomorphism
\begin{equation}
K^0(\overline{\cM}_{g,n}) \otimes K^0(\BGL_N)^{\otimes n} \to \bZ.
\end{equation}
The insight that the $K$-theoretic index is
well-defined for certain perfect complexes on unbounded stacks originated in
\cite{teleman-woodward}. Our work here extends, from $\GL_1$ to $\GL_N$, a project initiated by
Teleman-Tolland-Frenkel \cite{frenkel-teleman-tolland-ggw} to use this idea to
define $K$-theoretic ``Gromov-Witten gauge theory.'' It is expected that the
homomorphisms contain information about the Gromov-Witten invariants of any GIT
quotient of a projective manifold by $\GL_N$.

To our surprise, though, an apparently new algebraic structure jumped out of
purely moduli-theoretic considerations. For the moduli of vector bundles on
nodal curves, we introduce a new ``boundary
condition'' at marked points on the
curve. We use this to define a more general algebraic stack
$\StMap_{g,n,p}(\BGL_N,d)$ parameterizing bundles on nodal curves with $n$
negative and $p$ positive marked points. The stack $\StMap_{g,n,p}(\BGL_N,d)$ now has two evaluation morphisms: $\ev^+$ to $\BGL_N^p$ and $\ev^-$ to $\BGL_N^n$. The new stack admits a $\Theta$-stratification \emph{relative} to the morphism
$\ev^+ \times \stabi : \StMap_{g,n,p}(\BGL_N,d) \to (\BGL_N)^p \times \StMap_{g,n+p}$, and the semistable locus and the centers of the unstable
strata are cohomologically proper over $(\BGL_N)^p \times \StMap_{g,n+p}$.

Our main result applies to the Fourier-Mukai functor
$\transfer^a_{g,n,p}(N,d) : \Dqc((\BGL_N)^n)  \to \Dqc(
(\BGL_N)^p \times \StMap_{g,n+p})$ on derived categories of quasi-coherent sheaves for any $a>0$ given
by
\[
    \transfer^a_{g,n,p}(N,d) := R(\ev^+ \times c)_\ast \left( \cL^a_{\lev} \otimes (\ev^-)^\ast(-) \right).
\]

\begin{introthm}[= \Cref{T:admissibility}]\label{T:main_intro} In characteristic $0$, the
    functor $\transfer^a_{g,n,p}(N,d)$ takes perfect complexes to perfect
    complexes.
\end{introthm}

\begin{defn}
For any level $a>0$, genus $g\geq 0$, and degree $d$, the $K$-theoretic Gromov-Witten gauged invariants for $\GL_N$ are the homomorphisms
\[
    \ggw^a_{g,n,p}(N,d) : K^\ast(\overline{\cM}_{g,n+p}) \otimes \Rep(\GL_N)^{\otimes n} \to \Rep(\GL_N)^{\otimes p}
\]
taking $V \otimes F$ to the derived pushforward of $V \boxtimes \transfer_{g,n,p}^a(N,d)(F)$ along the morphism $\StMap_{g,n+p} \times (\BGL_N)^p \to (\BGL_N)^p$.
\end{defn}

These operations are natural from a moduli theoretic perspective, but it is
not yet clear what algebraic structure governs them. In \Cref{subsection: explicit
rank 1 and genus 0} we describe the binary operations on $\Rep(\mathbb{G}_m) \cong \bZ[x^{\pm 1}]$
corresponding to the case when $N=1$ and $g=0$. This is enough to show that
these operations do \emph{not} satisfy the WDVV equations, so it is not clear
exactly how to interpret these operations from the perspective of oriented topological quantum field theories. However,
we give evidence in favor of such an interpretation at the moduli-theoretic level in \Cref{subsection: gluing and
forgetful morphisms} by constructing gluing morphisms that are compatible with
the level line bundle $\cL_{\lev}$.

\subsection*{Structural results}  We work over a fixed Noetherian base
scheme $S$. Consider a triple $(\varphi : \tilde{C} \to C,
\sigma_\bullet^\pm,E)$ where $C$ and $\tilde{C}$ are flat families of prestable
curves over an $S$-scheme $T$, $E$ is a locally free sheaf on $\tilde{C}$, and
$\sigma_1^-,\ldots,\sigma_n^-,\sigma_1^+,\ldots,\sigma_p^+$ are pairwise
disjoint sections lying in the smooth locus of $\tilde{C} \to T$.

\begin{defn}[Pure Gieseker bundles]
    We say that $(\varphi : \tilde{C} \to C, \sigma_\bullet^\pm,E)$ is a
    $T$-family of pure Gieseker bundles if the following hold:
    \begin{enumerate}
        \item $\varphi : \tilde{C} \to C$ is a log-crepant contraction of curves
        (see \Cref{defn: contraction}).
        \item $\det(E)$ is $\varphi$-ample and the counit
        $\varphi^\ast(\varphi_\ast(E)) \to E$ is surjective.
        \item the $T$-fibers of $\varphi_\ast(E(-\sum_{i=1}^p \sigma_i^+))$ are
        pure of dimension $1$.
    \end{enumerate}
\end{defn}

We let $\GBun_{g,n,p}(N,d)$ denote the stack of pure Gieseker bundles in which
$C$ has genus $g$ fibers, and $E(-\sum_{i=1}^p \sigma_i^+)$ has rank $N$ and
total degree $d$.  The stack $\StMap_{g,n,p}(\BGL_N,d)$ discussed above is the
substack of $\GBun_{g,n,p}(N,d)$ of points whose underlying target curve $C$, along with the markings $\varphi
\circ \sigma_\bullet^{\pm}$, is a family of $(n+p)$-pointed stable curves. The
morphism $\stabi : \StMap_{g,n,p}(\BGL_N,d) \to \StMap_{g,n,p}$ assigns a family of
Gieseker bundles to this underlying target stable curve. Restricting
the universal bundle $E$ along the sections gives evaluation morphisms $\ev^- :
\StMap_{g,n,p}(\BGL_N,d) \to (\BGL_N)^n$ and $\ev^+ : \StMap_{g,n,p}(\BGL_N,d)
\to (\BGL_N)^p$.

Let $(\varphi: \widetilde{C} \to C, \sigma_{\bullet}^{\pm}, E)$ denote the
universal Gieseker bundle over $\StMap_{g,n,p}(\BGL_N,d)$, and let $\pi:
\widetilde{C} \to \StMap_{g,n,p}(\BGL_N,d)$ denote the structure morphism. We
define the level line bundle $\cL_{\lev}$ on $\StMap_{g,n,p}(\BGL_N,d)$ by
$\cL_{\lev} := \det\left(R\pi_*\left( E \oplus E \otimes
\omega_{\pi}\right)\right)^{\vee}$.

In order to show \Cref{T:main_intro}, we construct a $\Theta$-stratification for the moduli stack
$\GBun_{g,n,p}^{\pol}(N,d)$ parameterizing pure Gieseker bundles $(\varphi:
\widetilde{C} \to C, \sigma_{\bullet}^{\pm}, E)$ along with a relatively ample
line bundle $\cH$ on the target family of curves $C$. Using the morphism $c
\times \ev_+$, we may view $\GBun_{g,n,p}^{\pol}(N,d)$ as a stack over the
moduli $\prestable_{g,n,p}^{\pol} \times (\BGL_N)^p$ of polarized $(n,p)$-marked
prestable curves along with a $p$-tuple of vector bundles on the base. Our main
structural result for the stack of polarized pure Gieseker bundles in
characteristic $0$ is as follows (see also \Cref{thm: theta stratification for
Gieseker bundles} for a more general statement).
\begin{introthm} \label{introthm: B} Suppose that the base scheme $S$ is a
$\mathbb{Q}$-scheme. Then,
    \begin{enumerate}
        \item $\nu_\epsilon$ defines a $\Theta$-stratification on
        $\GBun_{g,n,p}^{\pol}(N,d)$ relative to $\prestable_{g,n,p}^{\pol}
        \times (\BGL_N)^p$ in the sense of \cite{halpernleistner2018structure}.
        \item The open semistable locus
        $\GBun_{g,n,p}^{\pol}(N,d)^{\nu_{\epsilon} \dash ss}$ admits a proper
        and flat relative good moduli space over $\prestable_{g,n,p}^{\pol}
        \times (\BGL_N)^p$.
        \item The center of every stratum admits another $\Theta$-stratification
        such that its centers have proper relative good moduli spaces over
        $\prestable_{g,n,p}^{\pol} \times (\BGL_N)^p$.
    \end{enumerate}
\end{introthm}
Base-changing via the morphism $\overline{\cM}_{g,n+p} \times (\BGL_N)^p \to
\prestable_{g,n,p}^{\pol} \times (\BGL_N)^p$ given by the log-canonical
polarization, we obtain a $\Theta$-stratification for $\StMap_{g,n,p}(\GL_N, d)$
relative to $\overline{\cM}_{g,n+p} \times (\BGL_N)^p$, as well as a proper and
flat relative good moduli space of $\nu_{\epsilon}$-semistable pure Gieseker
bundles over $\overline{\cM}_{g,n+p} \times (\BGL_N)^p$ (see \Cref{coroll: theta
stratification for Gieseker bundles over stable curves}).

The key idea in the proof of \Cref{introthm: B} is an interpretation of
$\GBun_{g,n,p}^{\pol}(N)$ in terms of Kontsevich stable maps to a universal Quot
scheme over a moduli of decorated pure sheaves (see \Cref{subsection: properness
over the stack of pure marked sheaves}). We believe that this interpretation is
of independent interest and offers a conceptually new approach to understanding
properness of the moduli of Gieseker bundles. The techniques in this paper have
already been used to construct compactifications of the Hitchin fibration over
stable curves \cite{donagi2024meromorphichitchinfibrationstable}; we expect that
they will have further applications to other moduli problems, even in settings
involving higher dimensional varieties.

\subsection*{Projective targets} In \Cref{section: projective targets}, we
explain how to adapt the ideas in the proof of \Cref{introthm: B} to the more
general setting of gauged maps to a projective $\GL_N$-scheme $Z$. We define
stacks of Gieseker maps $\GMap_{g,n,p}(Z/\GL_N)$ (see \Cref{defn: stack of
gieseker stable maps}), and construct a $\Theta$-stratification and moduli space
for the moduli of Gieseker maps over polarized prestable curves in \Cref{thm:
theta stratification for projective target}.

\subsection*{Concrete description of pure Gieseker bundles} The following is an
informal discussion of the definition of pure Gieseker bundles, which may help
the reader to get intuition on the nature of the moduli problem. Let $k$ be an
algebraically closed field. A pure Gieseker bundle over $k$ consists of the data
of a log crepant contraction $\varphi: (\widetilde{C}, \sigma_{\bullet}^{\pm})
\to (C, \overline{\sigma}_{\bullet}^{\pm})$ of $(n,p)$-marked nodal curves of
genus $g$ over $k$, and a rank $N$ vector bundle $E$ on the source curve
$\widetilde{C}$ satisfying some open conditions.

By \Cref{lemma: control of bubblings}, the contraction $\varphi: \widetilde{C}
    \to C$ is an isomorphism over the open complement $U \subset C$ of the union
    of nodes $\{ \nu_1, \nu_2, \ldots, \nu_m\} \subset C$ and positive markings
    $\{\overline{\sigma}_1^+, \ldots, \overline{\sigma}_p^+\} \subset C$.
    Furthermore, every connected subcurve $R \subset \widetilde{C}$ that is
    contracted to a point by the morphism $\varphi$ is isomorphic to a chain of
    smooth rational curves, as pictured below.
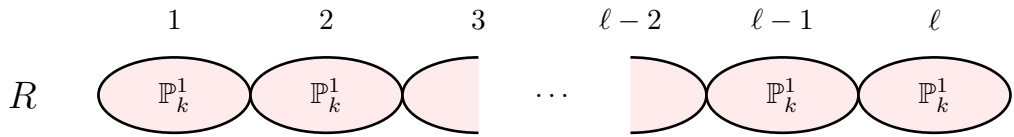
\begin{figure}[!ht]
    \centering
    \begin{circuitikz}

        \node [font=\Large] at (-2,0) {$R$}; \draw [ fill={rgb,255:red,255;
        green,235; blue,235} , line width=1pt ] (0,0) ellipse (1 and 0.5) node
        {$\mathbb{P}^1_k$} ; \node  at (0,1) {$1$}; \draw [
        fill={rgb,255:red,255; green,235; blue,235} , line width=1pt ] (2,0)
        ellipse (1 and 0.5) node {$\mathbb{P}^1_k$} ; \node  at (2,1) {$2$};

        \draw[ fill={rgb,255:red,255; green,235; blue,235} , line width=1pt ]
        (4,0.5) arc(90:270: 1 and 0.5); \node  at (4,1) {$3$};

        \node  at (5,0) {$\dots$};

        \draw[ fill={rgb,255:red,255; green,235; blue,235} , line width=1pt ]
        (6,-0.5) arc(-90:90: 1 and 0.5); \node  at (6,1) {$\ell-2$};

        \draw [ fill={rgb,255:red,255; green,235; blue,235} , line width=1pt ]
        (8,0) ellipse (1 and 0.5) node {$\mathbb{P}^1_k$} ; \node  at (8,1)
        {$\ell-1$};

        \draw [ fill={rgb,255:red,255; green,235; blue,235} , line width=1pt ]
        (10,0) ellipse (1 and 0.5) node {$\mathbb{P}^1_k$} ; \node  at (10,1)
        {$\ell$};
    \end{circuitikz}
\end{figure}

Altogether, this implies that $\widetilde{C}$ is obtained from $C$ by adding
rational bridges separating some of the nodes of $C$, and also attaching
rational tails at some of the positive markings $\overline{\sigma}_i^+$.
Moreover, the log-crepant condition forces the marking $\sigma_i^+$ to lie in
the terminal component of the rational tail. In Figure 1 we provide an example
of such log-crepant contraction of $(1,2)$-marked curves of genus $4$ with one
rational tail of length $2$ containing the first positive marking $\sigma_1^+$
in its terminal component, and one rational bridge of length $1$ separating the
second node $\nu_2$ of $C$.

\begin{figure}[!ht] \label{fig: figure 2}
    \centering
    \begin{circuitikz}

        \draw [ fill={rgb,255:red,255; green,235; blue,235} , line width=1pt]
        (-1.5,-1.5) ellipse (1 and 0.5); \node at (-1.5,-1.5) {$\bullet
        \sigma_1^+$};

        \draw [ fill={rgb,255:red,255; green,235; blue,235} , line width=1pt , ]
        (0,-1.5) ellipse (0.5 and 1);

        \draw [ fill={rgb,255:red,238; green,236; blue,236} , line width=1pt ]
        (0,0) ellipse (1 and 0.5) ;

        \node at (-2,0) [rectangle,draw] {$g_1=1$};

        \draw [ fill={rgb,255:red,238; green,236; blue,236} , line width=1pt ]
        (0,1.5) ellipse (0.5 and 1) ; \node at (0,1.75) {$\bullet \sigma_1^-$};
        \node at (-1.5,1.5) [rectangle,draw] {$g_2=2$};

        \draw [ fill={rgb,255:red,255; green,235; blue,235} , line width=1pt]
        (1.5,1.5) ellipse (1 and 0.5);

        \draw [ fill={rgb,255:red,238; green,236; blue,236} , line width=1pt,
        rotate around={315:(2.25,0)}] (2.25,0.5) ellipse (1 and 0.5) ; \node at
        (2.75,0.25) {$\bullet \sigma_2^+$}; \node at (2.75,-1) [rectangle,draw]
        {$g_3=1$};

        \node [font=\LARGE] at (0,3) {$\widetilde{C}$};

        \node [font=\LARGE] at (4.5,0) {$\xrightarrow{\varphi}$};

        \node [font=\LARGE] at (7,3) {$C$};

        \draw [ fill={rgb,255:red,238; green,236; blue,236} , line width=1pt ]
        (7,1.5) ellipse (0.5 and 1) ; \node at (7,2) {$\bullet
        \overline{\sigma}_1^-$}; \node at (5.5,1.5) [rectangle,draw] {$g_2=2$};

        \node at (7,0.5) {$\bullet$}; \node at (6.5,0.5) {$\nu_1$};

        \draw [ fill={rgb,255:red,238; green,236; blue,236} , line width=1pt ]
        (7,-0.5) ellipse (0.5 and 1) ; \node at (7,-1) {$\bullet
        \overline{\sigma}_1^+$}; \node at (7,-2) [rectangle,draw] {$g_1=1$};

        \node at (7.5,1.35) {$\bullet$}; \node at (7.75,0.85) {$\nu_2$};

        \draw [ fill={rgb,255:red,238; green,236; blue,236} , line width=1pt ]
        (8.5,1.5) ellipse (1 and 0.5) ; \node at (9,1.5) {$\bullet
        \overline{\sigma}_2^+$}; \node at (9,0.5) [rectangle,draw] {$g_3=1$};

    \end{circuitikz}
    \caption{Figure 1}
\end{figure}

Next, let us describe the conditions on the vector bundle $E$. Let $R \subset
\widetilde{C}$ be a maximal connected subcurve contracted by $\varphi$, which in
view of the above is either a rational bridge or a rational tail containing a
positive marking in its terminal component. Given a tuple of integers $(m_1,
m_2, \ldots, m_\ell)$, let us denote by $\cO(m_1 \mid  m_2 \mid \ldots \mid
m_\ell) = \cO(m_{\bullet})$ the unique isomorphism class of line bundle on $R$
whose restriction to the $i^{th}$ rational component of $R$ is
$\cO_{\mathbb{P}^1_k}(m_i)$ as pictured below.
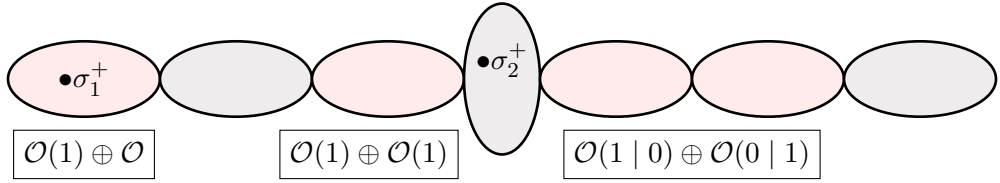
\begin{figure}[!ht]
    \centering
    \begin{circuitikz}

        \node  at (-2.5,-0.5) [font=\Large] {$\cO(m_{\bullet})=$};

        \draw [ fill={rgb,255:red,255; green,235; blue,235} , line width=1pt ]
        (0,0) ellipse (1 and 0.5) node {$\mathbb{P}^1_k$} ; \node  at (0,-1)
        {$\cO_{\mathbb{P}^1_k}(m_1)$}; \draw [ fill={rgb,255:red,255; green,235;
        blue,235} , line width=1pt ] (2,0) ellipse (1 and 0.5) node
        {$\mathbb{P}^1_k$} ; \node  at (2,-1) {$\cO_{\mathbb{P}^1_k}(m_2)$};

        \draw[ fill={rgb,255:red,255; green,235; blue,235} , line width=1pt ]
        (4,0.5) arc(90:270: 1 and 0.5);

        \node  at (5,0) {$\dots$};

        \draw[ fill={rgb,255:red,255; green,235; blue,235} , line width=1pt ]
        (6,-0.5) arc(-90:90: 1 and 0.5);

        \draw [ fill={rgb,255:red,255; green,235; blue,235} , line width=1pt ]
        (8,0) ellipse (1 and 0.5) node {$\mathbb{P}^1_k$} ; \node  at (8,-1)
        {$\cO_{\mathbb{P}^1_k}(m_{\ell-1})$};

        \draw [ fill={rgb,255:red,255; green,235; blue,235} , line width=1pt ]
        (10,0) ellipse (1 and 0.5) node {$\mathbb{P}^1_k$} ; \node  at (10,-1)
        {$\cO_{\mathbb{P}^1_k}(m_{\ell})$};
    \end{circuitikz}
\end{figure}

By \cite[Sect. 4]{martens_thaddeus_variations_grothendieck}, the restriction
$E|_R$ is isomorphic to a direct sum of $N$ line bundles $E|_R \cong
\bigoplus_{j=1}^N \cO(m_{1,j} \mid m_{2,j} \mid \ldots \mid m_{\ell, j})$. The
condition that the counit $\varphi^*\varphi_*(E) \to E$ is surjective is
equivalent to imposing that for all such maximal rational chains $R$ contracted
by $\varphi$ the numbers $m_{j,i}$ are nonnegative (in the language of
\cite[Defn.-Notn. 1]{nagaraj-seshadri-2}, we want each restriction $E|_{R}$ to
be positive). The condition that $\det(E)$ is $\varphi$-ample amounts to
imposing in addition that for all such $R$ and all $1 \leq i \leq \ell$, there
exists some $1 \leq j \leq N$ such that $m_{i,j}>0$ (we require that $E|_R$ is
strictly positive in the sense of \cite[Defn.-Notn. 1]{nagaraj-seshadri-2}).
Finally, the additional condition that the pushforward $\varphi_*(E(-
\sum_{i=1}^p \sigma_i^+))$ is pure holds if and only if for all $R$ and for all
$1 \leq j \leq N$, we have either $m_{\bullet, j}=0$ or $m_{\bullet, j} =
\delta_{i,j}$ for some $1 \leq i \leq \ell$, where $\delta$ denotes the
Kronecker delta function. In other words, each line bundle direct summand of
$E|_R$ is either isomorphic to the trivial bundle $\cO_R$, or is of the form
$\cO(0 \mid \ldots  \mid 0 \mid 1 \mid 0 \mid \ldots \mid 0)$ where all the
$m_{\bullet, j}$ are $0$ except for $m_{i,j}=1$. We note that all three
conditions on $E$ jointly imply that the length $\ell$ of each rational chain
$R$ contracted by $\varphi$ is bounded by the rank $N$. In Figure 2  below, we
depict a pure Gieseker vector bundle of rank $2$ on a $(0,2)$-marked curve, where
the contracting rational irreducible components are pictured in red and we omit
the target curve.

\begin{figure}[!ht]
    \centering
    \begin{circuitikz}

        \draw [ fill={rgb,255:red,255; green,235; blue,235} , line width=1pt]
        (-1.5,0) ellipse (1 and 0.5); \node at (-1.5,0) {$\bullet \sigma_1^+$};

        \node at (-1.75,-1) [rectangle,draw] {$\cO(1) \oplus \cO(1)$};

        \draw [ fill={rgb,255:red,238; green,236; blue,236} , line width=1pt ]
        (0,0) ellipse (0.5 and 1) ; \node at (0,0.25) {$\bullet \sigma_2^+$};

        \draw [ fill={rgb,255:red,255; green,235; blue,235} , line width=1pt , ]
        (1.5,0) ellipse (1 and 0.5);

        \node at (2.5,-1) [rectangle,draw] {$\cO(1\mid 0) \oplus \cO(0\mid 1)$};

        \draw [ fill={rgb,255:red,255; green,235; blue,235} , line width=1pt , ]
        (3.5,0) ellipse (1 and 0.5);

        \draw [ fill={rgb,255:red,238; green,236; blue,236} , line width=1pt]
        (5.5,0) ellipse (1 and 0.5) ;
    \end{circuitikz}
    \caption{Figure 2}
\end{figure}

One important idea in these pictures is that we allow the formation of rational
tails over the positive markings $\overline{\sigma}_i^+ \in C$ as a means of
``compactifying'' the ``output'' evaluation morphism $\ev_+$. The operation of
gluing a positive ``output'' marking to a negative ``input'' marking respects
the property of being a pure Gieseker bundle, and hence it induces well-defined
gluing morphisms at the level of moduli stacks. As an example, we may glue the
positive marking $\sigma_1^+$ in the pure Gieseker vector bundle pictured in
Figure 2 with the negative marking $\sigma_1^-$ in the following pure Gieseker
bundle on a $(1,1)$-marked curve:

\begin{figure}[!ht]
    \centering
    \begin{circuitikz}

        \draw [ fill={rgb,255:red,255; green,235; blue,235} , line width=1pt]
        (-1.5,0) ellipse (1 and 0.5); \node at (-1.5,0) {$\bullet \sigma_1^+$};

        \node at (-1.75,-1) [rectangle,draw] {$\cO(1) \oplus \cO$};

        \draw [ fill={rgb,255:red,238; green,236; blue,236} , line width=1pt ]
        (0.5,0) ellipse (1 and 0.5) ; \node at (0.5,0) {$\bullet \sigma_1^-$};
    \end{circuitikz}
\end{figure}

in order to obtain the following pure Gieseker bundle on a $(0,2)$-marked curve:

\begin{figure}[!ht]
    \centering
    \begin{circuitikz}

        \draw [ fill={rgb,255:red,255; green,235; blue,235} , line width=1pt]
        (-5.5,0) ellipse (1 and 0.5); \node at (-5.5,0) {$\bullet \sigma_1^+$};

        \node at (-5.5,-1) [rectangle,draw] {$\cO(1) \oplus \cO$};

        \draw [ fill={rgb,255:red,238; green,236; blue,236} , line width=1pt ]
        (-3.5,0) ellipse (1 and 0.5) ;

        \draw [ fill={rgb,255:red,255; green,235; blue,235} , line width=1pt]
        (-1.5,0) ellipse (1 and 0.5);

        \node at (-1.75,-1) [rectangle,draw] {$\cO(1) \oplus \cO(1)$};

        \draw [ fill={rgb,255:red,238; green,236; blue,236} , line width=1pt ]
        (0,0) ellipse (0.5 and 1) ; \node at (0,0.25) {$\bullet \sigma_2^+$};

        \draw [ fill={rgb,255:red,255; green,235; blue,235} , line width=1pt , ]
        (1.5,0) ellipse (1 and 0.5);

        \node at (2.5,-1) [rectangle,draw] {$\cO(1\mid 0) \oplus \cO(0\mid 1)$};

        \draw [ fill={rgb,255:red,255; green,235; blue,235} , line width=1pt , ]
        (3.5,0) ellipse (1 and 0.5);

        \draw [ fill={rgb,255:red,238; green,236; blue,236} , line width=1pt]
        (5.5,0) ellipse (1 and 0.5) ;
    \end{circuitikz}
\end{figure}

The operation of gluing of two positive markings does not respect the pure
Gieseker bundle property; this is an important reason for the different roles
played by negative and positive in our moduli problem.

\subsection*{Comparison with previous work in the literature} \;

\subsubsection*{The case without positive marked points.} In the
case when $n=p=0$, the moduli of semistable pure Gieseker bundles in
\Cref{coroll: theta stratification for Gieseker bundles over stable curves}
recovers the moduli space whose study was initiated by Gieseker
\cite{gieseker-degeneration}, refined by Nagaraj-Seshadri
\cite{nagaraj-seshadri-2} and constructed in full generality by Schmitt
\cite{schmitt-hilbert-compactification} using GIT. In this special case, the
main novelty of this work is the development of the $\Theta$-stratification of the unstable locus.

As explained in Pandharipande's work \cite{pandharipande-compactification}, the
moduli of Gieseker bundles with $n=p=0$ agrees with Caporaso's compactified
Jacobian \cite{caporaso_compactification} when $N=1$, but it differs from the
moduli of torsion-free sheaves constructed in
\cite{pandharipande-compactification} when $N>1$. As a consequence, the
semistable locus of the moduli problem constructed in \Cref{section: moduli of
gieseker bundles} when $N=1$ can be related to the moduli space of torsion-free
sheaves on stable maps considered by Cooper \cite{cooper_compactified_jac},
which can be thought of as an avatar of the stack of stable maps to a quotient
stack $X/\GL_N$ when the action of $\GL_N$ on $X$ is trivial. In comparison, our
moduli problem has the advantages of (1) allowing the case when the action of
$\GL_N$ is nontrivial, (2) admitting a natural quasi-smooth enhancement when $X$
is smooth for all $N \geq 1$, and (3) admitting well-defined gluing and
evaluation morphisms.

A stratification of the stack $\GBun_{g,n,0}(1)$ (with $p=0$ and $N=1$) was
considered by Frenkel-Teleman-Tolland in \cite{frenkel-teleman-tolland-ggw}.
They defined $K$-theoretic Gromov-Witten invariants for $B\mathbb{G}_m$ and
suggested that a generalization should hold for targets of the form
$X/\mathbb{G}_m$ with $X$ projective. We note that the level line bundle
$\cL_{\lev}$ we define differs from the one used in
\cite{frenkel-teleman-tolland-ggw}; this is necessary to ensure compatibility
with the gluing morphisms defined in \Cref{subsection: gluing and forgetful
morphisms}.

For $N>1$, Kausz took steps towards developing Gromov-Witten invariants for $B\GL_N$ in \cite{kausz_stable_maps}. His stack of Gieseker bundles with markings agrees
with our $\GBun_{g,n,0}(N)$, without positive markings, but his gluing morphisms differ from ours: he compactifies the space of identifications
between the fibers of the vector bundles, as opposed to using genuine
isomorphisms of the fibers as in \Cref{subsection: gluing and forgetful
morphisms}. Kausz also identified the lack of a
Harder-Narasimhan stratification as one of the main obstacles in defining
$K$-theoretic Gromov-Witten invariants in this setting, so our results fully realize Kausz's
proposal.

Finally, when $g>0$ we expect our $K$-theoretic Gromov-Witten invariants to be
substantially different from the gauged Gromov-Witten invariants defined in the
work of Woodward, Gonzalez-Woodward and Gonzalez-Solis-Woodward
\cite{woodward-quantum-quotients, gonzalez-solis-woodward-properness,
gonzalez-solis-woodward-stable-gauged, gonzalez-woodward-quantum}, and also
considered by the authors in \cite{gauged_theta_stratifications}. The reason is
that in the Gonzalez-Solis-Woodward setting the moduli problem is constructed
relative to a fixed smooth curve of genus $g$, whereas one should think of our
moduli problem as a more direct analog of stable maps, where we allow the curve
to genuinely vary.

\subsubsection*{The case with positive marked points.} One of the
main novelties in this work is the introduction of positive marked points to
compactify the evaluation morphisms. To our knowledge, this version of the
moduli problem has not been considered in the literature, even in the case when the rank
$N=1$.

The closest antecedent, which helped inspire our construction, are the stable quasi-maps from a fixed smooth curve $C$ to a quotient stack \emph{relative} to a marked point $p \in C$, studied by Okounkov in \cite[Sect.~6.3.1]{okounkovlectures}. For that construction, $G$ is a reductive group and $V$ is an affine $G$-variety with a GIT-semistable locus $V^{\rm{ss}} \subset V$. A stable quasi-map relative to $p$ is a morphism of marked prestable curves $(C',p') \to (C,p)$ that contracts a chain of rational curves over $p$ and a morphism $f : C' \to V/G$ that sends $p'$, every node, and every generic point of $C'$ to $V^{\rm{ss}}/G$, and has a finite automorphism group. This defines a separated DM stack $\operatorname{QM}_{{\rm{relative}}\,p}$ such that the evaluation morphism $\operatorname{QM}_{{\rm{relative}}\,p} \to V^{\rm{ss}}/G$ is proper. The goal is similar to ours, but for the construction of $\operatorname{QM}_{{\rm{relative}}\,p}$ to work, it is crucial that $G$ acts with finite stabilizers on $V^{\rm{ss}}$, so in a technical sense it is the opposite extreme from our situation.

\subsection*{Notation}
In this paper, we work over a fixed Noetherian base scheme $S$. Given any two
stacks $\cX, \cY$ over $S$, the undecorated product $\cX \times \cY$ should be
understood as the fiber product over $S$, unless stated otherwise. Throughout
the paper, we fix integers $g,n,p \geq 0$ (the genus, number of positive
markings, and number of negative markings).

For convenience, we list here the notation for some stacks considered in this
paper:
\begin{itemize}
    \item $\prestable_{g,n,p}$: The stack of marked prestable curves with
    positive and negative markings (see \Cref{D:prestable}).
    \item $\prestablepol_{g,n,p}$: The stack of polarized marked prestable
    curves (see \Cref{defn: stack of polarized curves}).
    \item Several variants of the stack of marked coherent sheaves on prestable
    curves (see \Cref{defn: stack of marked sheaves}, \Cref{defn: stack of pure
    marked sheaves}):
          \begin{itemize}
              \item $\Coh_{g,n,p}$: the stack of all marked coherent sheaves on
              prestable curves
              \item $\Coh_{g,n,p}(N)$: marked coherent sheaves whose restriction
              to every generic point has rank $N$ and that are marked by rank
              $N$ locally free sheaves
              \item $\Pur_{g,n,p} \subset \Coh_{g,n,p}$: the substack of pure
              marked sheaves
              \item $\Pur_{g,n,p}(N) := \Coh_{g,n,p}(N)\cap \Pur_{g,n,p}$.
          \end{itemize}
    \item $\CBun_{g,n,p}(N,d)$: The stack of Gieseker bundles of rank $N$ and
    degree $d$ (see \Cref{D:GBun}, \Cref{notn: degree of a Gieseker bundle}).
    \item $\GBun_{g,n,p}(N,d)$: The stack of pure Gieseker bundles of rank $N$
    and degree $d$ (see \Cref{defn: stack of pure gieseker bundles}, \Cref{notn:
    degree of a Gieseker bundle}).
    \item $\StMap_{g,n,p}(\BGL_N, d)$: The substack of $\GBun_{g,n,p}(N,d)$
    where the underlying target curve is an $(n+p)$-pointed stable curve (see
    \Cref{notn: gieseker bundles over stable curves}).
    \item $\GMap_{g,n,p}(Z/\GL_N, d)$ the stack of Gieseker maps to a quotient
    stack $Z/\GL_N$ with degree $d \in H_2(Z/\GL_N)$ (see \Cref{defn: stack of
    gieseker stable maps}, \Cref{defn: degree of Gieseker maps}).
\end{itemize}

We also consider variants of the stacks above in which the underlying nodal
family of curves is equipped with a relatively ample line bundle. For example,
$\Coh^{\rm{\pol}}_{g,n,p} := \Coh_{g,n,p} \times_{\prestable_{g,n,p}}
\prestablepol_{g,n,p}$, and the superscript $\rm{\pol}$ will always denote a
fiber product with $\prestablepol_{g,n,p}$ over $\prestable_{g,n,p}$. The
presence of the polarization on a prestable curve allows us to define a notion
of (global) rank of a coherent sheaf on that curve (\Cref{notn: Hilbert
polynomial and rank of pure sheaves}), and we let $\Pur^{\pol,r,d}_{g,n,p}
\subset \Pur^{\pol}_{g,n,p}$ denote the substack of pure marked sheaves of
degree $d$ and rank $r$ with respect to the polarization on the curve.

Let us also record some of the relevant morphisms in this paper:
\begin{itemize}
    \item $\Psi: \Pur_{g,n,p} \to \Pur_{g,n+p,0}$ (see \Cref{defn: psi})
    \item $c: \CBun_{g,n,p}(N) \to \prestable_{g,n,p}$ (see \Cref{L:cBun
    algebraic}).
    \item $\Phi : \CBun_{g,n,p}(N) \to \Coh_{g,n,p}(N)$ (see \Cref{lemma: Phis
    is well-defined}).
    \item $\ev_-: \GBun_{g,n,p}(N) \to (\BGL_N)^n$ and $\ev_+: \GBun_{g,n,p}(N)
    \to (\BGL_N)^p$ (see \Cref{defn: evaluation morphismsm}).
    \item $\Phi_Z : \GMap_{g,n,p}(Z/\GL_N) \to \Pur_{g,n,p}(N)$ (see
    \Cref{lemma: PhiZ is well-defined}).
\end{itemize}

\medskip

\noindent \textbf{Intrinsic GIT.} In this paper we use the notions of
$\Theta$-stratifications and numerical invariants for stacks in the sense of
\cite{halpernleistner2018structure}. We appeal to the intrinsic GIT criteria in
\cite[Thm. B]{halpernleistner2018structure} to construct
$\Theta$-stratifications and moduli spaces. We refer the reader to \cite[Sect.
2.4]{gauged_theta_stratifications} for an overview of the necessary background,
and to \cite{halpernleistner2018structure} for a more detailed discussion.

\subsection*{Acknowledgements}
We thank Jarod Alper, Younghan Bae, Ron Donagi, Eduardo Gonzalez, Giovanni Inchiostro,
Davesh Maulik, Andrei Okounkov, Nick Rozenblyum, and Constantin Teleman for discussions related to
the contents of this article. The first author was supported by NSF grants DMS-1945478(CAREER) and FRG-DMS-2052936(CAREER), by the Sloan Foundation research fellowship FG-2022-18834. The
second author was partially supported by a Junior Trimester Program fellowship
at the HIM (Bonn) and an AMS-Simons travel grant during the preparation of this
paper.

\addtocontents{toc}{\protect\setcounter{tocdepth}{2}}

\section{Marked sheaves on prestable curves} \label{section: stack of marked
sheaves} In this section we study the moduli of pure marked sheaves, which
parameterizes pure sheaves on marked prestable curves equipped with certain
extra decorations (see \Cref{defn: stack of marked sheaves}). We define a $\Theta$-stratification for this moduli problem and
construct proper good moduli spaces in the characteristic $0$ setting (see
\Cref{thm: weak theta stratification pure marked sheaves}). A lot of the results
in this section are rather technical in nature; the reader may choose to defer a
careful reading of this section and refer to it as needed.

\subsection{Prestable curves and polarizations}

\begin{defn}[Stack of marked prestable curves] \label{D:prestable} We denote by
    $\prestable_{g,n,p}$ the stack that sends an $S$-scheme $T$ to the groupoid
    of flat, finitely presented, proper families of connected nodal curves $\pi
    : C \to T$ of genus $g$ along with $n$ ``negative" sections
    $\sigma^-_1,\ldots, \sigma^-_n : T \to C$ and $p$ ``positive" sections
    $\sigma^+_1,\ldots,\sigma^+_p : T \to C$ such that all sections lie in the
    smooth locus of $\pi$ and no two distinct sections meet. We call these
    $T$-families of $(n,p)$-marked genus $g$ prestable curves.
\end{defn}

We have a smooth representable morphism $\prestable_{g,n,p} \to
\prestable_{g,0,0}$ that forgets the markings. By
\cite[\href{https://stacks.math.columbia.edu/tag/0E6W}{Tag
0E6W}]{stacks-project}, $\prestable_{g,0,0}$ is a smooth algebraic $S$-stack,
and hence the same holds for $\prestable_{g,n,p}$.

\begin{defn}[Stack of polarized marked curves] \label{defn: stack of polarized
    curves} We denote by $\prestablepol_{g,n,p}$ the stack that sends an
    $S$-scheme $T$ to the groupoid of tuples $(C, \sigma_{\bullet}^{\pm}, \cL)$,
    where $(C, \sigma^{\pm}_{\bullet})$ is a $T$-family $(n,p)$-marked prestable
    curves of genus $g$, and $\cL$ is a $T$-ample line bundle on $C$.
\end{defn}

\begin{prop} \label{prop: properties of the stack of polarized curves} The stack
    $\prestablepol_{g,n,p}$ of polarized $(n,p)$-marked prestable curves of
    genus $g$ is an algebraic stack with affine relative diagonal and locally of
    finite type over $S$. The morphism $\prestablepol_{g,n,p} \to
    \prestable_{g,n,p}$ that forgets the polarization is smooth and surjective.
\end{prop}
\begin{proof}
    The stack $\prestablepol_{g,n,p}$ is smooth and separated over
    $\prestablepol_{g,0,0}$ via the morphism that forgets the markings. On the
    other hand, $\prestablepol_{g,0,0}$ is an open substack (
    \cite[\href{https://stacks.math.columbia.edu/tag/0D5B}{Tag
    0D5B}]{stacks-project} and
    \cite[\href{https://stacks.math.columbia.edu/tag/0D5B}{Tag
    0E6U}]{stacks-project}) of the stack of polarized schemes, which is an
    algebraic stack of locally of finite type
    \cite[\href{https://stacks.math.columbia.edu/tag/0DPS}{Tag
    0DPS}]{stacks-project} and with affine relative diagonal \cite[Sect.
    2.1]{starr-dejong-amost-proper} over $S$. The morphism
    $\prestablepol_{g,n,p} \to \prestable_{g,n,p}$ that forgets the polarization
    is smooth and surjective by
    \cite[\href{https://stacks.math.columbia.edu/tag/0DQ0}{Tag
    0DQ0}]{stacks-project}.
\end{proof}

\subsection{Marked sheaves on prestable curves}

\begin{defn}[Stack of marked sheaves] \label{defn: stack of marked sheaves} Let
    $\Coh_{g,n,p}$ denote the stack whose $T$-points consist of an
    $(n,p)$-marked family of prestable curves $(C,\sigma_\bullet^{\pm})$ over
    $T$ along with
    \begin{enumerate}
        \item a $T$-flat locally finitely presented quasicoherent sheaf $F$ on
        $C$, and
        \item for each $1 \leq i \leq p$, a locally free quotient
        $(\sigma^+_i)^\ast(F) \twoheadrightarrow W_i$.
    \end{enumerate}
    We let $\Coh_{g,n,p}(N) \subset \Coh_{g,n,p}$ denote the substack of
    families such that each $W_i$ has rank $N$ and $F$ has rank $N$ at every
    generic point of every $T$-fiber.

\end{defn}

Consider a $T$-point of $\Coh_{g,n,p}$ as in \Cref{defn: stack of marked
sheaves}. For each $i$, adjunction induces a surjection $F \twoheadrightarrow
(\sigma^+_i)_\ast (W_i)$. Putting them together, we obtain a surjection $F
\twoheadrightarrow \bigoplus_{i=1}^p (\sigma^+_i)_\ast (W_i)$ of $T$-flat
sheaves.

\begin{defn}[Stack of pure marked sheaves] \label{defn: stack of pure marked
    sheaves} Let $\Pur_{g,n,p} \subset \Coh_{g,n,p}$ be the substack of points
    for which the $T$-flat kernel of the surjection $F \twoheadrightarrow
    \bigoplus_{i=1}^p (\sigma^+_i)_\ast (W_i)$ is pure of dimension $1$ on each
    $T$-fiber. In this case, we say that the marked sheaf $(C,
    \sigma_{\bullet}^{\pm}, F, (\sigma_{\bullet}^+)^*(F) \to W_{\bullet})$ is
    pure. We let $\Pur_{g,n,p}(N) \subset \Pur_{g,n,p}$ denote the intersection
    with the substack $\Coh_{g,n,p}(N) \subset \Coh_{g,n,p}$.
\end{defn}

\begin{lemma} \label{lemma: stack of pure sheaves is algebraic} The following
hold:
    \begin{enumerate}
        \item $\Pur_{g,n,p} \subset \Coh_{g,n,p}$ is an open substack.
        \item $\Coh_{g,n,p}(N) \subset \Coh_{g,n,p}$ is an open and closed
        substack.
        \item $\Coh_{g,n,p}$ is an algebraic stack, locally of type over the
        base $S$, and the forgetful morphism $\Coh_{g,n,p} \to
        \prestable_{g,n,p}$ has affine diagonal.
    \end{enumerate}
\end{lemma}
\begin{proof}
    Part (1) follows from \cite[Thm. 12.2.1 (iii)]{egaiv}.

    The stack $\Coh_{\cC/\prestable_{g,n,p}}$ of coherent sheaves on the fibers
    of the universal $(n,p)$-marked curve $\cC \to \prestable_{g,n,p}$ is
    locally of finite type and has relative affine diagonal over
    $\prestable_{g,n,p}$
    \cite[\href{https://stacks.math.columbia.edu/tag/08KA}{Tag
    08KA}]{stacks-project}. To show Part (2), it suffices to prove that the
    locus of geometric points of $\Coh_{\cC/\prestable_{g,n,p}}$ where the
    coherent sheaf has rank $N$ at every generic point of every fiber is open
    and closed. The proof of openness is standard \endnote{Let $T$ be a
    Noetherian scheme, and choose $T \to \Coh_{\cC/\prestable_{g,n,p}}$
    corresponding to a family $(C, \sigma_{\bullet}^{\pm}, F)$. Let $U \subset
    C$ denote the largest open subscheme where $F$ is locally free of rank $N$,
    and let $Z \subset C$ denote the closed complement $C \setminus U$ equipped
    with the reduced subscheme structure. By upper semicontinuity of fiber
    dimension, the locus of points of $Z$ where the $T$-fibers have dimension at
    least $1$ is a closed subset $Z^{\geq 1}$ inside $Z$. Since the projection
    $\pi: Z \to T$ is proper, the image of $\pi(Z^{\geq 1}) \subset T$ is
    closed. The locus of points in $T$ where $F$ has rank $N$ generically is the
    open complement of $\pi(Z^{\geq 1})$, as desired.}; we explain how to show
    this open subset is closed under specialization. Let $R$ be a discrete
    valuation ring and choose $\Spec(R) \to \Coh_{\cC/\prestable_{g,n,p}}$
    corresponding to a family $(C, \sigma_{\bullet}^{\pm}, F)$. We denote by
    $\eta$ (resp. $s$) the generic (resp. special) point of $\Spec(R)$. Suppose
    that the restriction of $F$ to the generic fiber $C_{\eta}$ has rank $N$ at
    every generic point. By the fiberwise criterion for flatness
    \cite[\href{https://stacks.math.columbia.edu/tag/039C}{Tag
    039C}]{stacks-project}, the maximal open subset $U \subset C$ where $F$ is
    locally free contains all the generic points of the special fiber $C_s$.
    Since $F|_{U}$ has locally constant rank and every generic point of $C_{s}$
    generalizes to a generic point of $C_{\eta}$, we conclude that $F|_{C_s}$
    has rank of $N$ at all generic points of $C_s$.

    Part (3) follows from the representability of Hom functors between sheaves
    by affine finite type schemes
    \cite[\href{https://stacks.math.columbia.edu/tag/08K6}{Tag
    08K6}]{stacks-project}.\endnote{Set $\SVect := \bigsqcup_{i \geq 0} \BGL_i$.
    Consider the product $\cS := \Coh_{\cC/\prestable_{g,n,p}} \times \SVect^p$,
    which is a stack locally of finite type and with relative affine diagonal
    over $\prestable_{g,n,p}$. By definition, there is a universal family $(\cC,
    \sigma_{\bullet}^{\pm}) \to \cS$ of $(n,p)$-marked prestable curves, an
    $\cS$-flat universal sheaf $F_{univ}$ on $\cC$, and a universal $p$-tuple
    $(W_i)_{i=1}^p$ of locally free sheaves on $\cS$. By
    \cite[\href{https://stacks.math.columbia.edu/tag/08K6}{Tag
    08K6}]{stacks-project}, for each $1 \leq i \leq p$ there is an affine
    finite type morphism $\Hom((\sigma_i^+)^*(F_{univ}), W_i) \to \cS$
    parametrizing homomorphisms from $(\sigma_i^+)^*(F_{univ})$ to $W_i$. The
    fiber product $\cW:= \prod_{i=1}^p \Hom((\sigma_i^+)^*(F_{univ}), W_i)$ over
    $\cS$ is locally of finite type and with relative affine diagonal over
    $\prestable_{g,n,p}$. Note that $\cW$ parameterizes tuples $(C,
    \sigma_{\bullet}^{\pm}, F, (\sigma_{\bullet}^{+})^*(F) \to W_{\bullet})$,
    and $\cS$ is the open substack of $\cW$ where each of the morphisms
    $(\sigma_{i}^{+})^*(F) \to W_{i}$ is surjective.}
\end{proof}

\subsubsection{An alternative description of pure marked sheaves}

\begin{defn}[The morphism $\Psi$] \label{defn: psi} Let $(\cC,
    \sigma_{\bullet}^{\pm}, \cF, (\sigma_{\bullet}^+)^*\cF \to \cW_{\bullet})
    \to \Pur_{g,n,p}$ be the universal family of pure marked sheaves. We denote
    by $\Psi: \Pur_{g,n,p} \to \Pur_{g,n+p,0}$ the morphism classifying
    $(\cC,\sigma^\pm_\bullet, \ker(\cF \to
    \bigoplus_{i=1}^p(\sigma_i^+)_*\cW_i)$.
\end{defn}

Let $(\cC, \sigma_{\bullet}, \cK)$ denote the universal family of sheaves on
$\Pur_{g,n+p,0}$. For for $1 \leq i \leq p$, set $\cK_i := (\sigma_{n+i})^\ast
(\cK(\sigma_{n+i}))$, which is a locally free sheaf on $\Pur_{g,n+p,0}$. For all
$N \geq 0$ let $\cO^{\oplus N}$ denote the trivial rank $N$ vector bundle on
$\Pur_{g,n+p,0}$. By \cite[\href{https://stacks.math.columbia.edu/tag/08K6}{Tag
08K6}]{stacks-project} there is a finite type relatively affine morphism
$\Hom(\cO^{\oplus N}, \cK_i) \to \Pur_{g,n+p,0}$ parameterizing homomorphisms
$\cO^{\oplus N} \to \cK_i$. We consider the quotient stack $\Hom(\cO^{\oplus N},
\cK_i)/\GL_N$, where $\GL_N = \Aut(\cO^{\oplus N})$ acts on the source
$\cO^{\oplus N}$, and we set
\[\Hom(\cO^{ *}, \cK_i)/\GL_* := \bigsqcup_{N \geq 0} \Hom(\cO^{\oplus N},
\cK_i)/\GL_N.\]

\begin{prop}\label{proposition: marked sheaves affine over nonmarked sheaves}
    Regarding $\Pur_{g,n,p}$ as a stack over $\Pur_{g,n+p,0}$ via the morphism
    $\Psi$, Grothendieck duality along the sections $\sigma_i^+$ gives an
    isomorphism of $\Pur_{g,n+p,0}$-stacks
    \[
        \Pur_{g,n,p} \cong (\Hom(\cO^\ast,\cK_1) / \GL_\ast) \underset{\Pur_{g,n+p,0}}{\times} \cdots \underset{\Pur_{g,n+p,0}}{\times} (\Hom(\cO^\ast,\cK_p)/\GL_\ast).
    \]
\end{prop}

\begin{proof}
    Setting $p=1$ for simplicity and letting $\sigma$ denote the positive marked
    point, $\Pur_{g,n,p}$ parameterizes families of extensions $0 \to K \to F
    \to \sigma_\ast(E) \to 0$, where $E$ is a vector space and $K$ is a pure
    sheaf of dimension $1$. Grothendieck duality for the morphism $\sigma$ gives
    canonical isomorphisms
    \begin{gather*}
        \Ext^1\left(\sigma_*(E), K \right) = \Ext^1\left(E, \sigma^! K \right) =  \Ext^1\left(E, \sigma^\ast(K(\sigma))[-1] \right) = \Hom(E,\sigma^\ast(K(\sigma))).
    \end{gather*}
    The statement of the proposition is just the natural extension of these
        isomorphisms for families of sheaves and multiple positive marked
        points.\endnote{ Let $\cX$ denote the right-hand side of the desired
        isomorphism. By definition $\cX$ sends a scheme $T$ to the groupoid of
        tuples $(C, \sigma_{\bullet}, K, W_{\bullet}, \phi_{\bullet})$, where
        $(C, \sigma_{\bullet})$ is an $(n+p)$-marked curve, $K$ is $T$-family of
        pure sheaves on $C$, $W_{\bullet}$ is a $p$-tuple of locally free
        sheaves on $T$, and $\phi_{\bullet}$ is a $p$-tuple of morphisms
        $\phi_{i}:  W_i \to \sigma_{n+i}^*(K(D^+))$, where we set $D^+ :=
        \sum_{i=n+1}^{n+p} \sigma_i$. For any such $T$-point, let $(C,
        \sigma_{\bullet}^{\pm})$ be the $(n,p)$-marked curve where $\sigma_i^-
        := \sigma_i$ and $\sigma_i^+ := \sigma_{n+i}$. Furthermore, let $F$ be
        the finitely presented $T$-flat sheaf on $C$ given as the pullback,
        \[
            \begin{tikzcd}[ampersand replacement = \&]
                F   \ar[r] \ar[d] \& K(D^+) \ar[d]  \\
                \bigoplus_{i=1}^p (\sigma_i^+)_\ast(W_i)\ar[r] \& K(D^+) \otimes \cO_{D^+}
            \end{tikzcd}
        \]
        where each homomorphism $(\sigma_i^+)_\ast(W_i) = (\sigma_{n+i})_*(W_i)
        \to K(D^+) \otimes \cO_{D^+}$ is induced by $\phi_i$. The sheaf $F$ fits
        into a short exact sequence
        \begin{equation} \label{equation: extensions}
            0 \to K \to F \to \bigoplus_{i=1}^p (\sigma_i^+)_\ast(W_i) \to 0,
        \end{equation}
        which gives $F$ the structure of an $(n,p)$-marked family of pure
        sheaves on $(C, \sigma_{\bullet}^{\pm})$. This defines a morphism $f:\cX
        \to \Pur_{g,n,p}$ of stacks over $\Pur_{g,n+p,0}$. To conclude, we shall
        show that $f$ is an isomorphism. This amounts to proving that for all
        Noetherian schemes $T$ and $T$-points $(C, \sigma_{\bullet}, K,
        W_{\bullet})$ of $\Pur_{g,n+p,0} \times \left(\bigsqcup_{N\geq 0}
        \BGL_N\right)^p$, the assignment defined above induces a bijection
        between the set $\prod_{i=1}^p \Hom\left(W_i,
        (\sigma_{n+i})^*(K(\sigma_{n+i})\right)$ and the set
        $\Ext^1\left(\bigoplus_{i=1}^p (\sigma_{n+i})_*(W_i), K \right)$ of
        extensions as in \eqref{equation: extensions} up to isomorphism. Let
        $\sigma: D^+ \hookrightarrow C$ denote the embedding of the Cartier
        divisor $D^+$. We can view the tuple $(W_i)_{i=1}^p$ as a locally free
        sheaf $E^+$ on $D^+$, so that $\sigma_*(E^+) \cong \bigoplus_{i=1}^p
        (\sigma_{n+i})_*(W_i)$. This way we get natural identifications
        $\prod_{i=1}^p \Hom\left(W_i, (\sigma_{n+i})^*(K(\sigma_{n+i})\right) =
        \Hom(E^+,\sigma^\ast(K(D^+)))$ and $\Ext^1\left(\bigoplus_{i=1}^p
        (\sigma_{n+i})_*(W_i), K \right) = \Ext^1\left(\sigma_*(E^+), K
        \right)$. Using Grothendieck duality for the embedding of the Cartier
        divisor $D^+ \subset C$, we have a chain of canonical identifications
        \begin{gather*}
            \Ext^1\left(\sigma_*(E^+), K \right) = \Ext^1\left(E^+, \sigma^! K \right) =  \Ext^1\left(E^+, \sigma^\ast(K(D^+))[-1] \right) = \Hom(E^+,\sigma^\ast(K(D^+))),
        \end{gather*}
        which can be checked to recover the assignment we defined above.}
\end{proof}

If we denote by $\SVect := \bigsqcup_{N \geq 0} \BGL_N$ the classifying stack of
locally free sheaves, then we have a morphism $\ev^+: \Pur^{\pol}_{g,n,p} \to
\SVect^p$ keeping track of the $p$-tuple of locally free sheaves.
\begin{prop} \label{prop: affineness of pur and valuative criterion for
    properness} The following hold:
    \begin{enumerate}
        \item The morphism $\Psi\times \ev^+: \Pur_{g,n,p} \to
        \Pur_{g,n+p,0}\times \SVect^p$ is affine and of finite type.
        \item The morphism $\Pur_{g,n,p} \to \prestable_{g,n,p} \times \SVect^p$
        satisfies the existence part of the valuative criterion for properness.
    \end{enumerate}
\end{prop}
\begin{proof}
    Part (1) is an immediate consequence of \Cref{proposition: marked sheaves
    affine over nonmarked sheaves}. For Part (2), let $R$ be a discrete
    valuation ring with generic point $\eta$, and choose a morphism $h: \Spec(R)
    \to \prestable_{g,n,p} \times \SVect^p$ along with a lift $f: \eta \to
    \Pur_{g,n,p}$. We may think of $h$ as parameterizing a tuple $(C,
    \sigma_{\bullet}^{\pm}, W_{\bullet})$ where $(C, \sigma_{\bullet}^{\pm})$ is
    an $R$-family of $(n,p)$-marked curves and $W_{\bullet}$ is a $p$-tuple of
    locally-free sheaves on $\Spec(R)$. By \Cref{proposition: marked sheaves
    affine over nonmarked sheaves}, we may think of the lift $f$ as a the data
    of a tuple $(K, \phi_{\bullet})$, where $K$ is a pure sheaf on the generic
    fiber $C_{\eta}$ and $\phi_{\bullet}$ is a $p$-tuple of homomorphisms
    $\phi_i: W_i|_{\eta} \to (\sigma_{i,\eta}^+)^*(K(\sigma_{i,\eta}^+))$. By
    \cite[Prop. 6.16]{torsion-freepaper}, there is an $R$-family $\widetilde{K}$
    of pure sheaves on $C$ along with an identification
    $\widetilde{K}|_{C_{\eta}} \cong K$. After replacing this identification
    with a scalar multiple, we may assume that the homomorphisms $\phi_i$ extend
    to homomorphisms $\widetilde{\phi}_i: W_i \to
    ((\sigma_{i}^+)^*(\widetilde{K}(\sigma_{i}^+)))$. \qedhere

\end{proof}

\subsection{The numerical invariant}

In order to define a notion of semistability for bundles on prestable curves, we
will need our curves to be polarized. Recall our notation $\Pur^{\pol}_{g,n,p}
:= \Pur_{g,n,p} \times_{\prestable_{g,n,p}} \prestablepol_{g,n,p}$.

\subsubsection{Numerical invariant for (non-marked) pure sheaves.}
\begin{notn}[Universal family]
    Let $\pi: (\cC, \sigma_{\bullet}^{\pm}, \cH) \to \Pur^{\pol}_{g,n+p,0}$
    denote the universal polarized curve, and let $F_{univ}$ denote the
    universal family of pure sheaves on $\cC$.
\end{notn}
\begin{defn}\label{defn: giesker line bundle coherent sheaves} We define the
    line bundle $\cL_{\deg}:= \det(R\pi_*(F_{univ} \otimes \omega_{\pi} \oplus
    F_{univ}))^{\vee}$ on $\Pur^{\pol}_{g,n+p,0}$. Similarly, we define
    $\cL_{\rk}:= \det(R\pi_*(F_{univ} \otimes\cH)) \otimes
    \det(R\pi_*(F_{univ}))^{\vee}$ on $\Pur^{\pol}_{g,n+p,0}$.
\end{defn}

\begin{remark}[Comparison with \cite{torsion-freepaper}]
    For the purposes of constructing numerical invariants, the line bundle
    $\cL_{\deg}$ in \Cref{defn: giesker line bundle coherent sheaves} could be
    replaced with a power of the line bundle $M_1^{\vee}$ defined in \cite[Defn.
    2.6]{torsion-freepaper}. On the other hand, we have $\cL_{\rk} := M_1
    \otimes M_0^{\vee} = b_1$ in the notation of \cite[Defn.
    2.9]{torsion-freepaper}.
\end{remark}

\begin{notn}[Hilbert polynomials and ranks] \label{notn: Hilbert polynomial and
    rank of pure sheaves} For any geometric point $(C, \sigma_{\bullet}^{\pm},
    \cH, F)$ of $\Pur^{\pol}_{g,n+p,0}$, we write $P_F(n)$ for the corresponding
    Hilbert polynomial with respect to $\cH$. This is a degree $1$ polynomial
    with integer coefficients; the degree $1$ coefficient of $P_F(n)$ is called
    the rank of $F$ (with respect to the polarization $\cH$), and we denote it
    by $\rk(F)$.
\end{notn}

Given a field $K$ and an integer $q \geq 1$, a graded point $(B\mathbb{G}_m^q)_K
\to \Pur^{\pol}_{g,n+p,0}$ relative to $\prestablepol_{g,n+p,0}$ amounts to the
data of a polarized $(n+p)$-pointed curve $(C, \sigma_{\bullet}, \cH)$ over $K$,
and a pure sheaf $F = \bigoplus_{\vec{m} \in \mathbb{Z}^q} F_{\vec{m}}$ of
dimension $1$ on $C$ equipped with a $\mathbb{Z}^q$-grading. For the following
definition, we refer to the notion of quadratic form on graded points in
\cite[Defn. 4.1.12]{halpernleistner2018structure}.
\begin{defn}[{\cite[Defn. 4.2]{torsion-freepaper}}] \label{defn: quadratic norm
    on coherent sheaves} We denote by $b$ the quadratic form on graded points of
    $\Pur^{\pol}_{g,n+p,0}$ relative to $\prestablepol_{g,n+p,0}$ such that, for
    every field $K$ and relative graded point $h: (B\mathbb{G}_m^q)_K \to
    \Pur^{\pol}_{g,n+p,0}$ corresponding to a tuple $(C, \sigma_{\bullet}, \cH,
    F = \bigoplus_{\vec{m} \in \mathbb{Z}^q} F_{\vec{m}})$, assigns the
    quadratic form $b(h)$ on $\mathbb{R}^q$ given by
    \[ b(h): \; \; \vec{r} \mapsto \sum_{\vec{m} \in \mathbb{Z}}
    \rk(F_{\vec{m}})^2 (\vec{r} \cdot_{std} \vec{m}),\] where $(-) \cdot_{std}
    (-)$ is the standard inner product on $\mathbb{R}^q$.
\end{defn}

The following simplifying notation will be useful for describing the numerical
invariant.
\begin{defn}[Slope and degree of a sheaf] \label{defn: degree and slope of
    sheaves on prestable curves} Let $C$ be a prestable curve over an
    algebraically closed field, and let $F$ be a coherent sheaf on $C$. We set
    $\deg(F):= \frac{1}{2} \cdot \left(\chi(F\otimes \omega_C) +
    \chi(F)\right)$. If the curve $C$ is equipped with a polarization $H$ and
    the sheaf $F$ satisfies $\rk(F) \neq 0$, then we set $\mu(F):= 2
    \deg(F)/\rk(F).$
\end{defn}

\begin{defn} \label{defn: rational line bundle L mu} Given a nonnegative
    rational number $\alpha \geq 0$, let $\cL_{\alpha}$ denote the rational line
    bundle on $\Pur^{\pol}_{g,n+p,0}$ given by $\cL_{\alpha} := -\cL_{\deg} -
    (\mu(-) + \alpha)\cdot \cL_{\rk}$, where $\mu(-): |\Pur^{\pol}_{g,n+p,0}|
    \to \mathbb{Q}$ is the locally constant function that sends a geometric
    point $(C, \sigma_{\bullet}, F)$ to $\mu(F)$.
\end{defn}

\begin{defn}
    We define $\nu_{\alpha}$ to be the $\mathbb{R}$-valued numerical invariant
    on $\Pur^{\pol}_{g,n+p,0}$ relative to $\prestablepol_{g,n+p,0}$ given by
    $\nu_{\alpha} := \wt(\cL_{\alpha})/\sqrt{b}$.
\end{defn}

For any given field $K$ and a $K$-point $x=(C, \sigma_{\bullet}, \cH, F)$ of
$\Pur^{\pol}_{g,n+p,0}$, a filtration $f:\Theta_K \to \Pur^{\pol}_{g,n+p,0}$ of
$x$ relative to $\prestablepol_{g,n+p,0}$ amounts to the data of a
$\mathbb{Z}$-indexed filtration $(F_m)_{m \in \mathbb{Z}}$ of subsheaves $F_i
\subset F$ satisfying the following
\begin{enumerate}[(1)]
    \item $F_{m+1} \subset F_m$ for all $m$; $F_m=0$ for $m \gg0$ and $F_m = F$
    for $m \ll0$.
    \item $F_m/F_{m+1}$ is pure of dimension $1$ for all $m$.
\end{enumerate}
See \cite[Prop. 1.0.1]{halpernleistner2018structure} or \cite[Prop.
3.8]{rho-sheaves-paper} for explanations of this using the Rees construction. A
direct computation, similar to the one at the beginning of Section 4.1 in
\cite{torsion-freepaper} shows that the value of the numerical invariant
$\nu_{\alpha}$ at such a filtration $f$ is given by
\begin{equation} \label{eqn: numerical invariant for pure sheaves}
    \nu_{\alpha}(f) = \frac{ \sum_{m \in \bZ} m \cdot \left(2 \deg(F_m/F_{m+1}) - \rk(F_m/F_{m+1}) \cdot (\mu(F)+ \alpha)\right)}{\sqrt{\sum_{m \in \bZ} m^2 \cdot \rk(F_m/F_{m+1}) }}
\end{equation}

\subsubsection{Numerical invariant for pure marked sheaves}

Recall from \Cref{prop: affineness of pur and valuative criterion for
properness} that the morphism $\Psi \times \ev^+: \Pur^{\pol}_{e,n,p} \to
\Pur^{\pol}_{g,n+p,0} \times \SVect^p$ is affine and of finite type. This makes
the following definition meaningful.
\begin{defn} \label{defn: rational norm on Pur} Let $b$ denote the quadratic
    norm on graded points of $\Pur^{\pol}_{g,n,p} \to \SVect^p$ relative to
    $\prestablepol_{g,n+p,0}\times\SVect^p$ induced by pulling from the first
    coordinate the quadratic norm defined in \Cref{defn: quadratic norm on
    coherent sheaves}. We also denote by $b$ the corresponding quadratic norm on
    graded points of $\Pur^{\pol}_{g,n,p}$ relative to $\prestablepol_{g,n,p}
    \times \SVect^p$ obtained by pulling back under the morphism $\Psi \times
    \ev^+$ described above.
\end{defn}

\begin{notn} \label{notn: line bundles on Pur} We denote by $\cL_{\deg},
    \cL_{\rk}$ and $\cL_{\alpha}$ the rational line bundles on
    $\Pur^{\pol}_{g,n,p}$ obtained by pulling back $\cL_{\deg}, \cL_{\rk}$ and
    $\cL_{\alpha}$ from \Cref{defn: giesker line bundle coherent sheaves} and
    \Cref{defn: rational line bundle L mu} via the morphisms $\Psi:
    \Pur^{\pol}_{g,n,p} \to \Pur^{\pol}_{g,n+p,0}$.
\end{notn}

\begin{defn} \label{defn: numerical invariants in marked pure sheaves and
    Gbunpol} We define the $\mathbb{R}$-valued numerical invariant $\nu$ on
    $\Pur^{\pol}_{g,n,p}$ relative to $\prestablepol_{g,n,p}\times\SVect^p$
    given by $\nu = \wt(\cL_{0})/\sqrt{b}$.
\end{defn}

\subsection{Monotonicity}

We now recall several codimension-two filling conditions on a morphism $\cM \to
\cY$, which were introduced in \cite[Sect. 5.2]{halpernleistner2018structure}.
All three will be formulated in terms of a regular 2-dimensional affine
Noetherian scheme $Y$ with an action of a torus $T$ (either $\bG_m$ or
$\bG_m^2$) and a unique $T$-fixed point $\mathfrak{o} \in Y$:
\begin{notn}\label{N:monotonicity_notation} \;

    \begin{enumerate}
        \item For strict $\Theta$-monotonicity: $Y = \Spec(R[t])$ for a complete
        discrete valuation ring $R$. $T = \bG_m$ acts on $Y$ by letting $t$ have
        weight $-1$.
        \item For strict $S$-monotonicity: $Y = \Spec(R[s,t]/(st-\pi))$ with $R$
        as in (1) and $\pi \in R$ a uniformizer. $T=\bG_m$ acts by giving $t$
        weight $-1$ and $s$ weight $1$.
        \item For strict $\Theta^2$-monotonicity: $Y = \bA^2_k = \Spec(k[s,t])$
        for a field $k$. $T = \bG_m^2$ acts on $Y$ by giving $s$ bidegree
        $(-1,0)$ and $t$ bidegree $(0,-1)$.
    \end{enumerate}
\end{notn}

In cases (1) and (2), the affine GIT quotient $Y /\!/ T \cong \Spec(R)$, and in
case (3) $Y/\!/T \cong \Spec(k)$.

\begin{defn}[Relative monotonicity conditions]\label{D:monotonicity} If $\cM \to
    \cY$ is a morphism of algebraic stacks, then we say that a numerical
    invariant $\mu$ (resp. line bundle $\cL$) on $\cM$ is strictly
    $\Theta$-monotone \emph{relative} to $\cY$ if for all $Y$ as in
    \Cref{N:monotonicity_notation}(1) and any commutative diagram
    \[
        \xymatrix{(Y \setminus \mathfrak{o}) / T \ar[d] \ar[r]^-f & \cM \ar[d] \\ Y/\!/T \ar[r] & \cY},
    \]
    after possibly replacing $R$ with an \'etale extension of discrete valuation
    rings, there is a $T$-equivariant proper morphism $\Sigma \to Y$ from an
    integral DM stack $\Sigma$ that is an isomorphism over $Y \setminus
    \mathfrak{o}$, and an extension of $f$ to a morphism $\tilde{f} : \Sigma /
    \bG_m \to \cM$ over $\cY$ such that if $x_1,\ldots,x_n$ are the $T$-fixed
    points in the fiber $|\Sigma_{\mathfrak{o}}|$ with their canonical ordering
    (see \cite[Lemma 5.2.6]{halpernleistner2018structure}), then
    $\mu(\{x_i\}/\bG_m \to \cM)$ is strictly monotone increasing in $i$ (resp.
    the weights $\wt(\cL|_{\{x_i\}/\mathbb{G}_m})$) are strictly monotone
    increasing in $i$).

    Strict $S$-monotonicity of $\cL$ or $\mu$ is defined by the same condition
    but for all $Y$ as in \Cref{N:monotonicity_notation}(2).

    Strict $\Theta^2$-monotonicity of $\cL$ or $\mu$ is defined as the existence
    of an extension $\tilde{f} : \Sigma / \bG_m^2 \to \cM$ for all diagrams with
    $Y$ as in \Cref{N:monotonicity_notation}(3), such that the monotonicity
    condition holds after composing $\tilde{f}$ with the canonical morphism
    $\Sigma'/\bG_m \to \Sigma / \bG_m^2$ from the base change $\Sigma' = \Sigma
    \times_{\bA^2_k} \Spec(k[\![s]\!][t])$, equipped with the action of $\{1\}
    \times \bG_m \subset \bG_m^2$ that gives $t$ weight $-1$.
\end{defn}

\begin{remark}
    In \Cref{D:monotonicity}, the $\bG_m$-fixed points of $\Sigma$ are ordered
    such that the condition that $\wt(\cL|_{\{x_i\}/\bG_m})$ is strictly
    monotone increasing in $i$ is equivalent to the condition that
    $\tilde{f}^\ast(\cL)$ is anti-ample on $\Sigma$.
\end{remark}

The strict monotonicity conditions are smooth-local on $\cY$. Note also that if
$\cY$ is an algebraic space, then the morphism $Y/T \to \cY$ factors uniquely
through $Y/\!/T$, so strict monotonicity can be reformulated as an ``absolute"
condition in this case.

To prove monotonicity results for the moduli of pure sheaves, we will use the
technique of infinite dimensional GIT from \cite{torsion-freepaper,
gauged_theta_stratifications}, which yields the following.
\begin{lemma} \label{lemma: rational filling} Let $R$ be a discrete valuation
    ring, and let $Y$ be as in \Cref{N:monotonicity_notation}(1) or (2). Let
    $(C, \sigma_{\bullet}, \cH)$ be an $R$-point of
    $\prestable^{\pol}_{g,n+p,0}$, and consider a $\mathbb{G}_m$-equivariant
    pure sheaf $K$ on the base-change $C_{Y \setminus \mathfrak{o}}$
    corresponding to a morphism $h: (Y \setminus \mathfrak{o}) /\mathbb{G}_m \to
    \Pur^{\pol}_{g,n+p,0}$ such that the composition $(Y \setminus \mathfrak{o})
    /\mathbb{G}_m \to \Pur^{\pol}_{g,n+p,0}$ recovers the $R$-point $(C,
    \sigma_{\bullet}, \cH)$.

    Then, there exists a projective $\mathbb{G}_m$-equivariant morphism $\Sigma
    \to Y$ from an integral scheme $\Sigma$ that is an isomorphism over $Y
    \setminus \mathfrak{o}$ and an extension $\widetilde{h}: \Sigma/\mathbb{G}_m
    \to \Pur^{\pol}_{g,n+p,0}$ of $h$ corresponding to a
    $\mathbb{G}_m$-equivariant pure sheaf $\widetilde{K}$ on $C_{\Sigma}$
    satisfying the following:
    \begin{enumerate}
        \item There is an open dense subscheme $\overline{U} \subset
        C_{\mathfrak{o}}$ such that the restriction of
        $\widetilde{K}_{\mathfrak{o}}$ to $\overline{U} \times
        (\Sigma_{\mathfrak{o}}/\mathbb{G}_m)$ is isomorphic to the pullback of a
        sheaf on $\overline{U} \times (\mathfrak{o}/\mathbb{G}_m)$.

        \item The pullback $\widetilde{h}_{\mathfrak{o}}^*(\cL_{\alpha})$ under
        the restriction $\widetilde{h}_{\mathfrak{o}}: \Sigma_{\mathfrak{o}} \to
        \Pur_{g,n+p,0}^{\pol}$ is anti-ample.
    \end{enumerate}
\end{lemma}
\begin{proof}
    In the proof of \cite[Thm. 4.11]{torsion-freepaper}, one obtains the
    following data:
    \begin{itemize}
        \item A finite morphism $q: C \to \mathbb{P}^1_R$.
        \item A $\mathbb{G}_m$-equivariant effective $Y$-flat Cartier divisor $D
        \subset \mathbb{P}^1_Y$ such that its preimage $q_Y^{-1}(D)$ under the
        base-change morphism $q_Y: C_Y \to \mathbb{P}^1_Y$ is cut out by a
        $\mathbb{G}_m$-equivariant section of a power of the ample line bundle
        $\cH|_{C_Y}$. For simplicity of notation, we set $V := \mathbb{P}^1_Y
        \setminus D$ and let $U:= (q_Y)^{-1}(V) \subset C_Y$.
        \item A $\mathbb{G}_m$-equivariant vector bundle $E$ on
        $\mathbb{P}^1_R$, and a $\mathbb{G}_m$-equivariant
        $(q_Y)_*(\cO_{U})$-module structure on the restriction $E|_V$.
        \item A $\mathbb{G}_m$-equivariant isomorphism of restrictions $E|_{V_{Y
        \setminus \mathfrak{o}}} \xrightarrow{\sim} (q_Y)_*(K|_{U_{Y \setminus
        \mathfrak{o}}})$ that is compatible with the corresponding
        $(q_Y)_*(\cO_{U_{Y \setminus \mathfrak{o}}})$-module structures.
    \end{itemize}

    Using $E$ and the $(q_Y)_*(\cO_{U})$-module structure on $E|_V$, one defines
    an ind-projective affine grassmannian $\Gr \to Y$ for the moduli of
    $(q_Y)_*(\cO_{C_Y})$-modules on $C_Y \to Y$ as in \cite[Sect.
    3.5]{torsion-freepaper}, which parametrizes families $G$ of pure sheaves on
    $C_Y$ along with an isomorphism of $E|_V \xrightarrow{\sim} (q_Y)_*(G|_{U})$
    which is compatible with the corresponding $(q_Y)_*(\cO_U)$-module
    structures. The data of $K$ and the isomorphism $E|_{V_{Y \setminus
    \mathfrak{o}}} \xrightarrow{\sim} (q_Y)_*(K|_{U_{Y \setminus
    \mathfrak{o}}})$ yields a section $Y \setminus \mathfrak{o} \to \Gr$ defined
    over $Y\setminus \mathfrak{o}$. If we set $\Sigma \hookrightarrow \Gr$ to be
    the scheme-theoretic closure of this section, then there is a
    $\mathbb{G}_m$-equivariant $\Sigma$-family $\widetilde{K}$ of pure sheaves
    on $C_{\Sigma}$ along with an equivariant isomorphism $E|_{V_{\Sigma}}
    \xrightarrow{\sim} (q_Y)_*(K|_{U_{\Sigma}})$. In particular, after we
    further restrict to $Y \setminus \mathfrak{o} = \Sigma_{Y \setminus
    \mathfrak{o}}$, we obtain an equivariant isomorphism $\widetilde{K}|_{C_{Y
    \setminus \mathfrak{o}}} \cong K$.

    As a summary, the output of the construction in the proof of \cite[Thm.
    4.11]{torsion-freepaper} is an integral $\mathbb{G}_m$-scheme $\Sigma$, a
    proper $\mathbb{G}_m$-equivariant morphism $\Sigma \to Y$ that is an
    isomorphism over $Y \setminus \mathfrak{o}$, and a
    $\mathbb{G}_m$-equivariant $\Sigma$-family of sheaves $\widetilde{K}$ on
    $C_{\Sigma}$ along with an identification $\widetilde{K}|_{C_{Y \setminus
    0}} \cong K$. We may view this as a morphism $\widetilde{h}:
    \Sigma/\mathbb{G}_m \to \Pur^{\pol}_{g,n+p,0}$ extending $h$. We claim that,
    up to possibly replacing $D$ with a larger $Y$-flat divisor in the
    construction above, we may arrange (1) and (2) as in the statement of the
    lemma.

    To see part (1), note that this is equivalent to saying that there is an
    open dense $\overline{U} \subset C_{\mathfrak{o}}$ and a graded sheaf
    $\bigoplus_{w \in \mathbb{Z}} F_w$ on $\overline{U}$ such that
    $\widetilde{K}|_{\overline{U} \times (\Sigma_{\mathfrak{o}}/\mathbb{G}_m)}$
    is the constant family of equivariant sheaves corresponding to the graded
    sheaf $\bigoplus_{w \in \mathbb{Z}} F_{w}|_{\overline{U} \times
    \Sigma_\mathfrak{o}}$. This follows directly from the construction if we
    take $\overline{U} := U_{\mathfrak{o}}$.

    To prove part (2), we first note that there is an $R$-flat effective Cartier
    divisor $D' \subset \mathbb{P}^1_R$ such that the restriction of the
    relative canonical line bundle $\omega_{C/R}$ to $C \setminus q^{-1}(D')$ is
    trivializable. \endnote{Indeed, the consider the limit $S$ of the
    complements of all $R$-flat effective Cartier divisors on $\mathbb{P}^1_R$
    (the inclusions are affine, so $S$ is a scheme equipped with a monomorphism
    to $\mathbb{P}^1_R$). One may think of the limit $S$ by first removing the
    infinity section from $\mathbb{P}^1_R$ to get $\mathbb{A}^1_R$, and then
    further removing the remaining $R$-flat Cartier divisors amounts to the
    localization $R[t]_{(\varpi)}$ at the prime $(\varpi) \subset R[t]$, where
    $\varpi$ is a uniformizer. Hence, the limit $S$ is the spectrum of a local
    ring, and the fiber product of $S$ with $C$ is the spectrum of a semilocal
    ring. In particular, it follows that the restriction of $\omega_{C/R}$ to
    the limit $S \times_{\mathbb{P}^1_R} C$ is trivializable. A standard
    spreading out argument then yields the existence of the desired $R$-flat
    effective Cartier divisor $D' \subset \mathbb{P}^1_R$ such that the
    restriction of the relative canonical line bundle $\omega_{C/R}$ to $C
    \setminus q^{-1}(D')$ is trivializable.}

    We may replace the effective Cartier divisor $D$ in the infinite dimensional
    GIT construction with $D \cup D'_Y \subset \mathbb{P}^1_Y$ in order to
    ensure that the relative canonical bundle $\omega_{C_Y/Y}$ is
    $\mathbb{G}_m$-equivariantly trivializable when restricted to $U \subset
    C_Y$. Fixing such a trivialization, we get a $\mathbb{G}_m$-equivariant
    isomorphism $\theta: \Gr \to \Gr$ which sends a pair $(G, \psi)$ of a family
    $G$ with isomorphism $\psi: E|_V \xrightarrow{\sim} (q_Y)_*(G|_{U})$ to the
    pair $(G \otimes \omega_{C_Y/Y}, \psi')$, where $\psi'$ is induced from
    $\psi$ and the chosen trivialization of $\omega_{C_Y/Y}$ to $U$. Consider
    the families of line bundles $M_m$ indexed by integers $m \in \mathbb{Z}$
    and the line bundles $b_0,b_1$ on $\Gr$ defined in the first paragraph of
    \cite[Sect. 4]{torsion-freepaper}. Then the pullback
    $\widetilde{h}_{\mathfrak{o}}^*(\cL_{\rk})$ coincides with the restriction
    $b_1|_{\Sigma_{\mathfrak{o}}}$ of the line bundle $b_1$, which is torsion
    when restricted to any projective stratum of $\Gr$ by \cite[Lemma
    3.23]{torsion-freepaper} combined with \cite[Prop. 3.27]{torsion-freepaper}.
    We conclude that the desired anti-ampleness of
    $\widetilde{h}_{\mathfrak{o}}^*(\cL_{\alpha}) =
    \widetilde{h}_{\mathfrak{o}}^*(\cL_{\deg})^{\vee} \otimes
    \widetilde{h}_{\mathfrak{o}}^*(\cL_{\rk})^{- \otimes (\mu(-)+ \alpha)}$ is
    equivalent to the anti-ampleness of
    $\widetilde{h}_{\mathfrak{o}}^*(\cL_{\deg})^{\vee}$. By definition
    $\widetilde{h}_{\mathfrak{o}}^*(\cL_{\deg})^{\vee}$ coincides with the
    restriction $(\theta^*(M_0) \otimes M_0)|_{\Sigma_{\mathfrak{o}}})$ of the
    line bundle $\theta^*(M_0) \otimes M_0$ on $\Gr$. Hence, it is enough to
    show that $M_0$ is anti-ample on each projective stratum of $\Gr$. This
    follows by the anti-ampleness of $M_m$ on each projective stratum of $\Gr$
    for $m \gg 0$ (see the proof of \cite[Lemma 3.23]{torsion-freepaper}), the
    equation $M_n = b_2^{\otimes \binom{m}{2}} \otimes b_1^{\otimes m} \otimes
    M_0$ from \cite[Defn. 2.9]{torsion-freepaper}, and the fact that both $b_1$
    and $b_2$ are torsion on each projective stratum of $\Gr$ (by the proof
    \cite[Lemma 3.23]{torsion-freepaper}).
\end{proof}

\begin{prop} \label{prop: monotonicity for nu on marked sheaves} Let $\alpha
    \geq 0$ be a rational number. The numerical invariant $\nu$ of \Cref{defn:
    numerical invariants in marked pure sheaves and Gbunpol} and the line bundle
    $\cL_{\alpha}$ of \Cref{defn: rational line bundle L mu} are strictly
    $\Theta$-monotone, strictly $S$-monotone, and strictly $\Theta^2$-monotone
    on $\Pur^{\pol}_{g,n,p}$ relative to $\prestablepol_{g,n,p} \times
    \SVect^p$.
\end{prop}
\begin{proof}
    We omit the proof of $\Theta^2$-monotonicity, as it is completely analogous
    to the proof of $\Theta$- and $S$- monotonicity below. The ranks
    $N_1,\ldots,N_p$ of the locally free sheaves $W_i$ in \Cref{defn: stack of
    marked sheaves} are locally constant, so we shall fix a choice of these
    ranks for the remainder of the proof, and temporarily use
    $\Pur^{\pol}_{g,n,p}$ to denote the corresponding open and closed preimage
    of $\prod_{i=1}^p \BGL_{N_i} \subset \SVect^p$ under $\ev^+$. We let $S \to
    \prod_{i=1}^p \BGL_{N_i}$ be the morphism classifying the trivial locally
    free sheaves of ranks $N_1,\ldots,N_p$, and set $\Pur^\ast_{g,n,p} := S
    \times_{\prod_{i=1}^p \BGL_{N_i}} \Pur^{\pol}_{g,n,p}$. Then
    \Cref{proposition: marked sheaves affine over nonmarked sheaves} gives an
    isomorphism
    \[
        \Pur^{\ast}_{g,n,p} \cong \Hom(\cO^{\oplus N_1},\cK_1) \times_{\Pur^{\pol}_{g,n+p,0}} \cdots \times_{\Pur^{\pol}_{g,n+p,0}} \Hom(\cO^{\oplus N_p}, \cK_p)
    \]
    relative to $\Pur^{\pol}_{g,n+p,0}$. It suffices to show that the
    corresponding pulled back numerical invariant $\nu$ and line bundle
    $\cL_{\alpha}$ on $\Pur^\ast_{g,n,p}$ are strictly $\Theta$-monotone and
    strictly $S$-monotone relative to $ \prestable^{\pol}_{g,n,p}$.

    Let $R$ be a discrete valuation ring and set $Y$ as in
    \Cref{N:monotonicity_notation}(1) or (2). Choose $f: (Y \setminus
    \mathfrak{o})/\mathbb{G}_m \to \Pur^{\ast}_{g,n,p}$ such that the
    composition $(Y \setminus \mathfrak{o})/\mathbb{G}_m \to \Pur^{\ast}_{g,n,p}
    \to \prestable_{g,n,p}^{\pol}$ factors through a morphism $\Spec(R) \to
    \prestable_{g,n,p}^{\pol}$, which corresponds to an $R$-family $(C,
    \sigma_{\bullet}^{\pm}, \cH)$ of polarized $(n,p)$-marked curves. We may
    think of $f$ as a pair $(K, \phi)$, where $K$ is a
    $\mathbb{G}_m$-equivariant $(Y \setminus \mathfrak{o})$-family of pure
    sheaves on the base-changed family of $(n,p)$-marked curves $(C_{Y \setminus
    \mathfrak{o}}, \sigma^{\pm}_{\bullet, Y \setminus \mathfrak{o}})$ and $\phi$
    is a $\mathbb{G}_m$-equivariant homomorphism $\phi: \cO_{Y \setminus
    \mathfrak{o}} \to \bigoplus_{i=1}^p K_i$, where we define $K_i :=
    (\sigma_{i, Y \setminus \mathfrak{o}}^+)^*K(\sigma_{i, Y \setminus
    \mathfrak{o}}^+)$. Consider the composition $h: (Y \setminus
    \mathfrak{o})/\mathbb{G}_m \xrightarrow{f} \Pur^{\ast}_{g,n,p}
    \xrightarrow{\Psi} \Pur^{\pol}_{g,n+p,0}$, and let $\widetilde{h}:
    \Sigma/\mathbb{G}_m \to \Pur^{\pol}_{g,n+p,0}$ denote an extension as in
    \Cref{lemma: rational filling}. Let $\bigoplus_{w \in \mathbb{Z}} F_w$
    denote the graded sheaf on $\overline{U}$ as in part (1) \Cref{lemma:
    rational filling}. Set $\rk(F_w)$ to be the rank with respect to
    $\cH_{\mathfrak{o}}$ of any coherent extension of $F_w$ to
    $C_{\mathfrak{o}}$, which is independent of the chosen extension because
    $\overline{U} \subset C_{\mathfrak{o}}$ is dense.

    \noindent $\bullet$ \textbf{Monotonicity for $\nu$.} Let $x_0, x_1, \ldots,
    x_{\ell}$ denote the $\mathbb{G}_m$-fixed points in
    $|\Sigma_{\mathfrak{o}}|$, ordered as in \cite[Lemma
    5.2.6]{halpernleistner2018structure}. Let us denote by $a_i$ the
    $\mathbb{G}_m$-weight of the fiber
    $\widetilde{h}_{\mathfrak{o}}^*(\cL_{0})|_{x_i}$. Then it follows from part
    (2) in \Cref{lemma: rational filling} that the sequence $a_0, a_1, \ldots,
    a_\ell$ is monotone increasing. On the other hand, for each index $i$ the
    graded point $g_i: x_i/\mathbb{G}_m \to \Sigma_{\mathfrak{o}}/\mathbb{G}_m
    \xrightarrow{\widetilde{h}_{\mathfrak{o}}} \Pur^{\pol}_{g,n+p,0}$
    corresponds to a graded pure sheaf $\bigoplus_{w \in \mathbb{Z}}
    \widetilde{F}^i_{w}$ on $C_{\mathfrak{o}}$ extending the graded sheaf
    $\bigoplus_{w \in \mathbb{Z}} F_w$ on $\overline{U}$ as above. It follows
    that the value of the norm $b(x_i)$ as in \Cref{defn: quadratic norm on
    coherent sheaves} is given by $\sum_{w \in \mathbb{Z}} w^2 \cdot \rk(F_w)$,
    and hence it does not depend on $i$.

    Let $E_1,\ldots, E_{\ell}$ be the irreducible components of
    $\Sigma_{\mathfrak{o}}$, where $E_i$ connects $x_i$ with $x_{i-1}$. Up to
    blowing up $\Sigma$, we may assume that each $E_i$ is a Cartier divisor on
    $\Sigma$. For any sequence of integers $r_1, r_2, \ldots, r_\ell$, we may
    replace $\widetilde{K}$ with $\widetilde{K} \otimes \cO_{\Sigma}(r_1
    E_1+\cdots+r_n E_n)$ to obtain a new morphism $\widetilde{h}_{r_1, \ldots,
    r_\ell}: \Sigma / \bG_m \to \Pur^{\pol}_{g,n+p,0}$ whose restriction to $Y
    \setminus \mathfrak{o}$ agrees with $h$. If we choose all $r_i$ to be
    sufficiently large, then the morphism $\phi: \cO_{Y \setminus \mathfrak{o}}
    \to \bigoplus_{i=1}^p K_i$ extends to a morphism $\phi: \cO_{\Sigma} \to
    \widetilde{K}_i$, where the $\widetilde{K}_i= (\sigma_{i,
    \Sigma}^+)^*\widetilde{K}(\sigma_{i, \Sigma}^+)$. In other words, we obtain
    a lift $\widetilde{f}_{r_1, \ldots, r_\ell}: \Sigma/\mathbb{G}_m \to
    \Pur^{\ast}_{g,n,p}$ for all choices of sufficiently large integers $r_i$.

    To conclude the proof, we analyze how the choice of modification
    $h':=\widetilde{h}_{r_1, \ldots, r_i}$ affects the weights $a_i$ and the
    norms $b_i$.  If we denote by $g_i': x_i/\mathbb{G}_m \to
    \Sigma_{\mathfrak{o}}/\mathbb{G}_m
    \xrightarrow{\widetilde{h}'_{\mathfrak{o}}} \Pur^{\pol}_{g,n+p,0}$ the
    corresponding modification, then $g_i'$ corresponds to the graded sheaf
    $\bigoplus_{w \in \mathbb{Z}} \widetilde{G}^{i}_w$, where we set
    $\widetilde{G}^{i}_w := \widetilde{F}_{w+(r_i-r_{i+1})}$ under the
    convention that $r_0 = r_{\ell+1} = 0$. If we set $\mu =
    \mu(\widetilde{K}_{\mathfrak{o}})$ and $\rk =
    \rk(\widetilde{K}_{\mathfrak{o}})$, then by Equation \eqref{eqn: numerical
    invariant for pure sheaves}, the $\mathbb{G}_m$-weight $a_i'$ of
    $(g_i')^*(\cL_{0})$ is
    \[ a_i' = \sum_{w \in \mathbb{Z}} w \cdot ( 2\deg(\widetilde{G}^{i}_w) - \mu
            \cdot \rk(\widetilde{G}^{i }_w)).\] Using $\widetilde{G}^{i}_w =
            \widetilde{F}_{w+(r_i-r_{i+1})}$ and additivity of ranks and
            degrees, we see that
    \[ a_i' = \sum_{w \in \mathbb{Z}} w \cdot ( 2\deg(\widetilde{F}^{i}_w) - \mu
    \cdot \rk(\widetilde{F}^{i}_w)) = a_i.\] On the other hand, we have
    \[ b_i' = \sum_{w \in \mathbb{Z}} w^2 \cdot \rk(\widetilde{G}^{i}_w)=
    \left(\sum_{w \in \mathbb{Z}} w^2 \cdot \rk(\widetilde{F}^{i}_w)\right) -
    2(r_i - r_{i+1}) \cdot \left(\sum_{w \in \mathbb{Z}} w \cdot
    \rk(\widetilde{F}^{i}_w)\right) + \rk \cdot (r_i - r_{i+1})^2.\] If we set
    $b= b_0 = b_2 = \ldots = b_\ell$ and $A = -\sum_{w \in \mathbb{Z}} w \cdot
    \rk(\widetilde{F}^{i}_w) = -\sum_{w \in \mathbb{Z}} w \cdot \rk(F_w)$ (note
    that this does not depend on $i$) then we may rewrite this as follows:
    \[b_i' = b + 2A(r_i - r_{i+1}) + \rk \cdot (r_i - r_{i+1})^2.\] Hence, the
    modification $\widetilde{h}'$ satisfies the monotonicity condition for the
    numerical invariant $\nu$ if the sequence $\nu_0', \nu_1', \ldots,
    \nu_\ell'$ defined by
    \[ \nu_i' := \frac{a_i}{\sqrt{b + 2A(r_i - r_{i+1}) + \rk \cdot (r_i -
    r_{i+1})^2}} \] satisfies $\nu_0' < \nu_1' < \ldots < \nu_\ell'$. Let us
    denote by $0 \leq j \leq \ell$ the index such that $a_j <0 \leq a_{j+1}$
    (with the convention that $a_{-1}\ll 0$ and $a_{\ell+1} \gg0$). We
    distinguish two cases:

    \noindent \textit{Case 1.} If $j \neq \ell$, then we have $a_{\ell} \geq 0$.
    We can set $r_i := (i+1)R$ for a sufficiently large integer $R \gg 0$. Then,
    we have $b_0 = b_1 = \ldots = b_{\ell-1} < b_{\ell}$, and in view of the
    monotonicity of $a_i$, this guarantees $\nu_0' < \nu_1' < \ldots <
    \nu_\ell'$.

    \noindent \textit{Case 2.} If $j \neq \ell$, then we have $a_i <0$ for all
    $0 \leq i \leq \ell$. We may choose $r_i:= (\ell+1-i) R$ for a sufficiently
    large integer $R \gg 0$. Then, we have $b_0 > b_1 = \ldots = b_{\ell-1} =
    b_{\ell}$. Again, these inequalities for the norms $b_i$ and the
    monotonicity of the $a_i$ imply the desired monotonicity $\nu_0' < \nu_1' <
    \ldots < \nu_\ell'$.

    \noindent $\bullet$ \textbf{Monotonicity for $\cL_{\alpha}$.} Similarly as
    before, the sequence of $\mathbb{G}_m$-weights
    $a_i:=\widetilde{h}_{\mathfrak{o}}^*(\cL_{\alpha})|_{x_i}$ is monotone
    increasing by part (2) of \Cref{lemma: rational filling}. Just as in the
    proof of monotonicity for $\nu$, for all choices of sufficiently large
    integers $r_i$ we obtain a lift $\widetilde{f}_{r_1, \ldots, r_\ell}:
    \Sigma/\mathbb{G}_m \to \Pur^{\ast}_{g,n,p}$, and we need to analyze the new
    weights $a_i':= \wt((g_i')^*(\cL_{\alpha}))$. By \Cref{eqn: numerical
    invariant for pure sheaves}, we have
    \[a_i' = \sum_{w \in \mathbb{Z}} w \cdot ( 2\deg(\widetilde{G}^{i}_w) -
    (\mu+ \alpha) \cdot \rk(\widetilde{G}^{i}_w)),\] where $\widetilde{G}^{i}_w
    := \widetilde{F}_{w+(r_i-r_{i+1})}$ , where we recall that by convention we
    set $r_0 = r_{\ell+1} =0$. Using additivity of ranks and degrees, we get
    \[ a_i' = a_i + \rk \cdot \alpha \cdot (r_i - r_{i+1}),\] where $\rk :=
    \rk(\widetilde{K}_{\mathfrak{o}})>0$. To conclude the proof of monotonicity,
    we need to show that we can find $r_i$ arbitrarily large such that $a_0' <
    a_2'< \ldots < a_\ell'$. In other words, we want that for all $1 \leq i \leq
    \ell$ we have $a_{i}-a_{i-1} + \rk\cdot \alpha \cdot (2 r_{i} - r_{i+1}-
    r_{i-1}) >0$. Since the $a_i$ are already strictly monotone, we just need to
    arrange $2 r_{i} - r_{i+1}- r_{i-1} \geq 0$, which is satisfied as long as
    we pick the sequence $r_i$ to be convex. This can be arranged with endpoints
    $r_0 = r_{\ell+1} =0$ and $r_i$ arbitrarily large for $1 \leq i \leq \ell$.
\end{proof}

\subsection{HN boundedness}

\begin{defn}
    Given a morphism of quasi-separated locally Noetherian algebraic stacks $\cM
    \to \cY$ with affine relative automorphism groups, a numerical invariant
    $\mu$ on $\cM$ satisfies HN boundedness relative to $\cY$ if for any scheme
    $T$ of finite type over $\cY$ and commutative diagram
    \[
        \xymatrix{\Theta_T \ar[r]^f \ar[d] & \cM \ar[d] \\
            T \ar[r] & \cY},
    \]
    there is another morphism $f' : \Theta_{T'} \to \cM$ such that for any
    finite type point $t \in T$, there is a point $t' \in T'$ satisfying the
    following:
    \begin{enumerate}
        \item $t$ and $t'$ map to the same geometric point $y \in |\cY|$,
        \item $f_t(1) = f_{t'}(1)$ in $|\cM \times_{\cY} \Spec(k(y))|$,
        \item $\mu(f'_{t'}) \geq \mu(f_t)$, and
        \item the composition $f_{t'} : \Theta_{k(t')} \to \cM \to \cY$ factors
        through the projection $\Theta_{k(t')} \to \Spec(k(t'))$.
    \end{enumerate}
\end{defn}
In this subsection, we prove the following:

\begin{prop} \label{prop: hn boundedness} The numerical invariant $\nu$ on
    $\Pur^{\pol}_{g,n,p}$ from \Cref{defn: numerical invariants in marked pure
    sheaves and Gbunpol} satisfies HN boundedness relative to
    $\prestablepol_{g,n,p} \times \SVect^p$.
\end{prop}

Similarly as in the beginning of the proof of monotonicity (\Cref{prop:
monotonicity for nu on marked sheaves}), we may fix a $p$-tuple of ranks
$N_{\bullet}$ and consider the base-change $\Pur^{\ast}_{g,n,p}:=
\Pur^{\pol}_{g,n,p} \times_{\SVect^p} S$ under $S \to \prod_{i=1}^p \BGL_{N_i}
\hookrightarrow \SVect^p$. It is sufficient to prove HN boundedness for the
corresponding pullback numerical invariant $\nu$ on $\Pur^{\ast}_{g,n,p}$
relative to $\prestablepol_{g,n,p}$. Using the alternative description
\[\Pur^{*}_{g,n,p} \cong \Hom(\cO^{\oplus N_1},\cK_1)
\times_{\Pur^{\pol}_{g,n+p,0}} \cdots \times_{\Pur^{\pol}_{g,n+p,0}}
\Hom(\cO^{\oplus N_p}, \cK_p)\] from \Cref{proposition: marked sheaves affine
over nonmarked sheaves}, for any field $k$ we regard a $k$-point of
$\Pur^{\ast}_{g,n,p}$ as a pure sheaf $K$ on a polarized $(n,p)$-marked
prestable curve $(C, \sigma_{\bullet}^{\pm}, \cH)$ over $k$, along with a
$p$-tuple $\phi_{\bullet}$ of homomorphisms of vector spaces $\phi_i : k^{\oplus
N_i} \to (\sigma^+_i)^\ast(K(\sigma^+_i)) \cong (\sigma^+_i)^\ast(K)$.
\begin{lemma}\label{L:describe_filtration}

    A filtration of $(C,\sigma_{\bullet}^{\pm}, \cH, K,W_{\bullet}, \phi_{\bullet}) \in
    \Pur^*_{g,n,p}(k)$ relative to $\prestablepol_{g,n,p}$ corresponds to a
    $\bZ$-weighted filtration $(K_w)_{w \in \mathbb{Z}}$ of coherent subsheaves
    $K_w \subset K$ satisfying the following:
    \begin{enumerate}
        \item $K_{w+1} \subset K_w$ for all $w \in \mathbb{Z}$,
        \item $K_w=0$ for all sufficiently large $w \gg0$, and $K_w = K$ for all
        sufficiently small $w \ll 0$,
        \item $K_w/K_{w+1}$ is pure of dimension $1$ for all $w \in \mathbb{Z}$,
        and
        \item The homomorphism $\oplus_i \phi_i: \bigoplus_{i=1}^p k^{\oplus
        N_i} \to \bigoplus_{i=1}^p (\sigma^+_i)^\ast(K)$ factors through the
        subspace $\bigoplus_{i=1}^p (\sigma^+_i)^\ast(K_0) \subset
        \bigoplus_{i=1}^p (\sigma^+_i)^\ast(K)$.
    \end{enumerate}
\end{lemma}
\begin{proof}
    This is a consequence of the Rees construction (see \cite[Prop.
    1.0.1]{halpernleistner2018structure} and \cite[Prop.
    3.8]{rho-sheaves-paper}).\endnote{A filtration $\Theta_k \to \Pur^*_{g,n,p}$
    of $(C, \sigma^{\pm}, \cH, K, \phi)$ relative to $\prestablepol_{g,n,p}$ is
    a $\Theta_k$-family of pure sheaves $\widetilde{K}$ on $C \times \Theta_k$
    along with an isomorphism $\psi: K \xrightarrow{\sim} \widetilde{K}|_{1}$
    over the open dense substack $1 \subset \Theta_k$ such that for all $1 \leq
    i \leq p$ the morphism
    \[k^{\oplus N_i} \xrightarrow{\phi_i} (\sigma_i^+)^*(K)
    \xrightarrow{(\sigma_i^+)^{*}(\psi)} (\sigma_i^+)^*(\widetilde{K}|_1)\] of
    vector bundles on $1 \subset \Theta_k$ extends (uniquely) to a morphism of
    vector bundles $\cO_{\Theta_k}^{\oplus N_i} \to (\sigma_{i,
    \Theta_k}^+)^*(\widetilde{K})$ on $\Theta_k$. The Rees construction (see
    \cite[Prop. 1.0.1]{halpernleistner2018structure} and \cite[Prop.
    3.8]{rho-sheaves-paper}) establishes a natural bijection between the set of
    $\mathbb{Z}$-filtrations $(K_w)_{w \in \mathbb{Z}}$ satisfying conditions
    (1) through (3) above and the set of isomorphism classes of pairs
    $(\widetilde{K}, \psi)$ of $\Theta_k$-families $\widetilde{K}$ of pure
    sheaves and isomorphisms $\psi: K \xrightarrow{\sim} \widetilde{K}|_{1}$. It
    follows directly from the construction in \cite[Prop.
    1.0.1]{halpernleistner2018structure} that the existence of the required
    extensions $\cO_{\Theta_k}^{\oplus N_i} \to (\sigma_{i,
    \Theta_k}^+)^*(\widetilde{K})$ for a given pair $(\widetilde{K}, \psi)$ is
    equivalent to condition (4) of the corresponding $\mathbb{Z}$-filtration
    $(K_w)_{w \in \mathbb{Z}}$.}
\end{proof}

\begin{remark} \label{remark: unweighted plus jumping weights} We may
    equivalently think of a filtration $(K_w)_{w \in \mathbb{Z}}$ as in
    \Cref{L:describe_filtration} in terms of the following data:
    \begin{enumerate}[(a)]
        \item An ``underlying" finite unweighted filtration $0 \subsetneq
        G_{(q)} \subsetneq G_{(q-1)} \subsetneq \cdots \subsetneq G_{(0)} = K$;
        and
        \item a sequence of integers $w_0 < w_1 <\ldots <w_q$ (the ``jumping
        weights").
    \end{enumerate}
    Here the underlying unweighted filtration $G_{\bullet}$ is obtained by
    forgetting the $\mathbb{Z}$-weights and only keeping track of the finitely
    many distinct subsheaves of $K$ appearing in the filtration $(K_w)_{w \in
    \mathbb{Z}}$. The ``jumping weight" $w_i$ is the smallest integer such that
    $G_i = K_{w_i}$.
\end{remark}

For the following result, we use \Cref{L:describe_filtration} to view
filtrations $f$ of a $k$-point $(C, \sigma_{\bullet}^{\pm}, \cH,K, W_{\bullet},
\phi_{\bullet})$ of $\Pur^{\ast}_{g,n,p}$ as a $\mathbb{Z}$-weighted filtration $(K_w)_{w
\in \mathbb{Z}}$. For simplicity, we use the description of filtrations in terms
of unweighted filtrations and jumping weights as in \Cref{remark: unweighted
plus jumping weights}.
\begin{lemma} \label{lemma: technical optimization lemma hn boundednes marked
    sheaves} Let $f$ be a nondegenerate filtration of a field-valued point
    $x=(C, \sigma_{\bullet}^{\pm}, \cH,K, W_{\bullet}, \phi_{\bullet})$ in
    $\Pur^{\ast}_{g,n,p}$ such that $\nu(f)>0$. Then, there is another
    nondegenerate filtration $f'$ of $x$ with $\nu(f')\geq \nu(f)$ such that the
    following hold:
    \begin{enumerate}
        \item The underlying unweighted filtration $0 \subsetneq G_{(q)}
        \subsetneq G_{(q-1)} \subsetneq \cdots \subsetneq G_{(0)} = K$ with
        associated graded sheaves $\overline{G}_i := G_{(i)}/G_{(i+1)}$ has
        slopes $\mu_i := \mu(\overline{G}_i)$ that satisfy
              \[
                  \mu_{max}(K) \geq \mu_q>
                  \cdots > \mu_{j+1} > \max(\mu_j,\mu(K)) > \mu_{j-1} > \cdots > \mu_{0} \geq \mu_{min}(K),
              \]
              where $j$ is the largest index such that $\im(\oplus _i \phi_i)
              \subset \bigoplus_{i=1}^p (\sigma^+_i)^*(G_j)$. Here
              $\mu_{max}(K)$ (resp.  $\mu_{min}(K)$) denotes the maximum (resp.
              minimum) slope, as in \Cref{defn: degree and slope of sheaves on
              prestable curves}, among the graded pieces of the
              Harder-Narasimhan filtration of $K$ with respect to the slope
              function $\mu$.
        \item We have
              \[
                  \mu_j = \frac{1}{\rk(\overline{G}_j)} \cdot \left( \rk(K) \cdot \mu(K) - \sum_{i \neq j}\rk(\overline{G}_i) \cdot \mu_i \right).
              \]
        \item There is some positive rational number $Q>0$ such that $w_i = Q
        \cdot (\mu_i - \mu(K))$ for all $i \neq j$ and $w_j = Q \cdot
        \max\left(\mu_j - \mu(K), 0\right)$.
    \end{enumerate}
\end{lemma}
\begin{proof}
    In view of the description of the filtrations in
    \Cref{L:describe_filtration} and the numerical invariant in \eqref{eqn:
    numerical invariant for pure sheaves}, this lemma follows from an identical
    optimization argument as in \cite[Lemma 6.9]{torsion-freepaper} with $d=1$
    and $\delta =0$. Indeed, the same argument holds by replacing the condition
    that $\im(\beta) \subset G_j$ in \cite[Lemma 6.9]{torsion-freepaper} with
    the condition that $\im(\oplus_i \phi_i) \subset \bigoplus_{i=1}^p
    (\sigma^+_i)^*(G_j)$, and replacing $\widehat{\mu}$ in \cite[Lemma
    6.9]{torsion-freepaper} with $\mu$. We note that case (1) in the statement
    of \cite[Lemma 6.9]{torsion-freepaper} cannot happen if $d=1$ because of the
    condition $\nu(f')>0$.
\end{proof}

\begin{proof}[Proof of \Cref{prop: hn boundedness}]
    We may replace $\Pur^{\pol}_{g,n,p} \to \prestablepol_{g,n,p} \times
    \SVect^p$ with the morphism $\Pur^{\ast}_{g,n,p} \to \prestablepol_{g,n,p}$
    and the pulled back numerical invariant $\nu$. Since $\Pur^{\ast}_{g,n,p}
    \to \prestablepol_{g,n,p}$ is affine and of finite type (\Cref{prop:
    affineness of pur and valuative criterion for properness}) this follows from
    a similar argument as in the proof of \cite[Prop. 6.10]{torsion-freepaper},
    where we replace the use of \cite[Lemma 6.9]{torsion-freepaper} with
    \Cref{lemma: technical optimization lemma hn boundednes marked sheaves}.
\end{proof}

\subsection{\texorpdfstring{$\Theta$}{Theta}-stratification and moduli spaces}

\begin{notn}[Components with fixed rank and degree] \label{notn: fixed degree
    and rank for pure sheaves} Given two integers $r, d$, we denote by
    $\Pur^{\pol, r, d}_{g,n+p,0} \subset \Pur^{\pol, r, d}_{g,n+p,0}$ the open
    and closed substack where the pure sheaf has rank $r$ and degree $d$ as in
    \Cref{defn: degree and slope of sheaves on prestable curves}. We denote by
    $\Pur^{\pol, r, d}_{g,n,p} \subset \Pur^{\pol}_{g,n,p}$ the corresponding
    open and closed preimage under $\Psi:\Pur^{\pol}_{g,n,p}\to
    \Pur^{\pol}_{g,n+p,0}$.
\end{notn}

\begin{thm} \label{thm: weak theta stratification pure marked sheaves} Fix some
    integers $r,d$. The $\mathbb{R}$-valued numerical invariant $\nu$ from
    \Cref{defn: numerical invariants in marked pure sheaves and Gbunpol} defines
    a weak $\Theta$-stratification on the stack $\Pur^{\pol,r,d}_{g,n,p}$
    relative to $\prestable^{\pol}_{g,n,p}\times \SVect^p$ such that the
    following hold:
    \begin{enumerate}
        \item For each $\gamma \in \mathbb{R}$, the corresponding open stratum
        $\left(\Pur^{\pol, r, d}_{g,n,p}\right)_{\nu \leq \gamma}$ is of finite
        type over $\prestable^{\pol}_{g,n,p}\times \SVect^p$. In particular the
        relative weak $\Theta$-stratification is well-ordered.

        \item If the Noetherian base scheme $S$ is a $\mathbb{Q}$-scheme, then
        $\nu$ defines a well-ordered $\Theta$-stratification, and the center of
        every stratum admits a proper relative good moduli space over
        $\prestable^{\pol}_{g,n,p}\times \SVect^p$.
    \end{enumerate}
\end{thm}
\begin{proof}
    The existence of the weak $\Theta$-stratification follows from \cite[Thm.
    B]{halpernleistner2018structure}, since we have checked the necessary
    requirements on $\Theta$-monotonicity (\Cref{prop: monotonicity for nu on
    marked sheaves}) and HN boundedness (\Cref{prop: hn boundedness}). The rest
    of the conditions can be checked locally on $\prestablepol_{g,n,p} \times
    \SVect^p$, and hence we may fix a $p$-tuple of ranks $N_{\bullet}$ and
    instead consider the morphism $\Pur^*_{g,n,p} \to \prestable^{\pol}_{g,n,p}$
    as in the proof of HN boundedness. We denote by $\Pur^{\ast, r,d}_{g,n,p}
    \subset \Pur^*_{g,n,p}$ the corresponding open and closed substack where the
    pure sheaf has rank $r$ and degree $d$.

    For Part (1), recall that $\left(\Pur^{\ast,r,d}_{g,n,p}\right)_{\nu \leq
    \gamma}$ denotes the open substack of $\Pur^{\ast, r,d}_{g,n,p}$ consisting
    of points whose maximal destabilizing filtration has numerical invariant at
    most $\gamma$. Since the morphism $\Pur^{\ast, r, d}_{g,n,p} \to
    \Pur_{g,n+p,0}^{\pol, r,d}$ is quasicompact by \Cref{prop: affineness of pur
    and valuative criterion for properness}, it suffices to show that the image
    of $\left(\Pur^{\ast,r,d}_{g,n,p}\right)_{\nu \leq \gamma}$ in
    $\Pur_{g,n+p,0}^{\pol, r,d}$ is quasicompact over $\prestablepol_{g,n+p,0}$.
    Given a geometric point $x=(C, \sigma_{\bullet}^{\pm}, \cH,K, W_{\bullet},
    \phi_{\bullet})$ in $\left(\Pur^{\ast,r,d}_{g,n,p}\right)_{\nu \leq
    \gamma}$, let $\mu_{max}(K)$ denote the maximal destabilizing slope (as in
    \Cref{lemma: technical optimization lemma hn boundednes marked sheaves}). In
    view the main result in \cite{langer-boundedness}, it is sufficient to find
    a uniform bound for the difference $\mu_{max}(K) - \mu(K)$ for any geometric
    point $x=(C, \sigma_{\bullet}^{\pm}, \cH,K, W_{\bullet}, \phi)$. If $U
    \subset K$ denotes the largest subsheaf of $K$ with maximally destabilizing
    slope $\mu_{max}(K)$, then the $\mathbb{Z}$-weighted filtration $(K_w)_{w
    \in \mathbb{Z}}$ defined by
    \[ K_w = \begin{cases} 0 & \text{ if $w > 1$ }   \\
            U & \text{ if $w =1$}     \\
            K & \text{ if $w \leq 0$} \end{cases}\] satisfies the conditions in
        \Cref{L:describe_filtration}, and hence it corresponds to a filtration
        $f$ of $x \in \Pur^{\ast,r,d}_{g,n,p}$ relative to
        $\prestable^{\pol}_{g,n,p}$. Using \Cref{eqn: numerical invariant for
        pure sheaves}, we get
    \[ \mu(f) = \sqrt{\rk(U)} \cdot( \mu(U) - \mu(K)) = \sqrt{\rk(U)}
    \cdot(\mu_{max}(K) - \mu(K)) \geq \mu_{max}(K) - \mu(K),\] where the last
    equality follows from the fact that $0<\rk(U)<r$ is a positive integer. If
    $x$ is a point in $\left(\Pur^{\ast,r,d}_{g,n,p}\right)_{\nu \leq \gamma}$,
    then we must have $\gamma \geq \mu(f) \geq \mu_{max}(K) - \mu(K)$, thus
    obtaining the desired uniform bound. This concludes the proof of Part (1).

    We are left to show Part (2). It follows from \cite[Cor.
    2.1.9]{halpernleistner2018structure} that the weak $\Theta$-stratification
    is automatically a $\Theta$-stratification under the assumption that $S$ is
    a scheme over $\mathbb{Q}$. By the semistable reduction theorem in
    \cite[Thm. 6.5]{alper2019existence} and the fact that
    $\Pur_{g,n,p}^{\pol,r,d} \to \prestable^{\pol}_{g,n,p} \times \SVect^p$
    satisfies the existence part of the valuative criterion for properness
    \Cref{prop: affineness of pur and valuative criterion for properness}, it
    follows that every open stratum $\left(\Pur^{\pol,r,d}_{g,n,p}\right)_{\nu
    \leq \gamma}$ also satisfies the existence part of the valuative criterion
    over $\prestable^{\pol}_{g,n,p} \times \SVect^p$. In particular, the same
    holds for the closed union of $\Theta$-strata
    $\left(\Pur^{\pol,r,d}_{g,n,p}\right)_{\nu = \gamma} \hookrightarrow
    \left(\Pur^{\pol,r,d}_{g,n,p}\right)_{\nu \leq \gamma}$. Let $\cZ_{\gamma}$
    be the union of the centers of the strata
    $\left(\Pur^{\pol,r,d}_{g,n,p}\right)_{\nu = \gamma}$. Then, there is a
    surjective morphism $\gr: \left(\Pur^{\pol,r,d}_{g,n,p}\right)_{\nu =
    \gamma} \to \cZ_{\gamma}$, which shows that $\cZ_{\gamma}$ also satisfies
    the existence part of the valuative criterion for properness over
    $\prestable^{\pol}_{g,n,p} \times \SVect^p$. In view of the quasicompactness
    proven in Part (1), the surjectivity of $\gr$ also shows that $\cZ_{\gamma}$
    is of finite type over $\prestable^{\pol}_{g,n,p} \times \SVect^p$. By
    \cite[Thm. A]{alper2019existence}, to conclude the proof of Part (2) it
    suffices to show that $\cZ_{\gamma}$ is $\Theta$-reductive and $S$-complete
    relative to $\prestable^{\pol}_{g,n,p} \times \SVect^p$. To check this, we
    may work instead with the $\Theta$-stratification defined on the morphism
    $\Pur^{\ast, r,d}_{g,n,p} \to \prestable^{\pol}_{g,n,p}$. We denote by
    $\cZ_{\gamma}^*$ the union of centers of strata in $\Pur^{\ast,
    r,d}_{g,n,p}$ corresponding to $\cZ_{\gamma}$.

    By definition, $\cZ_{\gamma}^*$ is an open substack of the stack of graded
    points
    \[\Gr_{\prestable^{\pol}_{g,n,p}}(\Pur^*_{g,n,p}) :=
    \Map_{\prestable^{\pol}_{g,n,p}}(B \mathbb{G}_m \times
    \prestable^{\pol}_{g,n,p}, \Pur^*_{g,n,p}).\] Similarly as in
    \Cref{L:describe_filtration}, for any $S$-scheme $T$ the $T$-points of
    $\Gr_{\prestable^{\pol}_{g,n,p}}(\Pur^*_{g,n,p})$ correspond to tuples $(C,
    \sigma_{\bullet}^{\pm}, \bigoplus_{w \in \mathbb{Z}}K_w, \phi_{\bullet})$,
    where $(C, \sigma_{\bullet}^{\pm})$ is a $T$-family of $(n,p)$-marked
    curves, $\bigoplus_{w \in \mathbb{Z}} K_w$ is a $\mathbb{Z}$-graded
    $T$-family of pure sheaves on $C$, and $\phi_{\bullet}$ is a tuple of
    homomorphisms $\phi_i: \cO^{\oplus N_i}_T \to
    (\sigma_i^+)^*(K_0(\sigma_i^+))$ (note that here the target involves the
    zeroth degree piece $K_0$ of the graded sheaf). For any fixed choice of a
    finite set of integer weights $w_1 < w_2 < \ldots< w_{\ell}$ and an
    $\ell$-tuple of pairs of integers $(r_i, d_i)$, there is an open and closed
    substack $\Gr_{w_{\bullet}}^{r_{\bullet}, d_{\bullet}} \subset
    \Gr_{\prestable^{\pol}_{g,n,p}}(\Pur^*_{g,n,p})$ where we require that
    $K_w=0$ for $w \neq w_1, \ldots, w_\ell$, and furthermore $\rk(K_{w_i}) =
    r_i$ and $\deg(K_{w_i}) = d_i$ for all $1 \leq i \leq \ell$. Note that
    $\Gr_{\prestable^{\pol}_{g,n,p}}(\Pur^*_{g,n,p})$ is a union of the open and
    closed substacks $\Gr_{w_{\bullet}}^{r_{\bullet}, d_{\bullet}}$ as we range
    over all tuples $(w_{\bullet}, r_{\bullet}, d_{\bullet})$, and without loss
    of generality we may always assume that $w_j =0$ for some $1 \leq j \leq
    \ell$ (by possibly setting $r_j=d_j =0$). There is an isomorphism to the
    following multiple fiber product over $\prestable^{\pol}_{g, n+p,0}$:
    \[\Gr_{w_{\bullet}}^{r_{\bullet}, d_{\bullet}} \xrightarrow{\sim}
    \Pur^{\ast, r_j, d_j}_{g,n,p} \times_{\prestable^{\pol}_{g, n+p,0}} \prod_{1
    \leq i \leq \ell, \; i\neq j} \Pur^{\pol, r_i, d_i}_{g,n+p,0} \] where we
    send a $T$-point $(C, \sigma_{\bullet}^{\pm}, \bigoplus_{w \in
    \mathbb{Z}}K_w, \phi_{\bullet})$ to the element $(C, \sigma_{\bullet}^{\pm},
    K_0, \phi_{\bullet})$ of $\Pur^{\ast, r_j, d_j}_{g,n,p}$ along with the
    tuple $(C, \sigma_{\bullet}, K_{w_i})$ of elements of $\Pur^{\pol, r_i,
    d_i}_{g,n+p,0}$ for $1 \leq i \leq \ell$ with $i \neq j$. It suffices to
    show that for any given tuple $(w_{\bullet}, r_{\bullet}, d_{\bullet})$, the
    open and closed substack $\cZ_{\gamma}^* \cap
    \Gr_{w_{\bullet}}^{r_{\bullet}, d_{\bullet}} \subset \cZ_{\gamma}^*$ is
    $\Theta$-reductive and $S$-complete. By definition $\cZ_{\gamma}^* \cap
    \Gr_{w_{\bullet}}^{r_{\bullet}, d_{\bullet}}$ is the open substack of points
    in $\Gr_{w_{\bullet}}^{r_{\bullet}, d_{\bullet}}$ that are isomorphic to the
    associated graded point of the Harder-Narasimhan filtration of a point in
    $\Pur^{\ast, r,d}_{g,n,p}$. In particular, such filtration must be maximally
    destabilizing for the numerical invariant $\nu$, and hence it must arise as
    the associated graded of a filtration satisfying the properties from
    \Cref{lemma: technical optimization lemma hn boundednes marked sheaves}.
    This means that there exists a fixed positive rational number $Q>0$ such
    that we have $w_i = Q \cdot \left(\frac{2d_i}{r_i} - \frac{2d}{r}\right)$
    for $i \neq j$, and that we are forced to have $0  \geq \frac{2d_j}{r_j} -
    \frac{2d}{r}$ whenever the rank $r_j$ of the $0^{th}$ weight piece is
    nonzero.

    The description of the open
    \[\cZ_{\gamma}^* \cap \Gr_{w_{\bullet}}^{r_{\bullet}, d_{\bullet}} \subset
    \Gr_{w_{\bullet}}^{r_{\bullet}, d_{\bullet}} \cong \Pur^{\ast, r_j,
    d_j}_{g,n,p} \times_{\prestable^{\pol}_{g, n+p,0}} \prod_{1 \leq i \leq
    \ell, \; i\neq j} \Pur^{\pol, r_i, d_i}_{g,n+p,0}\] as the substack of
    graded semistable points with respect to $\nu$ \cite[Thm.
    4.2.3]{halpernleistner2018structure} in this case identifies $\cZ_{\gamma}^*
    \cap \Gr_{w_{\bullet}}^{r_{\bullet}, d_{\bullet}}$ as a product of substacks
    $\cU_j \times_{\prestable^{\pol}_{g, n+p,0}} \prod_{1 \leq i \leq \ell, \;
    i\neq j} \cU_i$ where
    \begin{itemize}
        \item $\cU_i \subset \Pur^{\pol, r_i, d_i}_{g,n+p,0}$ is the
        $\Theta$-semistable locus with respect to the line bundle
        $\cL_{\alpha_i}$ from \Cref{defn: rational line bundle L mu} with
        $\alpha_i := \frac{w_i}{Q} + \frac{2d}{r} - \frac{2d_i}{r_i}$, with $Q$
        as above, and similarly
        \item $\cU_j \subset \Pur^{\ast, r_j, d_j}_{g,n,p}$ is the
        $\Theta$-semistable locus with respect to the line bundle
        $\cL_{\alpha_j}$ from \Cref{notn: line bundles on Pur} with $\alpha_j :=
        \frac{w_j}{Q} + \frac{2d}{r} - \frac{2d_i}{r_i} = \frac{2d}{r} -
        \frac{2d_i}{r_i}$.
    \end{itemize}
    (For complete precision, we may use the convention that $\alpha_i=0$
    whenever $r_i =0$; in that case the stack $\Pur^{ r_i,d_i}$ that we are
    dealing with is trivially isomorphic to the base stack $\prestable^{\pol}$).
    In view of the description of $w_i$ above, we conclude that $\alpha_i = 0$
    if $i \neq j$, and we always have $\alpha_j \geq 0$. By \Cref{prop:
    monotonicity for nu on marked sheaves}, the line bundle $\cL_{\alpha_i}$ for
    $1 \leq i \leq \ell$ is strictly $\Theta^2$-monotone, and hence the open
    substack $\cU_i$ defined as the semistable locus of $\cL_{\alpha_i}$ is
    $\Theta$-reductive and $S$-complete for all $1 \leq i \leq \ell$ by
    \cite[Rem.~5.5.9]{halpernleistner2018structure}. Hence, the same holds for
    the product $\cZ_{\gamma}^* \cap \Gr_{w_{\bullet}}^{r_{\bullet},
    d_{\bullet}} = \cU_j \times_{\prestable^{\pol}_{g, n+p,0}} \prod_{1 \leq i
    \leq \ell, \; i\neq j} \cU_i$.
\end{proof}

\section{The stack of Gieseker bundles} \label{section: moduli of gieseker
bundles} In this section we prove our main structural results on the existence
of $\Theta$-stratifications and moduli spaces for the moduli of pure Gieseker
bundles.

\subsection{Log-crepant contractions of prestable curves}

\begin{defn}[Degree 1 morphisms of curves] \label{defn: degree 1} Let $\varphi:
    \tilde{C} \to C$ be a morphism of prestable curves defined over a field $k$.
    We say that $\varphi$ has degree $1$ if for any choice of algebraically
    closed overfield $F \supset k$, the pushforward of the fundamental class
    $[\tilde{C}_{F}]$ equals the fundamental class $[C_{F}]$ as elements of the
    numerical cone of curves in $C_{F}$ (i.e. for all line bundles $\cL$ on
    $C_F$ we have $\deg(\pi^*(\cL)) = \deg(\cL)$).
\end{defn}

An equivalent way of phrasing \Cref{defn: degree 1} is that for any
algebraically closed overfield $F \supset k$ there is a dense open $U \subset
C_F$ such that the induced morphism $\varphi^{-1}(U) \to U$ is an isomorphism.
Note that it suffices to check the property in \Cref{defn: degree 1} for a
single choice of $F \supset k$, and then it will hold for all algebraically
closed field extensions.

\begin{lemma} \label{lemma: degree is open and closed} Let $T$ be an $S$-scheme,
    and let $\varphi: \tilde{C} \to C$ be a morphism between two $T$-families of
    genus $g$ prestable curves. Then, the set of points $t \in T$ such that the
    restriction $\varphi_{t}$ has degree $1$ is open and closed in $T$.
\end{lemma}
\begin{proof}
    Let $T_{1}$ denote the set of points $t \in T$ such that $\varphi_t$ has
    degree $1$, and set $T_{\neq 1} := T \setminus T_{1}$ to be its complement.
    We want to show that both $T_{1}$ and $T_{\neq 1}$ are open. This can be
    checked \'etale locally on $T$. For any geometric point $\overline{t} \to
    T$, the numerical cone of curves of $C_{\overline{t}}$ is isomorphic to
    $\mathbb{R}_{\geq 0}^{|I|}$, where $I$ is the set of irreducible components
    of $C_{\overline{t}}$. By \cite[Thm. 2.5.2(i)]{romagny-composantes}, we may
    pass to an \'etale cover of $T$ and assume that for every geometric point
    $\overline{t} \to T$ with image $t \in T$, the number of irreducible
    components of $C_t$ and $C_{\overline{t}}$ agree. Hence, we may check
    whether a point $t \in T$ is degree $1$ without passing to a geometric
    over-point $\overline{t}$. For any $t \in T$, the image of the Picard group
    $\Pic(C_t)$ in the real N\'eron-Severi space $\mathbb{R}^{|I|}$ spans it as
    an $\mathbb{R}$-vector space. If we choose a finite set
    $\{\cL_i\}_{i=1}^{|I|}\subset \Pic(C_t)$ whose image spans
    $\mathbb{R}^{|I|}$, then the morphism $\varphi_t$ has degree 1 if and only
    if $\deg(\cL_i) = \deg(\varphi_{\overline{t}}^*(\cL_i))$ for all $1 \leq i
    \leq |I|$.

    Choose a point $t_0 \in T$ and a spanning set $\{\cL_i\}_{i=1}^{|I|}\subset
    \Pic(C_{t_0})$. Since the relative Picard stack for the family of curves $C
    \to T$ is smooth, there is an \'etale neighborhood of $U \to T$ of $t_0$
    such that every $\cL_i$ extends to a line bundle $\tilde{\cL}_i$ on $C_U$.
    After shrinking $U$, we may assume that the dual graph of $C_t$ for any $t
    \in U$ is a generalization of the dual graph of $C_{t_0}$ as described in
    \cite[Lemma 2.5]{frenkel-teleman-tolland-ggw}. The local constancy of
    degrees of line bundles implies that the degree of the restriction of
    $\tilde{\cL}_i$ to any irreducible component of $C_t$ is always a sum of a
    subset of the degrees of $\cL_i$ on irreducible components of $C_{t_0}$ (cf.
    \cite[Lemma 2.5]{frenkel-teleman-tolland-ggw}). In other words, if we denote
    by $I'$ the irreducible components of $C_t$, then the image of
    $\widetilde{\cL}_i|_{C_t}$ in $\mathbb{R}^{|I'|}$ is obtained by taking the
    image of $\cL_i$ under a projection $\mathbb{R}^{|I|} \twoheadrightarrow
    \mathbb{R}^{|I'|}$. In particular, the images of $\widetilde{\cL}_i|_{C_t}$
    still span $\mathbb{R}^{|I'|}$. Hence, for all $t \in U$, we have that
    $\varphi_{t}$ is of degree $1$ if and only if we have
    $\deg(\tilde{\cL}_i|_{C_t}) = \deg(\varphi_y^*(\cL_i|_{C_t}))$ for every
    $i$. By local constancy of degrees of line bundles, it follows that both
    $U_{1}$ and $U_{\neq 1}$ are open.
\end{proof}

\begin{defn}[Contractions of marked curves] \label{defn: contraction} Let $T$ be
    an $S$-scheme, and let $(\widetilde{C}, \sigma_{\bullet}^{\pm})$ be a
    $T$-family of $(n,p)$-marked genus $g$ prestable curves. Let $C \to T$ be a
    family of genus $g$ prestable curves. A morphism $\varphi: \widetilde{C} \to
    C$ is called a contraction of marked curves if both of the following are
    satisfied:
    \begin{enumerate}[(1)]
        \item The tuple $(C, \varphi\circ \sigma_{\bullet}^{\pm})$ is a
        $(n,p)$-marked family of genus $g$ prestable curves.
        \item The morphism $\varphi: \tilde{C} \to C$ has degree $1$ on every
        fiber of $T$.
    \end{enumerate}
    We call $\varphi$ log-crepant if the induced morphism of relative
    log-canonical line bundles
    \[\varphi^*\omega_{C/T}\left(\sum_{i=1}^p \varphi \circ \sigma_i^+ +
    \sum_{i=1}^n \varphi \circ \sigma_i^-\right) \to \omega_{\tilde{C}/T}\left(
    \sum_{i=1}^p\sigma_i^+ + \sum_{i=1}^n \sigma_i^-\right)\] is an isomorphism.
    In this case, we say that $(\tilde{C}, \sigma_{\bullet}^{\pm})$ is a
    log-semistable model of the $(n,p)$-marked curve $(C, \varphi \circ
    \sigma_{\bullet}^{\pm})$.
\end{defn}

\begin{defn}[Stack of log-crepant contractions]
    We denote by $\precont_{g,n,p}$ the stack whose $T$-points consists of
    tuples $(\varphi: \tilde{C} \to C, \sigma_{\bullet}^{\pm})$, where
    $(\tilde{C}, \sigma_{\bullet}^{\pm})$ is a $T$-family of $(n,p)$-marked
    prestable curves of genus $g$, and $\varphi: \tilde{C} \to C$ is a
    log-crepant contraction.
\end{defn}

\begin{lemma} \label{lemma: properties of contraction stack} The following hold:
    \begin{enumerate}
        \item $\precont_{g,n,p}$ is a smooth algebraic stack over $S$.
        \item The forgetful morphism $f_{targ}: \precont_{g,n,p} \to
        \prestable_{g,n,p}$ given by $(\varphi: \tilde{C} \to C, \sigma_i^{\pm})
        \mapsto (C, \varphi \circ \sigma_i^{\pm})$ is flat.
        \item The forgetful morphism $f_{src}: \precont_{g,n,p} \to
        \prestable_{g,n,p}$ given by $(\varphi: \tilde{C} \to C, \sigma_i^{\pm})
        \mapsto (\tilde{C}, \sigma_i^{\pm})$ is smooth.
    \end{enumerate}
\end{lemma}
\begin{proof}

    Everything may be checked after smooth base-change on the target
    $\prestable_{g,n,p}$ of the morphism $f_{targ}: \precont_{g,n,p} \to
    \prestable_{g,n,p}$. Base-changing via a smooth forgetful morphism
    $\prestable_{g,n,p+i} \to \prestable_{g,n,p}$ for $i\gg0$, we may locally
    add enough markings to the target curve so that it becomes stable, and hence
    it is identified with the stabilization of the source marked curve. In
    particular, the data of the target curve becomes redundant. Hence,
    $\precont_{g,n,p} \times_{ \prestable_{g,n,p}} \prestable_{g,n,p+i}$ becomes
    locally isomorphic to the smooth stack of marked prestable curves
    $\prestable_{g,n,p+i}$, the base-changed source morphism $f_{src}$ is
    locally identified with the smooth forgetful morphism $\prestable_{g,n,p+i}
    \to \prestable_{g,n,p}$, and the base-changed target morphism $f_{targ}$ is
    locally identified with the stabilization morphism $st: \prestable_{g,n,p+i}
    \to \overline{\cM}_{g,n,p+i}$, which is known to be flat \cite[Prop
    3]{behrend-gromov-witten-theory}. \endnote{ For every integer $i \geq 0$,
    there is a forgetful morphism $\prestable_{g,n,p+i} \to \prestable_{g,n,p}$
    which is smooth and representable. We denote by $\overline{\cM}_{g,n,p+i}
    \subset \prestable_{g,n,p+i}$ the open substack where the marked prestable
    curve is stable as a $(n+p+i)$-pointed curve. Consider the fiber product
        \[
            \begin{tikzcd}[ampersand replacement = \&]
                \cF_{g,n,p+i} \ar[r] \ar[d] \&  \precont_{g,n,p} \ar[d, "f_{targ}"] \\  \prestable_{g,n,p+i} \ar[r] \& \prestable_{g,n,p}.
            \end{tikzcd}
        \]
        The stack $\cF_{g,n,p+i}$ parametrizes an $(n,p)$-marked curve
        $(\tilde{C}, \sigma_{\bullet}^{\pm})$, an $(n, p+i)$-marked curve $(C,
        \overline{\sigma}_{\bullet}^{\pm})$, and a log-crepant contraction
        $\varphi: \tilde{C} \to C$ that is compatible with the first $n+p$
        markings. We denote by $\cU_{g,n,p+i} \subset \cF_{g,n,p+i}$ the open
        substack consisting of families where
        $(C,\overline{\sigma}_{\bullet}^{\pm})$ is stable as a $(n+p+i)$-pointed
        curve, and where the markings $\overline{\sigma}_{p+1}^+, \ldots
        \overline{\sigma}_{p+i}^+$ on $C$ land on the open locus of $C$ where
        $\varphi: \tilde{C} \to C$ is an isomorphism. Hence, the markings
        $\overline{\sigma}_{p+1}^+, \ldots \overline{\sigma}_{p+i}^+$ lift
        uniquely to $\tilde{C}$, and it follows that $\cU_{g,n,p+i} \subset
        \precont_{g,n,p+i}$ is the open substack of log-crepant contractions
        $\varphi$ of marked $(n,p+i)$-prestable genus $g$ curves such that the
        target is stable and $\varphi$ is an isomorphism over the images of the
        last $i$ markings.

        The smooth representable morphism $\bigsqcup_{i\geq 0} \cU_{g,n,p+i} \to
        \precont_{g,n,p}$ is surjective, and fits into a commutative diagram
        \[
            \begin{tikzcd}[ampersand replacement = \&]
                \bigsqcup_{i\geq 0} \cU_{g,n,p+i} \ar[r] \ar[d] \&  \precont_{g,n,p} \ar[d,"f_{targ}"] \\  \bigsqcup_{i\geq 0} \overline{\cM}_{g,n,p+i} \ar[r] \& \prestable_{g,n,p}
            \end{tikzcd}
        \]
        where $\bigsqcup_{i\geq 0} \overline{\cM}_{g,n,p+i} \to
        \prestable_{g,n,p}$ is also representable and smooth. To conclude, it
        suffices to show that each $\cU_{g,n,p+i}$ is smooth, that the morphism
        $\cU_{g,n,p+i} \to \overline{\cM}_{g,n,p+i}$ is flat, and that the
        composition $\cU_{g,n,p+i} \to \precont_{g,n,p} \xrightarrow{f_{src}}
        \prestable_{g,n,p}$ is smooth. To show this, we consider the smooth
        stack $\prestable_{g,n,p+i}$, along with its stabilization morphism $st:
        \prestable_{g,n,p+i} \to \overline{\cM}_{g,n,p+i}$, which is flat by
        \cite[Prop 3]{behrend-gromov-witten-theory}. By construction, the
        composition $\cU_{g,n,p+i} \hookrightarrow \precont_{g,n,p+i}
        \xrightarrow{f_{src}} \prestable_{g,n,p+i}$ exhibits $\cU_{g,n,p+i}$ as
        the open substack of $\prestable_{g,n,p+i}$ consisting of
        $(n,p+i)$-marked curves $(\tilde{C}, \sigma_{\bullet}^{\pm})$ such that
        their marked stabilization morphism $\varphi: (\tilde{C},
        \sigma_{\bullet}^{\pm}) \to (C, \overline{\sigma}_{\bullet})$ induces an
        identification of log-canonical line bundles. In particular
        $\cU_{g,n,p+i} \subset \prestable_{g,n,p+i}$ is smooth over $S$. The
        composition $\cU_{g,n,p+i} \hookrightarrow \prestable_{g,n,p+i}
        \xrightarrow{st} \overline{\cM}_{g,n,p+i}$ is flat because of the
        flatness of the stabilization morphism $st$. On the other hand, the
        composition $\cU_{g,n,p+i} \to \precont_{g,n,p} \xrightarrow{f_{src}}
        \prestable_{g,n,p}$ factors as $\cU_{g,n,p+i} \hookrightarrow
        \prestable_{g,n,p+i} \to \prestable_{g,n,p}$, where
        $\prestable_{g,n,p+i} \to \prestable_{g,n,p}$ is the smooth forgetful
        morphism. Therefore, $\cU_{g,n,p+i} \to \precont_{g,n,p}
        \xrightarrow{f_{src}} \prestable_{g,n,p}$ is smooth.}
\end{proof}

\subsection{Gieseker bundles}

\begin{defn}[Stack of Gieseker bundles] \label{D:GBun} Let $\CBun_{g,n,p}(N)$
    denote the pseudofunctor that sends an $S$-scheme $T$ to the groupoid of
    tuples $(\varphi : \tilde{C} \to C, \sigma^\pm_{\bullet}, E)$, where
    \begin{enumerate}[(1)]
        \item $(\varphi : \tilde{C} \to C, \sigma^\pm_{\bullet})$ is a
        $T$-family of log-crepant contractions in $\precont_{g,n,p}(T)$.

        \item $E$ is a rank $N$ vector bundle on $\tilde{C}$ such that $\det(E)$
        is $\varphi$-ample and for every $t \in T$ the counit $\varphi_t^*
        (\varphi_t)_*(E|_{\tilde{C}_t}) \to E|_{\tilde{C}_t}$ is surjective.
    \end{enumerate}
\end{defn}

\begin{lemma} \label{lemma: pushforward lemma} Let $T$ be an $S$-scheme, and let
    $(\varphi : \tilde{C} \to C, \sigma^\pm_{\bullet})$ be a (not necessarily
    log-crepant) contraction of $T$-families of $(n,p)$-marked prestable curves
    of genus $g$ (as in \Cref{defn: contraction}). Let $E$ be a vector bundle on
    $\tilde{C}$ such that for every $t \in T$ the counit $\varphi_t^*
    (\varphi_t)_*(E|_{\tilde{C}_t}) \to E|_{\tilde{C}_t}$ is surjective. Then,
    the following hold:
    \begin{enumerate}[(a)]
        \item $R^1 \varphi_*\left(E(-\sum_{i=1}^p \sigma_i^{+})\right)= R^1
        \varphi_*(E)=0$, and this holds after base-change by any morphism $T'
        \to T$.
        \item Both $\varphi_*(E)$ and $\varphi_*\left(E(-\sum_{i=1}^p
        \sigma_i^{+})\right)$ are $T$-flat, and their formations commute with
        arbitrary base-change $T' \to T$.
        \item The formation of the counits $\varphi^*\varphi_*\left(E)\right)
        \to E$ and $\varphi^*\varphi_*\left(E(-\sum_{i=1}^p \sigma_i^{+})\right)
        \to E(-\sum_{i=1}^{p} \sigma_i^+)$ commutes with arbitrary base-change
        $T' \to T$.
    \end{enumerate}
\end{lemma}
\begin{proof}
    We can work Zariski locally on $C$, and hence we can assume without loss of
    generality that $C$ and $T=\Spec(R)$ are affine schemes. Then
    $R\varphi_*\left(E(-\sum_{i=1}^{p} \sigma_i^+)\right)$ is quasi-isomorphic
    to the Cech complex $\check{\cC}^{\bullet}\left(E(-\sum_{i=1}^{p}
    \sigma_i^+)\right)$ associated to a finite cover of $C$ by affines, which is
    concentrated in degree $\geq 0$ and consists of $R$-flat $\cO_{C}$-modules.
    Furthermore, the formation of the Cech complex is compatible with arbitrary
    base-change $T' = \Spec(R') \to \Spec(R) = T$. The proofs of (a), (b) and
    (c) follow for $E(-\sum_{i=1}^p \sigma_i^+)$ if we can show that the complex
    $\check{\cC}^{\bullet}\left(E(-\sum_{i=1}^{p} \sigma_i^+)\right)$ has Tor
    amplitude $[0,0]$. We may check this fiberwise on $T$, and hence we may
    assume that $T$ is the spectrum of an algebraically closed field. By working
    locally on $C$, we may choose a surjection from a free
    $\mathcal{O}_{C}$-module $\mathcal{G}$ onto $\varphi_*(E)$. Thus, we get a
    surjection $\varphi^*(\mathcal{G})(-\sum_{i=1}^{p} \sigma_i^+) \to
    E(-\sum_{i=1}^{p} \sigma_i^+)$ fitting into a short exact sequence
    \[ 0 \to \mathcal{K} \to \varphi^*(\mathcal{G})(-\sum_{i=1}^{p} \sigma_i^+)
    \to E(-\sum_{i=1}^{p} \sigma_i^+) \to 0\] Note that $\varphi: \tilde{C} \to
    C$ exactly contracts a disjoint union $\bigsqcup_i \tilde{C}_i \subset
    \tilde{C}$ of prestable subcurves $\tilde{C}_i$ consisting of rational
    components and containing no cycles of curves, where each $\tilde{C}_i$
    contains at most one marking, and if $\varphi(\tilde{C}_i)$ is a node of $C$
    then $\tilde{C}_i$ contains no markings. A local computation shows that
    $R^1\varphi_*(\mathcal{O}_{\tilde{C}}(-\sum_{i=1}^{p} \sigma_i^+)$ vanishes
    (cf. \cite[\href{https://stacks.math.columbia.edu/tag/0E7Q}{Tag
    0E7Q}]{stacks-project} and references therein), and in view of the short
    exact sequence above we conclude that $R^i \varphi_*(E(-\sum_{i=1}^{p}
    \sigma_i^+)) =0$ for $i>0$, as desired.

    The same argument applies if we replace $E(-\sum_{i=1}^{p} \sigma_i^+)$ with
    $E$.\qedhere
\end{proof}

The compatibility of $\varphi_*\left(E(-\sum_{i=1}^p \sigma_i^{+})\right)$ with
base-change in \Cref{lemma: pushforward lemma} makes the following meaningful.
\begin{defn}[Stack of pure Gieseker bundles] \label{defn: stack of pure gieseker
    bundles} The stack $\GBun_{g,n,p}(N)$ is the substack of $\CBun_{g,n,p}(N)$
    of $T$-points $(\varphi : \tilde{C} \to C, \sigma^\pm_{\bullet}, E)$ such
    that the restriction of $\varphi_\ast(E(-\sum_{i=1}^{p} \sigma^+_i))$ to
    every $T$-fiber of $C \to T$ is pure of dimension $1$.
\end{defn}

\begin{lemma}
    \label{L:cBun algebraic}
    The following statements hold.
    \begin{enumerate}[(a)]
        \item $\CBun_{g,n,p}(N)$ is an algebraic stack, and the following
        forgetful morphism is smooth
              \[\CBun_{g,n,p}(N) \to \precont_{g,n,p}, \; \; \; (\varphi :
              \tilde{C} \to C, \sigma^\pm_{\bullet}, E) \mapsto (\varphi :
              \tilde{C} \to C, \sigma^\pm_{\bullet})\] In particular,
              $\CBun_{g,n,p}(N)$ is smooth over $S$.
        \item $\GBun_{g,n,p}(N) \subset \CBun_{g,n,p}(N)$ is an open substack.
        \item The following forgetful morphism is flat with affine relative
        diagonal
              \[c: \CBun_{g,n,p}(N) \to \prestable_{g,n,p}, \; \; \; (\varphi:
              \tilde{C} \to C, \sigma^\pm_\bullet, E) \mapsto (C,\varphi \circ
              \sigma^\pm_\bullet).\]
    \end{enumerate}
\end{lemma}
\noindent \textit{Proof of (a).} Since the conditions on $E$ in the definition
of Gieseker bundle are open, the stack $\CBun_{g,n,p}(N) \to \precont_{g,n,p}$
is open inside the mapping stack $\Map_{\precont_{g,n,p}}(\tilde{C}, \BGL_N)$
for the universal source curve $\widetilde{C} \to \precont_{g,n,p}$, and hence
it is locally of finite type and with affine relative diagonal by \cite[Thm.
1.2]{hall-rydh-tannakahom}. The smoothness holds because deformations of vector
bundles on curves are unobstructed. \endnote{Let $\tilde{C}_u \to
\precont_{g,n,p}$ denote the universal prestable source curve. By \cite[Thm.
1.2]{hall-rydh-tannakahom}, the relative stack
$\Map_{\precont_{g,n,p}}(\tilde{C}_u, \BGL_N)$ of vector bundles on the fibers
of the flat proper morphism $\tilde{C}_u \to \precont_{g,n,p}$ is an algebraic
stack locally of finite type over $\precont_{g,n,p}$. Then $\CBun_{g,n,p}(N)$ is
the substack of $\Map_{\precont_{g,n,p}}(\tilde{C}_u, \BGL_N)$ parameterizing
vector bundles $E$ such that fiberwise on $T$ we have that $\det(E)$ is
$\varphi$-ample \cite[Cor. 9.6.5]{egaiv} and the counit $\varphi^*\varphi_*(E)
\to E$ is surjective. Both of these conditions are open (the ampleness of
$\det(E)$ by \cite[Cor. 9.6.4]{egaiv}, and the surjectivity of the counit by the
compatibility of the counit with base-change in \Cref{lemma: pushforward
lemma}). We conclude that $\CBun_{g,n,p}(N) \subset
\Map_{\precont_{g,n,p}}(\tilde{C}_u, \BGL_N)$ is also an algebraic stack locally
of finite type over $\precont_{g,n,p}$. The forgetful morphism $\CBun_{g,n,p}(N)
\to \precont_{g,n,p}$ given by $(\varphi : \tilde{C} \to C,
\sigma^\pm_{\bullet}, E) \mapsto (\varphi : \tilde{C} \to C,
\sigma^\pm_{\bullet})$ is smooth, because the deformation theory of vector
bundles on curves is unobstructed. Since $\precont_{g,n,p}$ is smooth over $S$
(\Cref{lemma: properties of contraction stack}), it follows that
$\CBun_{g,n,p}(N)$ is smooth over $S$.}

\medskip

\noindent \textit{Proof of (b).} By \Cref{lemma: pushforward lemma}, there is a
well-defined morphism
\[\CBun_{g,n,p}(N) \to \Coh_{g,n+p,0}(N), \; \; \; (\varphi : \tilde{C} \to C,
\sigma^\pm_{\bullet}, E) \mapsto (C, \varphi \circ \sigma_{\bullet},
\varphi_*(E(-\sum_{i=1}^p \sigma_i^+)).\] By definition, $\GBun_{g,n,p}(N)
\subset \CBun_{g,n,p}(N)$ is the preimage of the open substack
$\Pur_{g,n+p,0}(N) \subset \Coh_{g,n+p,0}(N)$ (see \Cref{lemma: stack of pure
sheaves is algebraic}).

\medskip

\noindent \textit{Proof of (c).} The morphism $c:\CBun_{g,n,p}(N) \to
\prestable_{g,n,p}$ is flat, because it is the composition of the smooth
forgetful morphism $\CBun_{g,n,p}(N) \to \precont_{g,n,p}$ from part (a)
followed by the flat morphism $f_{targ}: \precont_{g,n,p} \to
\prestable_{g,n,p}$ (\Cref{lemma: properties of contraction stack}). We are left
to show that $\CBun_{g,n,p}(N) \to \prestable_{g,n,p}$ has affine relative
diagonal, for which we may consider the base-change $c:\CBun_{g,n,p}^{\pol}(N)
\to \prestable_{g,n,p}^{\pol}$ to the stack of polarized curves. Then, for every
quasicompact open $ \cU \subset \CBun_{g,n,p}^{\pol}(N)$, we may use a high
enough power $\varphi^*(\cH)^{\otimes m}$ of the pullback of polarization on the
target curve twisted by the determinant $\det(E)$ in order to obtain a
polarization of the source curve, and hence obtain a morphism $\cU \to
(\prestable^{\pol}_{g,n,p})^2$ sending a tuple $(\varphi: \widetilde{C} \to C,
\sigma^{\pm}_{\bullet}, E, \cH)$ to the pair of polarized curves
$(\widetilde{C}, \sigma_{\bullet}^{\pm}, \varphi^*(\cH)^{\otimes m} \otimes
\det(E))$ and $(C, \varphi \circ \sigma_{\bullet}^{\pm},\cH)$. Since the second
projection $\left(\prestablepol_{g,n,p}\right)^2 \to \prestablepol_{g,n,p}$ has
affine relative diagonal (\Cref{prop: properties of the stack of polarized
curves}), it suffices to show that $\cU \to (\prestable^{\pol}_{g,n,p})^2$ has
affine relative diagonal. But the fibers of $\cU \to
(\prestable^{\pol}_{g,n,p})^2$ over a pair of curves $(\widetilde{C},
\sigma_{\bullet}^{\pm}, \cH')$ and $(C, \varphi \circ
\sigma_{\bullet}^{\pm},\cH)$ parameterize triples of vector bundles $E$,
log-crepant contractions $\varphi: \widetilde{C} \to C$ and isomorphisms
$\det(E) \xrightarrow{\sim} \varphi^*(\cH)^{-\otimes m} \otimes \cH'$, subject
to some open condition. By
\cite[\href{https://stacks.math.columbia.edu/tag/0DPN}{Tag
0DPN}]{stacks-project} and \Cref{lemma: degree is open and closed}, the stack
parameterizing log-crepant contractions $\varphi$ is separated over
$(\prestable^{\pol}_{g,n,p})^2$. We conclude because the stack of rank
$N$-vector bundles on the universal curve $\widetilde{C} \to
(\prestable^{\pol}_{g,n,p})^2$ from the first coordinate has affine relative
diagonal \cite[Thm. 1.2]{hall-rydh-tannakahom}, and the additional data of the
isomorphism $\det(E) \xrightarrow{\sim} \varphi^*(\cH)^{-\otimes m} \otimes
\cH'$ is parameterized by a relatively affine stack. \endnote{Consider the fiber
product $\CBun_{g,n,p}^{\pol}(N) := \CBun_{g,n,p}(N) \times_{\prestable_{g,n,p}}
\prestable_{g,n,p}^{\pol}$, which parametrizes points $(\varphi: \tilde{C} \to
C, \sigma^\pm_\bullet, E)$ of $\CBun_{g,n,p}(N)$ along with a relative
polarization $\cH$ of $C$. By \Cref{prop: properties of the stack of polarized
curves}, it suffices to prove that $\CBun_{g,n,p}^{\pol}(N) \to
\prestablepol_{g,n,p}$ has relative affine diagonal. It enough to prove this for
every quasicompact open substack $\cU \subset \CBun_{g,n,p}^{\pol}(N)$. By
quasicompactness of $\cU$, there exists a positive integer $m>0$ such that for
every scheme $T$ and every $T$-point $(\varphi: \tilde{C} \to C,
\sigma^\pm_\bullet, E, \cH)$ of $\cU$, the line bundle $\varphi^*(\cH)^{\otimes
m} \otimes \det(E)$ is relatively ample on $\tilde{C} \to T$. This induces a
morphism of stacks $\cU \to \left(\prestablepol_{g,n,p}\right)^2$ given by
\[ (\varphi: \tilde{C} \to C, \sigma^\pm_\bullet, E, \cH) \mapsto
\left[(\tilde{C}, \sigma_{\bullet}^{\pm},\varphi^*(\cH)^{\otimes m} \otimes
\det(E)), (C, \varphi \circ \sigma_{\bullet}^{\pm},\cH)\right].\] By \Cref{prop:
properties of the stack of polarized curves} the second projection
$\left(\prestablepol_{g,n,p}\right)^2 \to \prestablepol_{g,n,p}$ has affine
relative diagonal. Hence, it suffices to show that the morphism $\cU \to
\left(\prestablepol_{g,n,p}\right)^2$ defined above has relative affine
diagonal. By \cite[\href{https://stacks.math.columbia.edu/tag/0DPN}{Tag
0DPN}]{stacks-project}, there is a separated representable morphism locally of
finite type $\mathcal{M}or \to \left(\prestablepol_{g,n,p}\right)^2$, where
$\mathcal{M}or$ is an algebraic stack that parametrizes pairs of polarized
$(n,p)$-marked curves $(\tilde{C}, \sigma_{\bullet}^{\pm},\tilde{\cH}), (C,
\overline{\sigma}_{\bullet}^{\pm},\cH)$ along with a morphism $\varphi:
\tilde{C} \to C$. \Cref{lemma: degree is open and closed} implies there is a
locally closed substack $\mathcal{M}or_{=1} \subset \mathcal{M}or$ where the
morphism $\varphi$ is a log-crepant contraction of marked curves (note that the
condition on the equality of markings is closed, and the requirement on the
isomorphism of log-canonical bundles is an open condition). The morphism $\cU
\to \left(\prestablepol_{g,n,p}\right)^2$ factors as $\cU \to \mathcal{M}or_{=1}
\to \left(\prestablepol_{g,n,p}\right)^2$ by remembering also the morphism
$\varphi$. It suffices to show that $\cU \to \mathcal{M}or_{=1}$ has relative
affine diagonal.

Let $(\tilde{\cC}, \tilde{\cH}) \to \mathcal{M}or_{=1}$ denote the universal
source prestable curve equipped with its universal polarization. By \cite[Thm.
1.2]{hall-rydh-tannakahom}, the mapping stack $\cW :=
\Map_{\mathcal{M}or_{=1}}(\tilde{\cC}, \BGL_N)$ is an algebraic stack with
affine relative diagonal over $\mathcal{M}or_{=1}$. The morphism $h: \cU \to
\mathcal{M}or_{=1}$ factors as $\cU \to \cW \to \mathcal{M}or_{=1}$ by
remembering also the vector bundle $E$. By construction $\cU$ is the equalizer
of two morphisms $f, P: \cW \to \Map_{\mathcal{M}or_{=1}}(\tilde{\cC}, \BGL_1),$
where $P$ is induced by the universal polarization $\tilde{\cH}$, and $f$ sends
a $T$-family $\varphi:(\tilde{C}, \sigma_{\bullet}^{\pm},\tilde{\cH}) \to  (C,
\overline{\sigma}_{\bullet}^{\pm},\cH)$  in $\mathcal{M}or_{=1}$ along with a
vector bundle $E$ on $\tilde{C}$ to the same point of $\mathcal{M}or_{=1}$ and
the polarization $\varphi^*(\cH)^{\otimes m} \otimes \det(E)$ on $\tilde{C} \to
T$ as defined above. Since $\Map_{\mathcal{M}or_{=1}}(\tilde{\cC}, \BGL_1)$ has
affine relative diagonal over $\mathcal{M}or_{=1}$, it follows that the
equalizer $\cU \to \cW$ is an affine morphism.} \qed

\begin{notn} \label{notn: forgetful morphism c} We denote by $c:
    \GBun_{g,n,p}(N) \to \prestable_{g,n,p}$ the restriction of the flat
    forgetful morphism from \Cref{L:cBun algebraic}(c) to $\GBun_{g,n,p}(N)
    \subset \CBun_{g,n,p}(N)$.
\end{notn}

\begin{lemma} \label{lemma: torsion-free kontsevich stable is log-crepant} Let
    $(\varphi : \tilde{C} \to C,\sigma^\pm_\bullet)$ be a contraction of
    $(n,p)$-marked curves of genus $g$ defined over an algebraically closed
    field. Let $E$ be a vector bundle on $\tilde{C}$ such that
    $\omega_{\tilde{C}}(\sum \sigma_{\bullet}^{\pm}) \otimes \det(E)^{\otimes
    3}$ is $\varphi$-ample and $\varphi_*(E(-\sum_{i=1}^p \sigma_i^+))$ is pure
    of dimension $1$ on $C$. Then,
    \begin{enumerate}[(a)]
        \item There is an isomorphism $\varphi^\ast(\omega_C) \cong
        \omega_{\tilde{C}} (\sum_{i=1}^p \sigma^+_i)$. This implies that
        $\varphi$ is an isomorphism over an open subset $U \subset C$ whose
        complement consists only of nodes and positive marked points
        $\varphi(\sigma^+_i)$.
        \item If $x \in C \setminus U$ is a node, then $\varphi^{-1}(x)$ is a
        rational bridge.
        \item If $x \in C \setminus U$ is a positive marked point
        $\varphi(\sigma^+_i)$, then $\varphi^{-1}(x)$ is a rational tail
        consisting of a chain of $\bP^1$'s with $\sigma^+_i$ on the terminal
        component.
        \item The contraction $\varphi$ is log-crepant.
    \end{enumerate}
\end{lemma}
\begin{proof}
    The morphism $\varphi$ contracts a disjoint union of connected subcurves
    $\bigsqcup_i \tilde{C}_i \subset \tilde{C}$, where each $\tilde{C}_i$ is a
    prestable curve with smooth rational irreducible components and containing
    no cycles. Furthermore, each $\tilde{C}_i$ contains at most one marking, and
    if the image $\varphi(\tilde{C}_i)$ is a node of $C$, then $\tilde{C}_i$
    cannot contain any markings. A $\mathbb{P}^1$-component of $\tilde{C}_i$ is
    called terminal if the canonical bundle $\omega_{\tilde{C}}$ has negative
    degree on this component (i.e. it contains exactly one node of
    $\widetilde{C}$). All the items (a)-(d) follow if we can show that every
    terminal $\mathbb{P}^1$ on $\tilde{C}_i$ contains a positive marking
    $\sigma_j^+$.

    Suppose for the sake of contradiction that there is a terminal irreducible
    component $\mathbb{P}^1 \cong \Sigma \subset \tilde{C}_i$ that does not
    contain a positive marking. Since $\Sigma$ is terminal, it intersects the
    closure of its complement in $\tilde{C}$ at a single point $x$. Since
    $\omega_{\tilde{C}}(\sum \sigma_{\bullet}^{\pm}) \otimes \det(E)^{\otimes
    3}$ is ample on $\Sigma$ and $\omega_{\tilde{C}}(\sum
    \sigma_{\bullet}^{\pm})|_{\Sigma}$ has degree $-1$ or $0$ ($\Sigma$ cannot
    contain two markings), this forces $\det(E)$ to have positive degree on
    $\Sigma$. We know that $E|_{\Sigma} \cong \bigoplus_j
    \cO_{\mathbb{P}^1}(n_j)$ for some integers $n_j$. Since $\det(E)|_{\Sigma}$
    is ample, there is some index $j_0$ with $n_{j_0}>0$. Pick a nontrivial
    global section $s \in H^0(\Sigma, \cO_{\mathbb{P}^1}(n_{j_0}))$ which
    vanishes at the attaching point $x$. We can view $s$ as a section of
    $E|_{\Sigma}$ by setting it equal to $0$ on the other direct summands
    $\cO_{\mathbb{P}^1}(n_j)$ of $E|_{\Sigma}$. By the assumption that $\Sigma$
    contains no positive marking, we have $E(- \sum_{i=1}^p
    \sigma_i^+)|_{\Sigma} \cong E|_{\Sigma}$, so we can view $s$ as a section of
    $E(- \sum_{i=1}^p \sigma_i^+)|_{\Sigma}$.  We glue $s$ to the zero section
    of $E(- \sum_{i=1}^p \sigma_i^+)$ on the closure of the complement of
    $\Sigma$ in $\tilde{C}$ to obtain a nontrivial global section $\widetilde{s}
    \in H^0(\tilde{C}, E(- \sum_{i=1}^p \sigma_i^+))$ that vanishes on the
    complement of $\Sigma$, thus contradicting the assumption that
    $\varphi_*(E(-\sum_{i=1}^p \sigma_i^+))$ is pure of dimension $1$.
\end{proof}

\begin{lemma} \label{lemma: control of bubblings} If $(\varphi : \tilde{C} \to
    C,\sigma^\pm_\bullet, E)$ is a geometric point of $\GBun_{g,n,p}(N)$, then
    \begin{enumerate}[(a)]
        \item All the items in the statement of \Cref{lemma: torsion-free
        kontsevich stable is log-crepant} are satisfied.
        \item The length of each rational bridge and rational tail contracted by
        $\varphi$ is at most $N$.
        \item The restriction of $E$ to each $\mathbb{P}^1$-component contracted
        by the stabilization morphism is of the form
        $\mathcal{O}_{\mathbb{P}^1}(1)^{\oplus k} \oplus
        \mathcal{O}_{\mathbb{P}^1}^{\oplus N-k}$ for some $k>0$.
    \end{enumerate}
\end{lemma}
\begin{proof}
    For (a), note that the log-canonical bundle $\omega_{\tilde{C}}(\sum
    \sigma_{\bullet}^{\pm})$ is trivial on every fiber of $\varphi$, because
    $\varphi$ is log-crepant. In particular, the condition that $\det(E)$ is
    $\varphi$-ample is equivalent to $\omega_{\tilde{C}}(\sum
    \sigma_{\bullet}^{\pm}) \otimes \det(E)^{\otimes 3}$ being $\varphi$-ample.
    Therefore, the assumptions of \Cref{lemma: torsion-free kontsevich stable is
    log-crepant} are satisfied, and hence that the conclusions of \Cref{lemma:
    torsion-free kontsevich stable is log-crepant} hold for $\varphi$.

    Let $R \subset \widetilde{C}$ be a maximal connected subcurve that is
    contracted to a point under $\varphi$. By part (a), we have two
    possibilities: either $\varphi(R)$ is a node and $R$ is a rational bridge,
    or  $\varphi(R)$ is smooth and $R$ is a rational tail with a positive marked
    point $\sigma^+$ at the terminal irreducible component. In either case, $R$
    is a chain of rational curves of some length, say $\ell$. Choose an
    enumeration of the irreducible components $\Sigma_1, \Sigma_2, \ldots
    \Sigma_{\ell}$ such that $\Sigma_i$ intersects $\Sigma_{i+1}$ at exactly one
    point for all $i<\ell$. For each $\ell$-tuple $(n_1, n_2, \ldots, n_{\ell})$
    of integers, we denote by $\cO(\, n_1 \mid n_2 \mid \ldots \mid n_{\ell}\,)$
    the unique line bundle on $R$ which restricts to $\cO_{\mathbb{P}^1}(n_i)$
    on $\Sigma_i$. By \cite[Sect. 4]{martens_thaddeus_variations_grothendieck},
    there is an $N\times \ell$-tuple $((n_{i,j})_{i=1}^{\ell})_{j=1}^{N}$ such
    that $E|_R \cong \bigoplus_{j=1}^N \cO(\, n_{1,j} \mid n_{2,j}\mid \ldots
    \mid n_{\ell,j}\,)$. Such line bundle $\cO(\, n_1 \mid n_2 \mid \ldots \mid
    n_{\ell}\,)$ on $R$ is generated by global sections if and only if $n_i \geq
    0$ for all $i$. Hence, the surjectivity of $\varphi^*\varphi_*(E) \to E$
    implies that $n_{i,j} \geq 0$ for all $i,j$.

    In the case when $R$ is a rational bridge without marked points, then we
    have $E(-\sum_{i=1}^p \sigma_i^+))|_R \cong E|_R$. If we let $\{r_1, r_2\} =
    R \cap \overline{\widetilde{C} \setminus R}$, then
    $\varphi_*(E(-\sum_{i=1}^p \sigma_i^p))$ is pure of dimension 1 in an open
    neighborhood of $\varphi(R)$ if and only if $E|_R$ has no nonzero global
    sections that vanish at both $r_1$ and $r_2$. A direct computation shows
    that this holds if and only if for every separate index $j$ we have either
    \begin{itemize}
        \item $n_{i,j}= 0$ for all $i$, or
        \item there exists some $1 \leq q_j \leq \ell$ such that $n_{i,j} =
        \delta_{i,q_j}$ for all $i$, where $\delta$ is the Kronecker delta
        function.
    \end{itemize}
    Hence, if we define $\cL_{q}:= \cO(\,\delta_{1,q} \mid \delta_{2,q} \mid
    \ldots \mid \delta_{\ell,q}\,)$, then we have $E|_R \cong \cO_R^{\oplus a_0}
    \oplus \bigoplus_{q=1}^{\ell} \cL_q^{\oplus a_q}$ for some nonnegative
    integers $a_i \geq 0$. Since $\det(E)$ is $\varphi$-ample, we must have
    $a_q>0$ for all $q>1$, and hence parts (b) and (c) follow in the case when
    $R$ is a rational bridge.

    In the case when $R$ is a rational tail with a positive marking $\sigma^+$
    at the terminal component, then we have $E(-\sum_{i=1}^p \sigma_i^+))|_R
    \cong E(-\sigma^+)|_R$. If we let $\{r\} = R\cap \overline{\widetilde{C}
    \setminus R}$, then $\varphi_*(E(-\sum_{i=1}^p \sigma_i^p))$ is pure of
    dimension 1 in an open neighborhood of $\varphi(R)$ if and only if $E|_R$
    has no nonzero global sections that vanish at both $r$ and $\sigma^+$. The
    same reasoning as in the previous paragraph shows that (b) and (c) are
    satisfied in this case.
\end{proof}

\begin{defn}[Degree of Gieseker bundles] \label{notn: degree of a Gieseker
    bundle} We say that a $T$-family of Gieseker bundles $(\varphi : \tilde{C}
    \to C, \sigma^\pm_{\bullet}, E)$ has degree $d$ if the degree of the
    $T$-flat family of sheaves $\varphi_*(E(-\sum_{i=1}^p\sigma_i^+))$ on $C$ is
    $d$ as in \Cref{defn: degree and slope of sheaves on prestable curves}. For
    any integer $d$, we denote by $\CBun_{g,n,p}(N,d)$ the open and closed
    substack of $\CBun_{g,n,p}(N)$ consisting of those families of Gieseker
    bundles of degree $d$. Similarly, we set $\GBun_{g,n,p}(N,d) :=
    \CBun_{g,n,p}(N,d) \cap \GBun_{g,n,p}(N)$.
\end{defn}

\begin{lemma} \label{L: connected components of Gbun} The fibers of
    $\CBun_{g,n,p}(N,d) \to S$ are geometrically irreducible. In particular, the
    $S$-fibers of the open substack $\GBun_{g,n,p}(N,d) \subset
    \CBun_{g,n,p}(N,d)$ are also geometrically irreducible.
\end{lemma}
\begin{proof}
    After base-change, we may assume without loss of generality that $S =
    \overline{s}$ is the spectrum of an algebraically closed field. The
    forgetful morphism $\CBun_{g,n,p}(N,d) \to \prestable_{g,n,p}$ given by
    $(\varphi : \tilde{C} \to C, \sigma^\pm_{\bullet}, E) \mapsto (\tilde{C},
    \sigma^\pm_{\bullet})$ is smooth, because it is a composition of the smooth
    forgetful morphism $\CBun_{g,n,p}(N,d) \to \precont_{g,n,p}$ (\Cref{L:cBun
    algebraic}) and the smooth morphism $f_{src}: \precont_{g,n,p} \to
    \prestable_{g,n,p}$ (\Cref{lemma: properties of contraction stack}). Let
    $\cM_{g,n,p} \subset \prestable_{g,n,p}$ denote the open substack
    parameterizing smooth $(n,p)$-marked genus $g$ curves. Every point of
    $\prestable_{g,n,p}$ deforms to a $(n,p)$-marked curved in $\cM_{g,n,p}$
    \cite[\href{https://stacks.math.columbia.edu/tag/0E7Y}{Tag
    0E7Y}]{stacks-project}, and hence $\cM_{g,n,p} \subset \prestable_{g,n,p}$
    is dense. Since $\cM_{g,n,p}$ is known to be irreducible
    \cite{deligne_mumford}, it follows that $\prestable_{g,n,p}$ is also
    irreducible. To conclude the irreducibility of $\CBun_{g,n,p}(N,d)$, it
    suffices to show that a general fiber of the forgetful morphism
    $\CBun_{g,n,p}(N,d) \to \prestable_{g,n,p}$ is irreducible. Given a
    geometric point $\overline{t} = (\tilde{C}, \sigma_{\bullet}^{\pm})$ of
    $\cM_{g,n,p}$, the corresponding fiber
    $\left(\CBun_{g,n,p}(N,d)\right)_{\overline{t}}$ is isomorphic to the stack
    $\Bun_{N}(\tilde{C})_{d+pN}$ of vector bundles of rank $N$ and degree $d+pN$
    on $\tilde{C}$. Hence all geometric fibers over $\cM_{g,n,p}$ are
    irreducible, since $\Bun_{N}(\tilde{C})_{d+pN}$ is known to be irreducible
    \cite[Cor. A.5]{hoffmann-ks-rationality}.
\end{proof}

\subsection{Kontsevich stable maps to Grassmannians}

\begin{defn} \label{defn: quot scheme} Let $X$ be a scheme and let $\cF$ be a
    quasicoherent sheaf of finite type over $X$. We denote by $Q^N_{\mathcal{F}}
    \to X$ the Quot scheme of rank $N$ locally free quotients of $\mathcal{F}$.
    More formally, for any $f: T \to X$, the $T$-points of $Q^N_{\cF}$ are given
    by
    \[ Q^N_{\cF}(X \to T) = \left\{ \text{rank $N$ locally free quotients
    $f^*(\mathcal{F}) \twoheadrightarrow \mathcal{E}$ on $T$} \right\} / \cong
    \]
\end{defn}

\begin{lemma} \label{lemma: ample line bundle on quot} Let $X$ and $\mathcal{F}$
    as in Definition \ref{defn: quot scheme} above.
    \begin{enumerate}[(a)]
        \item $Q^N_{\cF}$ is represented by a projective scheme over $X$. The
        determinant bundle $\det(\mathcal{E}_{univ})$ of the universal quotient
        $\mathcal{E}_{univ}$ is $X$-very ample.
        \item Suppose that $\cF$ is finitely presented, let $X \to T$ be a
        projective morphism of finite presentation over a quasicompact scheme
        $T$, and let $\mathcal{O}_{X}(1)$ be a $T$-ample line bundle. Then, for
        all sufficiently large $m$, the line bundle
        $\text{det}(\mathcal{E}_{univ}) \otimes O_{X}(m)$ on $Q^N_{\cF}$ is
        $T$-very ample.
    \end{enumerate}
\end{lemma}
\begin{proof}
    \noindent (a). We may work Zariski locally on $X$, and choose a surjection
    $\mathcal{O}_{X}^{\oplus h} \twoheadrightarrow \mathcal{F}$. This induces a
    closed immersion $W^N_{\mathcal{F}} \hookrightarrow \text{Gr}(h,N) \times
    X$, and we conclude by standard properties of Grassmannians.\endnote{Since
    the functor of points of $Q^N_{\cF}$ is a Zariski sheaf, we can work Zariski
    locally on $X$ and assume that there is a surjection
    $\mathcal{O}_{X}^{\oplus h} \twoheadrightarrow \mathcal{F}$ for some natural
    number $h$. Let $\text{Gr}(h,N) \to S$ denote the Grassmannian of rank $N$
    quotients of $\cO_S^{\oplus h}$. Consider the morphism of functors
    \[Q^N_{\cF} \to \text{Gr}(h,N)\times_S X, \; \; \; \; \; \;
    [f^*(\mathcal{F}) \twoheadrightarrow \mathcal{E}] \mapsto
    [\mathcal{O}_T^{\oplus h} \twoheadrightarrow f^*(\mathcal{F})
    \twoheadrightarrow \mathcal{E}]\] This morphism is represented by a closed
    immersion, and hence it follows that $Q^N_{\cF}$ is $X$-projective. The
    determinant of the universal quotient on $\text{Gr}(h,N)$ is very ample, and
    so its pullback to $Q^N_{\cF}$ is $X$-very ample.}

    \noindent (b). By Noetherian approximmation, we may assume that $T$ is
    affine and Noetherian. Choose $m \gg 0$ such that $\mathcal{F}(m)$ is
    globally generated. Then, we get a closed immersion $Q^N_{\cF}
    \hookrightarrow \Gr(h,N) \times X$ as in part (a). We conclude by using the
    isomorphism $Q^N_{\cF} \xrightarrow{\sim} Q^N_{\cF(m)}$ induced by tensoring
    by $\cO_X(m)$. \endnote{By working locally on $T$ and using Noetherian
    approximation, we can reduce to the case when $T$ is affine and Noetherian.
    For sufficiently large $m$, the twist $\mathcal{F}(m)$ is globally
    generated. Tensoring by the pullback of $\mathcal{O}_X(m)$ induces an
    isomorphism of $T$-schemes $Q^N_{\cF} \xrightarrow{\sim} Q^N_{\cF(m)}$ given
    by $\mathcal{E} \mapsto \mathcal{E}(m)$. Choose a surjection
    $\mathcal{O}_{X}^{\oplus h} \to \mathcal{F}(m)$. By the argument in part
    (a), there is a closed immersion of $T$-schemes $Q^N_{\cF}
    \xrightarrow{\sim} Q^N_{\cF(m)} \hookrightarrow \text{Gr}(h,N) \times_S X$.
    Under this inclusion, the $T$-very ample line bundle
    $\mathcal{O}_{\text{Gr}(h,N)}(1) \boxtimes \mathcal{O}_{X}(1)$ pulls back to
    the $T$-very ample line bundle $\text{det}(\mathcal{E}_{univ}) \otimes
    \mathcal{O}_{X}(N\cdot m +1)$ on $Q^N_{\cF}$.}
\end{proof}

Let $\cY$ be a locally Noetherian stack over $S$, and fix a $\cY$-family $(C,
\overline{\sigma}_{\bullet}^{\pm})$ of $(n,p)$-marked genus $g$ prestable
curves. Let $\mathcal{F}$ be a sheaf of finite type over $C$.
\begin{defn} \label{defn: kontsevich stable maps} The stack
    $\cK_{g,n,p}(Q_{\cF}^N/\cY)$ over $\cY$ is the pseudofunctor that sends a
    $\cY$-scheme $T$ to the groupoid of $T$-families of $(n+p)$-pointed genus
    $g$ Kontsevich stable morphisms $(\tilde{C}, \sigma_{\bullet}^{\pm}) \to
    Q^N_{\mathcal{F}}\times_{\cY}T$ such that the composition $\varphi:
    (\tilde{C}, \sigma_{\bullet}^{\pm}) \to Q^N_{\mathcal{F}}\times_{\cY}T \to
    (C_T, \overline{\sigma}_{\bullet}^{\pm, T})$ is a contraction of
    $T$-families of $(n,p)$-marked curves as in \Cref{defn: contraction}.
\end{defn}

For the following definition, we denote by $\pi: \tilde{\cC} \to
\cK_{g,n}(Q_{\cF}^N/\cY)$ the universal source curve with markings
$\sigma_{\bullet}^{\pm}$, and let $E_{un}$ denote the vector bundle on
$\tilde{\cC}$ obtained by pulling back the universal quotient on $Q^N_{
\mathcal{F}}$.
\begin{defn}
    Let $\cH$ denote a line bundle on $C \to \cY$. For any $m \in \mathbb{Z}$,
    let $\mathcal{L}_{\Cor,\cH,m}$ denote the line bundle on
    $\cK_{g,n,p}(Q_{\cF}^N/B)$ given by
    \begin{gather*}
        \mathcal{L}_{\Cor,\cH,m} := \det\left(R\pi_*\left(\left(\left[\omega_{\pi}(\sum \sigma_{\bullet}^{\pm})\otimes \left(\det(E_{un}) \otimes \cH^{\otimes m}|_{\tilde{\cC}}\right)^{\otimes3}\right] - [\mathcal{O}_{\tilde{\cC}}]\right)^{\otimes 2}\right)\right),
    \end{gather*}
    where the pushforward and determinant are taken in the sense of K-theory.
\end{defn}

\begin{lemma} \label{lemma: stack kontsevich stable maps is union of projective
    DM stacks} The stack $\cK_{g,n,p}(Q_{\cF}^N/\cY)$ is smooth-locally on $\cY$
    a disjoint union of relatively DM $\cY$-projective stacks. Suppose
    furthermore that we are given a $\cY$-ample line bundle $\cH$ on $C$, that
    $\cY$ is quasicompact, and that $\mathcal{F}$ is finitely presented. Then
    for all sufficiently large $m$ the line bundle $\mathcal{L}_{\Cor, \cH, m}$
    restricts to a $\cY$-ample line bundle on any proper substack of
    $\cK_{g,n,p}(Q_{\cF}^N/\cY)$.
\end{lemma}
\begin{proof}
    Since the condition on the equality of markings $\varphi\circ \sigma_i^{\pm}
    = \overline{\sigma}_{\bullet}^{\pm}$ and the condition on having degree $1$
    are closed (\Cref{lemma: degree is open and closed}), it follows that
    $\cK_{g,n,p}(Q_{\cF}^N/\cY)$ is a closed substack of the stack
    $\cM_{g,n}(Q_{\cF}^N/\cY)$ of all Kontsevich stable maps. The desired
    assertions follow from the corresponding statements proven for
    $\cM_{g,n}(Q_{\cF}^N/\cY)$ in \cite[Prop. 2.9]{gauged_theta_stratifications}
    and the fact that $\det(\mathcal{E}_{univ}) \otimes O_{\tilde{C}}(m)$ is
    $\cY$-very ample for $m\gg0$ (\Cref{lemma: ample line bundle on
    quot}).\qedhere
\end{proof}

\begin{lemma} \label{lemma: criterion kontsevich stability} Let $T \to \cY$ be a
    scheme. Let $(\tilde{C}, \sigma_{\bullet}^{\pm})$ be a $T$-family of
    $(n,p)$-marked prestable genus $g$ curves equipped with a $T$-morphism
    $\tilde{C} \to Q^N_{\mathcal{F}}\times_{\cY}T $ such that the composition
    $\varphi: (\tilde{C}, \sigma_{\bullet}^{\pm}) \to (C_T,
    (\overline{\sigma}_{\bullet}^{\pm})_T)$ is a contraction of $(n,p)$-marked
    curves. Let $E$ be the vector bundle on $\tilde{C}$ obtained as the pullback
    of the universal quotient on $Q^N_{\mathcal{F}}$. Then $(\tilde{C},
    \sigma_{\bullet}^{\pm}) \to Q^N_{\mathcal{F}}\times_{\cY}T $ is a $T$-family
    of Kontsevich stable maps if and only if $\omega_{\tilde{C}}(\sum
    \sigma_{\bullet}^{\pm}) \otimes \det(E)^{\otimes 3}$ is $\varphi$-ample. In
    particular, if $\varphi$ is log-crepant, then Kontsevich stability is
    equivalent to $\det(E)$ being $\varphi$-ample.
\end{lemma}
\begin{proof}
    Since relative ampleness of $\omega_{\tilde{C}}(\sum \sigma_{\bullet}^{\pm})
    \otimes \det(E)^{\otimes 3}$ can be checked on fibers \cite[Cor.
    9.6.5]{egaiv}, we can reduce to the case when $T=\cY$ is the spectrum of an
    algebraically closed field. The composition $\varphi: (\tilde{C},
    \sigma_{\bullet}^{\pm}) \to (C, \overline{\sigma}_{\bullet}^{\pm})$ only
    contracts $\mathbb{P}^1$-components, and the group of automorphisms of
    $(\tilde{C}, \sigma_{\bullet}^{\pm}) \to Q_{\mathcal{F}}$ is finite if and
    only if $\tilde{C} \to Q_{\mathcal{F}}$ is a finite morphism when restricted
    to each marked-unstable $\mathbb{P}^1$-component of $(\tilde{C},
    \sigma_{\bullet}^{\pm})$ that is contracted by $\varphi$. This happens if
    and only if the pullback $\det(E)$ of the $C$-ample line bundle
    $\det(\cE_{univ})$ is ample on each marked-unstable $(\tilde{C},
    \sigma_{\bullet}^{\pm})$ contracted by $\varphi$. Since the log-canonical of
    a marked-unstable component has degree at least $-2$ and is positive on the
    marked-stable components, this is equivalent to the $\varphi$-ampleness of
    $\omega_{\tilde{C}}(\sum \sigma_{\bullet}^{\pm}) \otimes \det(E)^{\otimes
    3}$. Finally, if $\varphi$ is log-crepant, then the log-canonical bundle
    $\omega_{\tilde{C}}(\sum \sigma_{\bullet}^{\pm})$ is trivial on each fiber
    of $\varphi$, and therefore the condition is equivalent to $\det(E)$ being
    $\varphi$-ample.
\end{proof}

\subsection{Properness over the stack of pure marked sheaves} \label{subsection:
properness over the stack of pure marked sheaves} Let $(\varphi : \tilde{C} \to
C, \sigma^\pm_\bullet, E)$ be a $T$-point of $\CBun_{g,n,p}(N)$. Set $F:=
\varphi_*(E)$. For each $i$, we may push forward the morphism $E \to
(\sigma_i^+)_*(E|_{\sigma_i^+})$ to obtain a morphism $\varphi_*(E) \to
\varphi_*(\sigma_i^+)_*(E|_{\sigma_i^+})$. By adjunction, this yields a morphism
$(\varphi \circ \sigma_{i}^+)^*(\varphi_*(E)) \to E|_{\sigma_i^+}$. Note that
$W_i:= E|_{\sigma_i^+}$ is a rank $N$ vector bundle. Therefore, we get the data
$(\varphi_*(E), (\varphi \circ \sigma_{\bullet}^+)^*(\varphi_*(E)) \to
E|_{\sigma_{\bullet}^+})$ on the $(n,p)$-marked curve $(C, \varphi \circ
\sigma_{\bullet}^{\pm})$.
\begin{lemma} \label{lemma: Phis is well-defined} The assignment
    \[(\varphi : \tilde{C} \to C, \sigma^\pm_\bullet, E) \mapsto (C, \varphi
    \circ \sigma^\pm_\bullet, \varphi_*(E), (\varphi \circ
    \sigma_{\bullet}^+)^*(\varphi_*(E)) \to E|_{\sigma_{\bullet}^+})\] defines a
    morphism of algebraic stacks $\Phi : \CBun_{g,n,p}(N) \to \Coh_{g,n,p}(N)$.
    The substack $\GBun_{g,n,p}(N)$ is the preimage of the open substack
    $\Pur_{g,n,p}(N)$ under this morphism.
\end{lemma}
\begin{proof}
    For any $T$-family $(\varphi : \tilde{C} \to C, \sigma^\pm_\bullet, E)$ on
    $\GBun_{g,n,p}(N)$, \Cref{lemma: pushforward lemma} implies that
    $\varphi_*(E)$ is $T$-flat, and that its formation commutes with
    $T$-base-change. We have a short exact sequence
    \[0 \to E(-\sum \sigma_{\bullet}^{+}) \to E \to \bigoplus_{i=1}^p
    (\sigma_i^+)_*(E|_{\sigma_i^+}) \to 0 \] By applying $\varphi_*(-)$ and
    using the vanishing of $R\varphi_*(E(-\sum \sigma_{\bullet}^{+}))$
    (\Cref{lemma: pushforward lemma}), it follows that $\varphi_*(E) \to
    \bigoplus_{i=1}^p \varphi_*(\sigma_i^+)_*(E|_{\sigma_i^+})$ is surjective.
    In particular, for each $i$, the morphism $\varphi_*(E) \twoheadrightarrow
    \varphi_*(\sigma_i^+)_*(E|_{\sigma_i^+})$ is surjective. Hence each morphism
    $(\varphi \circ \sigma_{\bullet}^+)^*(\varphi_*(E)) \to
    E|_{\sigma_{\bullet}^+}$ obtained by adjunction is also surjective.
    Therefore, the data $(\varphi_*(E), (\varphi \circ
    \sigma_{\bullet}^+)^*(\varphi_*(E)) \to E|_{\sigma_{\bullet}^+})$ is a
    marked sheaf on the $(n,p)$-marked curve $(C, \varphi \circ
    \sigma_{\bullet}^{\pm})$, and its formation commutes with $T$-base-change.
    \qedhere

\end{proof}

Recall that there is a smooth surjective morphism $\prestablepol_{g,n,p} \to
\prestable_{g,n,p}$ from the stack $\prestablepol_{g,n,p}$ of polarized
$(n,p)$-marked prestable curves
\cite[\href{https://stacks.math.columbia.edu/tag/0DQC}{Tag
0DQC}]{stacks-project}. Recall that we have stacks $\GBun^{\pol}_{g,n,p}(N)$ and
$\Pur^{\pol}_{g,n,p}(N)$ given by the following diagram with Cartesian squares
\[\begin{tikzcd} \GBun^{\pol}_{g,n,p}(N) \ar[r, "\Phi"] \ar[d] &
        \Pur^{\pol}_{g,n,p}(N) \ar[r] \ar[d] & \prestablepol_{g,n,p} \ar[d] \\
        \GBun_{g,n,p}(N) \ar[r,"\Phi"] & \Pur_{g,n,p}(N) \ar[r]  &
    \prestable_{g,n,p} \end{tikzcd}\]

Let $(\cal{C}, \overline{\sigma}^{\pm}_{\bullet}, \cH) \to
\prestablepol_{g,n,p}$ be the universal polarized $(n,p)$-marked prestable
curve. We abuse notation and also denote by $(\cal{C},
\overline{\sigma}^{\pm}_{\bullet}, \cH) \to \GBun^{\pol}_{g,n,p}(N)$ the base
change to $\GBun^{\pol}_{g,n,p}(N)$ whenever the meaning is clear from context.
By definition, there is a family $\pi: (\widetilde{\cC}, \sigma_{\bullet}^{\pm})
\to \GBun^{\pol}_{g,n,p}(N)$ of $(n,p)$-marked curves equipped with a
log-crepant contraction $\varphi: (\widetilde{C}, \sigma^{\pm}_{\bullet}) \to
(\cal{C}, \overline{\sigma}^{\pm}_{\bullet})$, and a universal rank $N$ vector
bundle $E_{un}$ on $\widetilde{C}$.

\begin{defn} \label{defn: cornalba line bundle on gbunpol} For any $m \in
    \mathbb{Z}$, let $\mathcal{L}_{\Cor,m}$ denote the line bundle on
    $\GBun^{\pol}_{g,n,p}(N)$ given by
    \begin{gather*}
        \mathcal{L}_{\Cor,m} := \det\left(R\pi_*\left(\left(\left[\omega_{\pi}(\sum \sigma_{\bullet}^{\pm})\otimes \left(\det(E_{un}) \otimes \cH^{\otimes m}|_{\tilde{\cC}}\right)^{\otimes3}\right] - [\mathcal{O}_{\tilde{\cC}}]\right)^{\otimes 2}\right)\right).
    \end{gather*}
\end{defn}

\begin{remark}
    The bilinearity of the Deligne pairing \cite{dolce_deligne_pairing} implies that the assignment $m \mapsto
    \mathcal{L}_{\Cor, m}$ defines a polynomial of degree $2$ in the variable $m$ with values
    in the Picard group of $\GBun^{\pol}_{g,n,p}(N)$. In other words, we have
    $\mathcal{L}_{\Cor, m}  = A_2^{\otimes m^2} \otimes A_1^{\otimes m} \otimes A_0$ for some fixed line
    bundles $A_i$ on $\GBun^{\pol}_{g,n,p}(N)$.
\end{remark}

Our main result for this section is the following.
\begin{thm} \label{T:projective_morphism} The morphism $\Phi :
    \GBun_{g,n,p}^{\pol}(N) \to \Pur^{\pol}_{g,n,p}(N)$ is relatively DM and
    locally projective. For any given quasi-compact open substack $\cW \subset
    \Pur^{\pol}_{g,n,p}(N)$, there exists some $m\gg 0$ such that the
    restriction of the line bundle $\cL_{\Cor,m}$ to $\Phi^{-1}(\cW)$ is
    relatively ample on all $\cW$-fibers.
\end{thm}

We will need some preliminary preparation before the proof of
\Cref{T:projective_morphism}. Let $(\cC, \cF, (\sigma_{\bullet}^+)^*(\cF) \to
\cW_{\bullet})$ be the universal family of pure marked sheaves on
$\Pur^{\pol}_{g,n,p}(N)$. The corresponding Quot scheme $Q^N_{\cF}$ is schematic
and locally projective over $\Pur^{\pol}_{g,n,p}(N)$ (\Cref{lemma: ample line
bundle on quot}). Let $\cK_{g,n,p}(N):=
\cK_{g,n,p}(Q^N_{\cF}/\Pur^{\pol}_{g,n,p}(N))$ denote the substack of the stack
of $(n,p)$-marked stable Kontsevich morphisms as in \Cref{defn: kontsevich
stable maps}.

\begin{lemma} \label{lemma: open immersion of GBun into kontsevich stable maps}
    There is a locally closed immersion $\GBun_{g,n,p}^{\pol}(N) \hookrightarrow
    \cK_{g,n,p}(N)$ such that the composition $\GBun_{g,n,p}^{\pol}(N)
    \hookrightarrow \cK_{g,n,p}(N) \to \Pur^{\pol}_{g,n,p}(N)$ is isomorphic to
    $\Phi$.
\end{lemma}
\begin{proof}
    By \Cref{lemma: criterion kontsevich stability}, for any scheme $T$, the
    $T$-points of the stack $\cK_{g,n,p}(N)$ classify the data for a contraction
    of $(n,p)$-marked curves $(\varphi: \tilde{C} \to C,
    \sigma_{\bullet}^{\pm})$, a $T$-ample line bundle $\cH$ on $C$, a pure
    marked sheaf $(F, (\sigma_{\bullet}^+)^*(F) \twoheadrightarrow
    W_{\bullet})$, and a rank $N$ quotient $\varphi^*(F) \twoheadrightarrow E$
    such that $\omega_{\tilde{C}}(\sum \sigma_{\bullet}^{\pm}) \otimes
    \det(E)^{\otimes 3}$ is $\varphi$-ample. Note that in this case the counit
    $\varphi^*\varphi_*(E) \to E$ is surjective on very $T$-fiber, and its
    formation commutes with base-change by \Cref{lemma: pushforward lemma}.
    Adjunction induces a natural morphism $F \to \varphi_*(E)$ of $T$-flat
    sheaves. The requirement that this induces a morphism of short exact
    sequences of $T$-flat sheaves as below
    \[
        \begin{tikzcd}
            0 \ar[r] & K \ar[r] \ar[d, dashed] & F \ar[r] \ar[d] & \bigoplus \varphi_*(\sigma_{\bullet}^+)_*(W_i) \ar[r] \ar[d, dashed] & 0 \\
            0 \ar[r] & \varphi_*(E(-\sum \sigma_{\bullet}^+)) \ar[r]  & \varphi_*(E)  \ar[r] & \bigoplus \varphi_*(\sigma_{\bullet}^{+})_*(E|_{\sigma_{\bullet}^{+}})  \ar[r] & 0
        \end{tikzcd}
    \]
    is a closed condition (here the surjectivity of the bottom short exact
    sequence follows from the vanishing of $R^1\varphi_*(E(-\sum
    \sigma_{\bullet}^+))$ proven in \Cref{lemma: pushforward lemma}). We denote
    by $\cU \subset \cK_{g,n,p}(N)$ the locally closed substack where the
    morphism $F \to \varphi_*(E)$ induces an isomorphism of short exact sequence
    as above. In particular, this means that the data of the pure marked sheaf
    $(F, (\sigma_{\bullet}^+)^*(F) \twoheadrightarrow W_i)$ is redundant and can
    be recovered from $E$. Furthermore, note that for points in $\cU \subset
    \cK_{g,n,p}(N)$ the pushforward $\varphi_*(E(-\sum \sigma_{\bullet}^+))$ is
    fiberwise pure of dimension 1, and hence by \Cref{lemma: torsion-free
    kontsevich stable is log-crepant} the morphism $\varphi$ is forced to be
    log-crepant. Therefore, the condition on $\omega_{\tilde{C}}(\sum
    \sigma_{\bullet}^{\pm}) \otimes \det(E)^{\otimes 3}$ being $\varphi$-ample
    is equivalent to the condition that $\det(E)$ is $\varphi$-ample. We
    conclude that $\cU \subset \cK_{g,n,p}(N)$ is isomorphic to
    $\GBun_{g,n,p}^{\pol}(N)$.
\end{proof}

\vspace{-0.35cm}

\begin{lemma} \label{L: valuative criterion over torsion-free moduli} 
The morphism $\Phi: \GBun_{g,n,p}(N) \to \Pur_{g,n,p}(N)$ satisfies the
    existence part of the valuative criterion for properness for discrete
    valuation rings.
\end{lemma}
\begin{proof}
    Since $\prestablepol_{g,n,p} \to \prestable_{g,n,p}$ is smooth and
    surjective (\Cref{prop: properties of the stack of polarized curves}), we
    may base-change to $\prestablepol_{g,n,p}$ and check the valuative criterion
    for $\Phi: \GBun_{g,n,p}^{\pol}(N) \to \Pur^{\pol}_{g,n,p}(N)$ instead. Let
    $R$ be a discrete valuation ring. We denote by $\eta$ and $s$ the generic
    and special point of $\Spec(R)$ respectively. Let $\Spec(R) \to
    \Pur^{\pol}_{g,n,p}(N)$ be an $R$-point corresponding to a tuple $(C,
    \overline{\sigma}_{\bullet}^{\pm},\cH, F, (\overline{\sigma}_{\bullet}^+)^*F
    \twoheadrightarrow W_i)$ over $\Spec(R)$. Suppose that we are given a
    $\Phi$-lift defined over the generic point $\eta$. Since the morphism
    $\cK_{g,n,p}(N) \to \Pur^{\pol}_{g,n,p}(N)$ is represented by a union of
    projective relatively DM stacks smooth-locally on the target, it follows
    that, after perhaps passing to an extension $R$, we can extend the
    composition $\eta \to\GBun_{g,n,p}^{\pol}(N) \hookrightarrow \cK_{g,n,p}(N)$
    to lift $\Spec(R)\to \cK_{g,n,p}(N)$ corresponding to a contraction of
    $(n,p)$-marked curves $\varphi: (\tilde{C}, \sigma_{\bullet}^{\pm}) \to (C,
    \overline{\sigma}_{\bullet}^{\pm})$ and a rank $N$ quotient $\varphi^*(F)
    \twoheadrightarrow E$. To conclude the proof, it suffices to check that
    $\Spec(R) \to \cK_{g,n,p}(N)$ factors through the locally closed substack
    $\GBun_{g,n,p}^{\pol}(N)$.

    Adjunction induces a morphism $F \to \varphi_*(E)$. By the proof of
    \Cref{lemma: open immersion of GBun into kontsevich stable maps}, it
    suffices to show that this morphism induces an isomorphism of short exact
    sequences of $T$-flat sheaves
    \begin{equation} \label{diagram 1}
        \begin{tikzcd}[ampersand replacement=\&]
            0 \ar[r] \& K \ar[r] \ar[d, dashed] \& F \ar[r] \ar[d] \& \bigoplus \varphi_*(\sigma_{\bullet}^+)_*(W_i) \ar[r] \ar[d, dashed] \& 0 \\
            0 \ar[r] \& \varphi_*(E(-\sum \sigma_{\bullet}^+)) \ar[r]  \& \varphi_*(E)  \ar[r] \& \bigoplus \varphi_*(\sigma_{\bullet}^{+})_*(E|_{\sigma_{\bullet}^{+}})  \ar[r] \& 0
        \end{tikzcd}
    \end{equation}
    First, the existence of the vertical morphisms is equivalent to the
    composition $K \to F \to \varphi_*(E) \to \bigoplus
    \varphi_*(\sigma_{\bullet}^{+})_*(E|_{\sigma_{\bullet}^{+}})$ being the zero
    morphism. This, in turn, is equivalent to the adjoint morphism $(\varphi
    \circ \sigma_i^{+})^*(K) \to E|_{\sigma_i^+}$ being zero for all $i$. Note that this is
    a morphism of sheaves on $\Spec(R)$ whose codomain is a vector bundle. Since
    the morphism $\eta \to \cK_{g,n,p}(N)$ factors through
    $\GBun_{g,n,p}^{\pol}(N)$, the restriction of $(\varphi \circ
    \sigma_i^{+})^*(K) \to E|_{\sigma_i^+}$ to the generic fiber is zero. It follows that
    $(\varphi \circ \sigma_i^{+})^*(K) \to E|_{\sigma_i^+}$ must be zero, and hence we get
    the existence of the vertical arrows in the diagram \eqref{diagram 1} above.
    We are left to show that both dashed vertical arrows in diagram
    \eqref{diagram 1} are isomorphisms.

    Since $K \to \varphi_*(E(-\sum \sigma_{\bullet}^+))$ is a morphism of
    $R$-flat sheaves which fiberwise have the same Hilbert polynomial with
    respect to $\cH$, to conclude that it is an isomorphism it suffices to show
    that it is injective. This follows because the source $K$ is fiberwise pure,
    and the morphism $K \to \varphi_*(E(-\sum \sigma_{\bullet}^+))$ is
    generically an isomorphism on every fiber by construction. \endnote{We start
    by showing that $K \to \varphi_*(E(-\sum \sigma_{\bullet}^+))$ is an
    isomorphism. This is a morphism of $R$-flat sheaves on $C \to \Spec(R)$
    which fiberwise have the same Hilbert polynomial with respect to the
    polarization $\cH$. Therefore, it suffices to show that for each $R$-fiber
    the morphism $K \to \varphi_*(E(-\sum \sigma_{\bullet}^+))$ is injective. We
    already know that it is an isomorphism over the generic fiber $\eta \in
    \Spec(R)$, so we focus on the special fiber. Consider the restriction
    $K|_{C_s} \to \varphi_*(E(-\sum \sigma_{\bullet}^+))|_{C_s}$ to the special
    fiber. Note that $\varphi_s: \tilde{C}_s \to C_s$ is an isomorphism over a
    dense open subscheme $U \subset C_s$. By construction, the morphism
    $K|_{C_s} \to \varphi_*(E(-\sum \sigma_{\bullet}^+))|_{C_s}$ is an
    isomorphism over $U$. Since the family of marked sheaves $\Spec(R) \to
    \Pur^{\pol}_{g,n,p}(N)$ is pure, it follows that $K|_{C_s}$ is pure of
    dimension 1 on $C_s$. Therefore, the generically injective morphism
    $K|_{C_s} \to \varphi_*(E(-\sum \sigma_{\bullet}^+))|_{C_s}$ is forced to be
    injective, as desired.}

    We are left to show that $\bigoplus \varphi_*(\sigma_{\bullet}^+)_*(W_i) \to
    \bigoplus
    \varphi_*(\sigma_{\bullet}^{+})_*(E|_{\sigma_{\bullet}^{+}})|_{C_s}$ is an
    isomorphism. By applying adjunction to diagram \eqref{diagram 1} we obtain a
    diagram of exact sequences on $\tilde{C}$
    \[
        \begin{tikzcd}
            & \varphi^*(K) \ar[r] \ar[d] & \varphi^*(F) \ar[r] \ar[d] & \bigoplus (\sigma_{\bullet}^+)_*(W_i) \ar[r] \ar[d] & 0 \\
            0 \ar[r] & E(-\sum \sigma_{\bullet}^+) \ar[r]  & E  \ar[r] & \bigoplus (\sigma_{\bullet}^{+})_*(E|_{\sigma_{\bullet}^{+}})  \ar[r] & 0
        \end{tikzcd}
    \]
    Since both $\varphi^*(F) \to E$ and $E \to \bigoplus
    (\sigma_{\bullet}^{+})_*(E|_{\sigma_{\bullet}^{+}})$ are surjective, it
    follows that right-most vertical morphism $\bigoplus
    (\sigma_{\bullet}^+)_*(W_i) \to \bigoplus
    (\sigma_{\bullet}^{+})_*(E|_{\sigma_{\bullet}^{+}})$ is surjective. Since it
    is a morphism of pushforwards of vector bundles on $\Spec(R)$ of the same
    rank, we must have that $\bigoplus (\sigma_{\bullet}^+)_*(W_i) \to \bigoplus
    (\sigma_{\bullet}^{+})_*(E|_{\sigma_{\bullet}^{+}})$ is an isomorphism.
    Pushing forward via $\varphi$, we obtain that $\bigoplus
    \varphi_*(\sigma_{\bullet}^+)_*(W_i) \to \bigoplus
    \varphi_*(\sigma_{\bullet}^{+})_*(E|_{\sigma_{\bullet}^{+}})|_{C_s}$ is an
    isomorphism.
\end{proof}

\begin{notn} \label{notn: fixed degree pure equirank marked sheaves} For any
    given integer $d$, let $\Pur_{g,n,p}(N,d) \subset \Pur_{g,n,p}(N)$ denote
    the open and closed substack consisting of points $(\varphi: \tilde{C} \to
    C, \sigma_{\bullet}^{\pm}, F, (\sigma_{\bullet}^+)^*(F) \to W_{\bullet})$
    such that the kernel $K$ of the corresponding surjection $F
    \twoheadrightarrow \bigoplus_{i =1}^p(\sigma_i^+)_*(W_i)$ has degree $d$ as
    in \Cref{defn: degree and slope of sheaves on prestable curves}. Note that
    the morphism $\Phi: \GBun_{g,n,p}(N,d) \to \Pur_{g,n,p}(N)$ factors through
    $\Pur_{g,n,p}(N,d) \subset \Pur_{g,n,p}(N)$.
\end{notn}

\begin{proof}[Proof of \Cref{T:projective_morphism}]
    By \Cref{lemma: open immersion of GBun into kontsevich stable maps}, there
    is locally closed embedding $\GBun_{g,n,p}^{\pol}(N) \hookrightarrow
    \cK_{g,n,p}(N)$ into the stack $\cK_{g,n,p}(N)$. By \Cref{L: valuative
    criterion over torsion-free moduli}, the composition
    $\GBun_{g,n,p}^{\pol}(N) \hookrightarrow \cK_{g,n,p}(N) \to
    \Pur^{\pol}_{g,n,p}(N)$ satisfies the existence part of the valuative
    criterion for properness. Since the morphism $\cK_{g,n,p}(N) \to
    \Pur_{g,n,p}(N)$ is separated (\Cref{lemma: stack kontsevich stable maps is
    union of projective DM stacks}), it follows that $\GBun_{g,n,p}^{\pol}(N)
    \hookrightarrow \cK_{g,n,p}(N)$ is a closed immersion. In view of
    \Cref{lemma: stack kontsevich stable maps is union of projective DM stacks},
    in order to conclude the proof it suffices to show that $\Phi:
    \GBun_{g,n,p}^{\pol}(N) \to \Pur^{\pol}_{g,n,p}(N)$ is quasicompact.

    It suffices to show quasicompactness of the corresponding morphism of stacks
    without polarizations $\Phi: \GBun_{g,n,p}(N,d) \to \Pur_{g,n,p}(N,d)$ for
    each degree $d$ separately. By \Cref{L: connected components of Gbun}, the
    smooth $S$-stack $\GBun_{g,n,p}(N,d)$ has finitely many irreducible
    components. Since we may check quasicompactness smooth locally on the
    target, we may base-change to a quasicompact stack $\cU$ equipped with a
    smooth morphism $\cU \to \Pur_{g,n,p}(N,d)$. We may assume that the morphism
    factors through $\cU \to \Pur^{\pol}_{g,n,p}(N)$ (\Cref{prop: properties of
    the stack of polarized curves}) and that $\cK_{g,n,p}(N)
    \times_{\Pur^{\pol}_{g,n,p}(N)} \cU$ is a disjoint union of proper stacks
    over $\cU$ (\Cref{lemma: stack kontsevich stable maps is union of projective
    DM stacks}). Note that the morphism $\cU \to \Pur_{g,n,p}(N,d)$ is
    quasicompact, since $\Pur_{g,n,p}(N,d)$ is quasiseparated. The base change
    $\GBun_{g,n,p}(N,d)\times_{\Pur_{g,n,p}(N,d)} \cU$ is a stack which is
    smooth and quasicompact over the stack $\GBun_{g,n,p}(N,d)$, and hence it
    has finitely many irreducible components. It follows that the closed
    immersion $\GBun_{g,n,p}(N,d) \times_{\Pur_{g,n,p}(N,d)} \cU \hookrightarrow
    \cK_{g,n,p}(N) \times_{\Pur^{\pol}_{g,n,p}(N)} \cU$ must factor through
    finitely many of the closed and open substacks of $\cK_{g,n,p}(N)
    \times_{\Pur^{\pol}_{g,n,p}(N)} \cU$ that are proper over $\cU$. We conclude
    that $\GBun_{g,n,p}(N,d)\times_{\Pur_{g,n,p}(N,d)} \cU$ is quasicompact.
\end{proof}

\subsection{\texorpdfstring{$\Theta$}{Theta}-stratification for the moduli of
Gieseker bundles}
\begin{notn}
    By a formal number we mean a polynomial in $\mathbb{R}[\epsilon]$ with
    variable $\epsilon$. For any two numbers $f, g
    \in
    \mathbb{R}[\epsilon]$, we say that $f
    \leq
    g$ if there exists $H>0$ such that $f(\epsilon)
    \leq
    g(\epsilon)$ for all real numbers $0<\epsilon < H$.

\end{notn}

\begin{defn}[Evaluation morphisms] \label{defn: evaluation morphismsm} We define
    positive and negative evaluation morphisms as follows:
    \[ \ev_-: \GBun_{g,n,p}(N) \to (\BGL_N)^n, \; \; \; (\varphi : \tilde{C} \to
    C, \sigma^\pm_\bullet, E) \mapsto
    \left((\sigma_i^-)^*(E)\right)_{i=1}^n;\]
    \[\ev_+: \GBun_{g,n,p}(N) \to (\BGL_N)^p, \; \; \; (\varphi : \tilde{C} \to
    C, \sigma^\pm_\bullet, E) \mapsto
    \left((\sigma_i^+)^*(E)\right)_{i=1}^p.\]
\end{defn}

We view $\GBun_{g,n,p}(N)$ as a stack over $\prestable_{g,n,p} \times
\left(\BGL_N\right)^p$ via the morphism $c\times \ev_+: \GBun_{g,n,p} \to
\prestable_{g,n,p} \times \left(\BGL_N\right)^p$, where $c$ is the forgetful
morphism from \Cref{notn: forgetful morphism c}. By definition, we have a
commutative diagram
\[\begin{tikzcd} \Pur_{g,n,p}(N)   \ar[r] \ar[d, "\ev^+"] & \Pur_{g,n,p} \ar[d,
        "\ev^+"]  \\
        \left(\BGL_N\right)^p \ar[r, symbol = \hookrightarrow] & \SVect^p,
    \end{tikzcd} \] and the resulting map $\Pur_{g,n,p}(N) \to
    \prestable_{g,n,p} \times \left(\BGL_N\right)^p$ fits into a factorization
\[c \times \ev_+: \GBun_{g,n,p} \xrightarrow{\Phi} \Pur_{g,n,p}(N) \to
\prestable_{g,n,p} \times \left(\BGL_N\right)^p.\]

Just as in the case of marked sheaves, we work with the polarized version
$\GBun_{g,n,p}^{\pol}(N)$ of the moduli problem to define the numerical
invariant and the resulting stratification.

\begin{notn} \label{notn: line bundles and rational norm on Gieseker bundles} We
    abuse notation and denote by $\cL_0$ the line bundle on
    $\GBun_{g,n,p}^{\pol}(N)$ obtained by pulling back $\cL_0$ (\Cref{notn: line
    bundles on Pur}) via the relatively DM morphism $\GBun_{g,n,p}^{\pol}(N)
    \xrightarrow{\Phi} \Pur_{g,n,p}^{\pol}(N) \hookrightarrow
    \Pur_{g,n,p}^{\pol}$. Similarly, we denote by $b$ the rational nondegenerate
    norm on graded points of $\GBun_{g,n,p}^{\pol}(N)$ relative to
    $\prestable_{g,n,p}^{\pol} \times \left(\BGL_N\right)^p$ obtained by pulling
    back the norm $b$ on $\Pur_{g,n,p}$ as in \Cref{defn: rational norm on Pur}.
\end{notn}

We use the line bundles and the rational quadratic norm $b$ to define a
numerical invariant (cf. \cite[Sect. 2.4]{gauged_theta_stratifications}). Recall
the polynomial family of line bundles $\cL_{\Cor,m}$ on
$\GBun_{g,n,p}^{\pol}(N)$ that we introduced in \Cref{defn: cornalba line bundle
on gbunpol}.
\begin{defn} \label{defn: numerical invariants in Gbunpol} We define a formal
    numerical invariant $\nu_{\epsilon}$ on
    $\GBun_{g,n,p}^{\pol}(N)$ with values in $\bR[\epsilon]$
    relative to $\prestablepol_{g,n,p}\times(\BGL_N)^p$ given by $\nu_{\epsilon}
    = \left(\wt(\cL_{0}) - \epsilon^3 \wt(\cL_{\Cor,\epsilon^{-1}})\right) / \sqrt{b}$.
\end{defn}

Given a component $\gamma$ of $\Grad(\cX)$, if $u : \Grad(\cX)_\gamma \to \cX$
is the forgetful morphism, and $b$ is the norm on graded points, one has a
shifted numerical invariant
\[
    \nu^\gamma_\epsilon(f) = \nu_\epsilon(f) - \frac{\wt_\gamma(u^\ast(L_\epsilon))}{\lVert \gamma \rVert_b^2} (\gamma, \gr(f))_b
\]
The pairing $(\gamma, \gr(f))_b$ is interpreted as follows: for any filtration
$f : \Theta_k \to \Grad(\cX)_\gamma$, $f(0)$ has a canonical action of
$\bG_m^2$, where the first factor of $\bG_m$ acts by $\gamma$ and the second
acts by the cocharacter coming from $f : \Aut_\Theta(0) \to \Aut_\cX(f(0))$. The
quadratic norm on graded points

As a consequence of \Cref{thm: weak theta stratification pure marked sheaves}, we have the following:
\begin{thm} \label{thm: theta stratification for Gieseker bundles} For any $N > 0$ and $d$, the numerical invariant $\nu_{\epsilon}$
    defines a weak $\Theta$-stratification of $\GBun^{\pol}_{g,n,p}(N,d)$
    relative to $\prestable^{\pol}_{g,n,p} \times (\BGL_N)^p$ satisfying the
    following:
    \begin{enumerate}
        \item For every $\gamma(\epsilon) \in \bR[\epsilon]$, the open
        substack $\GBun^{\pol}_{g,n,p}(N,d)_{\nu_\epsilon \leq \gamma} \subseteq
        \GBun^{\pol}_{g,n,p}(N,d)$ is of finite type over
        $\prestable^{\pol}_{g,n,p} \times (\BGL_N)^p$. In particular, the
        relative weak $\Theta$-stratification is well-ordered.
        \item If the base $S$ is a $\mathbb{Q}$-scheme, then:
              \begin{enumerate}
                  \item The weak $\Theta$-stratification is a
                  $\Theta$-stratification;
                  \item The semistable locus
                  $\GBun^{\pol}_{g,n,p}(N,d)^{\nu_{\epsilon}\dash ss}$ admits a
                  proper and flat good moduli space relative to
                  $\prestable^{\pol}_{g,n,p} \times (\BGL_N)^p$;
                  \item On the center $\cZ$ of each $\nu_\epsilon$-unstable stratum,
                  for all rational values of $\epsilon_0$ with $0<\epsilon_0 \ll 1$, 
                  the numerical invariant  $-\epsilon_0^2 \wt(\cL_{\Cor,\epsilon_0^{-1}}) /
                  \sqrt{b}$ defines a $\Theta$-stratification of $\cZ$ whose centers have
                  proper good moduli spaces relative to
                  $\prestable^{\pol}_{g,n,p} \times (\BGL_N)^p$. In particular $\cZ$ is cohomologically proper over
                  $\prestable^{\pol}_{g,n,p} \times (\BGL_N)^p$ in the sense of
                  \cite[Defn. 2.4.1]{halpernleistner2019mapping}.
              \end{enumerate}
    \end{enumerate}
\end{thm}

\begin{proof}
    Since the polynomial sequence of line bundles $\cL_{\Cor, m}$ is relatively
    ample for the proper morphism $\phi: \GBun^{\pol}_{g,n,p}(N,d) \to
    \Pur_{g,n,p}^{\pol}$ (\Cref{T:projective_morphism}), and the numerical
    invariant $\nu$ is strictly $\Theta$-monotone (\Cref{prop: monotonicity for
    nu on marked sheaves}) and defines a weak $\Theta$-stratification of $\Pur_{g,n,p}^{\pol}$ (\Cref{thm: weak theta stratification pure marked
    sheaves}), it follows from \cite[Prop. 5.5.3]{halpernleistner2018structure} that $\nu_\epsilon$ is
    strictly $\Theta$-monotone and defines a $\Theta$-stratification of
    $\GBun^{\pol}_{g,n,p}(N,d)$. To
    prove the remaining assertions, we replace $\prestable^{\pol}_{g,n,p}$ with
    one of its connected components. Then there is a fixed rank $r$ such that
    the morphism $\phi: \GBun^{\pol}_{g,n,p}(N,d) \to \Pur_{g,n,p}^{\pol}$
    factors through the open and closed substack $\Pur_{g,n,p}^{\pol, r,d}
    \subset \Pur_{g,n,p}^{\pol}$.

    For Part (1), fix $\gamma(\epsilon) \in \mathbb{R}[\epsilon]$ with constant term
    $\gamma_0 \in \mathbb{R}$. By
    \cite[Prop.~5.5.3(3)]{halpernleistner2018structure}, the image of the open
    stratum $\GBun^{\pol}_{g,n,p}(N,d)_{\nu_\epsilon \leq \gamma(\epsilon)}$
    under $\phi$ is contained in $\left(\Pur_{g,n,p}^{\pol, r,
    d}\right)_{\nu \leq \gamma_0}$. Since $\phi$ is of finite type and $\left(\Pur_{g,n,p}^{\pol, r, d}\right)_{\nu \leq
    \gamma_0}$ is of finite type over $\prestable^{\pol}_{g,n,p} \times
    \SVect^p$ (\Cref{thm: weak theta stratification pure marked sheaves}), the
    same holds for $\GBun^{\pol}_{g,n,p}(N,d)_{\nu_\epsilon \leq \gamma(\epsilon)}$.

    For Part (2), the well-ordered relative weak $\Theta$-stratification is a
    $\Theta$-stratification by \cite[Cor. 2.1.9]{halpernleistner2018structure}.
    $\GBun^{\pol}_{g,n,p}(N,d) \to \prestable^{\pol}_{g,n,p} \times (\BGL_N)^p$
    satisfies the existence part of the valuative criterion for properness by
    \Cref{prop: affineness of pur and valuative criterion for properness} and
    \Cref{L: valuative criterion over torsion-free moduli}, so
    $\GBun^{\pol}_{g,n,p}(N,d)^{\nu_{\epsilon}\dash ss}$ admits a proper good
    moduli space relative to $\prestable^{\pol}_{g,n,p} \times (\BGL_N)^p$ by
    \cite[Thm.~5.5.10]{halpernleistner2018structure}. The relative moduli space
    morphism is flat by \Cref{L:cBun algebraic}(c), because if $\cX \to S$ is a flat morphism
    to a space $Y$ and $\cX \to X$ is a good moduli space morphism, then the
    universal morphism $X \to Y$ is also flat \cite[Thm. 4.16(ix)]{alper-good-moduli}.

    According to \cite[Thm.~5.5.10(3)]{halpernleistner2018structure}, there are
    two kinds of strata in $\GBun^{\pol}_{g,n,p}(N,d)$: (A) those that lie over an
    unstable stratum of $\Pur_{g,n,p}^{\pol, r, d}$; and (B) those that lie over the
    semistable locus of $\Pur_{g,n,p}^{\pol, r, d}$. For strata of type (A), the centers are the preimages of the centers of the strata of
    $\Pur_{g,n,p}^{\pol, r, d}$ under the induced morphism
    \[\Grad_{\prestable^{\pol}_{g,n,p} \times
    \left(\GL_N\right)^p} \left( \GBun^{\pol}_{g,n,p}(N,d)\right) \to \Grad_{\prestable^{\pol}_{g,n,p} \times
    \left(\GL_N\right)^p} \left( \Pur_{g,n,p}^{\pol, r, d}\right),\]
    which are quasi-compact closed substacks of
    \[\left(\Grad_{\prestable^{\pol}_{g,n,p} \times
    \left(\GL_N\right)^p} \left( \Pur_{g,n,p}^{\pol, r, d}\right)\right)
    \times_{\Pur_{g,n,p}^{\pol, r, d}} \GBun^{\pol}_{g,n,p}(N,d).\]
    Therefore each center of type (A) is projective over the corresponding
    center of $\Pur_{g,n,p}^{\pol, r, d}$, and $\epsilon^2
    \wt(\cL_{\Cor,\epsilon^{-1}})$ is a relatively ample class. It follows that
    after substituting a sufficiently small rational $\epsilon$, the numerical
    invariant $- \epsilon^2 \wt(\cL_{\Cor,\epsilon^{-1}}) / \sqrt{b}$ defines a
    $\Theta$-stratification of the center of each stratum of
    $\GBun^{\pol}_{g,n,p}(N,d)$ of type (A) whose centers admit good moduli
    spaces relative to $\prestable^{\pol}_{g,n,p} \times
    \left(\GL_N\right)^p$ that are projective over the relative good moduli
    spaces for the centers of $\Pur_{g,n,p}^{\pol, r, d}$, and hence projective
    over $\prestable^{\pol}_{g,n,p} \times
    \left(\GL_N\right)^p$.

    According to \cite[Thm.~5.5.10(3)]{halpernleistner2018structure}, the strata of type (B) in $\GBun^{\pol}_{g,n,p}(N,d)$ are strata for the formal
    numerical invariant $- \epsilon^2 \wt(\cL_{\Cor,\epsilon^{-1}}) / \sqrt{b}$ on the open
    substack $\cY:=\phi^{-1}(\Pur_{g,n,p}^{\pol, r, d}) \subset
    \GBun^{\pol}_{g,n,p}(N,d)$.
    If one substitutes a small positive rational value
    $\epsilon_0$ for $\epsilon$, then $- \epsilon_0^2
    \wt(\cL_{\Cor,\epsilon_0^{-1}})$ is relatively ample for the projective morphism $\phi :
    \cY \to (\Pur_{g,n,p}^{\pol,r,d})^{\rm{ss}}$, and the target admits a good
    moduli space, so the real-valued numerical invariant $- \epsilon_0^2
    \wt(\cL_{\Cor,\epsilon_0^{-1}}) / \sqrt{b}$ defines a $\Theta$ stratification
    whose centers admit relative good moduli spaces that are projective over the
    good moduli space of $(\Pur_{g,n,p}^{\pol,r,d})^{\rm{ss}}$. This
    $\Theta$-stratification refines the $\Theta$-stratification of $\cY$ from
    the formal numerical invariant $\nu_\epsilon$, and hence it endows the centers for the
    formal numerical invariant with $\Theta$-stratifications whose centers are
    centers of the $\epsilon_0$ stratification. This phenomenon will be discussed in more detail in \cite{hl_ibanez}.
\end{proof}

\begin{notn} \label{notn: gieseker bundles over stable curves} We denote by
    $\StMap_{g,n,p}(\BGL_N,d)$ the fiber product
    \[\begin{tikzcd} \StMap_{g,n,p}(\BGL_N,d)   \ar[r] \ar[d] &
            \GBun_{g,n,p}^{\pol}(N,d) \ar[d, "c"]  \\
            \StMap_{g,n+p} \ar[r] & \prestable^{\pol}_{g,n,p}, \end{tikzcd} \]
        where the right most vertical morphism $c: \GBun_{g,n,p}^{\pol}(N,d) \to
        \prestable^{\pol}_{g,n,p}$ is the forgetful morphism as in \Cref{notn:
        forgetful morphism c}, and $\StMap_{g,n+p}$ is the stack of stable
        $(n+p)$-pointed curves equipped with the morphism to
        $\prestable^{\pol}_{g,n,p}$ induced by the log-canonical bundle.
\end{notn}

We equip $\StMap_{g,n,p}(\BGL_N,d)$ with the formal numerical invariant
$\nu_\epsilon$ pulled back from $\GBun_{g,n,p}^{\pol}(N,d)$. The following is an
immediate consequence of \Cref{thm: theta stratification for Gieseker
bundles}.
\begin{coroll} \label{coroll: theta stratification for Gieseker bundles over
    stable curves} Fix a pair of integer $N$ and $d$. The $\nu_{\epsilon}$
    defines a weak $\Theta$-stratification of $\StMap_{g,n,p}(\BGL_N,d)$
    relative to $\StMap_{g,n+p} \times (\BGL_N)^p$ satisfying the following:
    \begin{enumerate}
        \item Every open $\StMap_{g,n,p}(\BGL_N,d)_{\nu_\epsilon \leq \gamma}$
        is of finite type over $\StMap_{g,n+p} \times (\BGL_N)^p$. In
        particular, the relative weak $\Theta$-stratification is well-ordered.
        \item If the base scheme $S$ is a $\mathbb{Q}$-scheme, then
              \begin{enumerate}
                  \item The weak $\Theta$-stratification is a
                  $\Theta$-stratification;
                  \item The semistable locus
                  $\StMap_{g,n,p}(\BGL_N,d)^{\nu_{\epsilon}\dash ss}$ admits a
                  proper and flat good moduli space relative to $\StMap_{g,n+p}
                  \times (\BGL_N)^p$;
                  \item The center of every stratum
                is cohomologically proper relative to $\StMap_{g,n+p} \times
                (\BGL_N)^p$.\qed
              \end{enumerate}
    \end{enumerate}
\end{coroll}

\begin{remark}
    If $p=n=0$ and $S$ is a $\mathbb{Q}$-scheme, then \Cref{coroll: theta
    stratification for Gieseker bundles over stable curves} shows that there
    exists a proper and flat relative good moduli space of semistable pure
    Gieseker bundles of rank $N$ and degree $d$ over the stack
    $\overline{\cM}_{g}$ of stable curves. This recovers the moduli space
    constructed by Schmitt in \cite{schmitt-hilbert-compactification} using GIT.
\end{remark}

\section{Quantum operations on \texorpdfstring{$\Rep(\GL_N)$}{Rep(GLN)}} In this
section, we prove \Cref{T:admissibility}, which asserts that the Fourier-Mukai
transform $\transfer^a_{g,n,p}(N,d)$ used to define $K$-theoretic Gromov-Witten
invariants preserves perfect complexes. In \Cref{subsection: explicit rank 1 and
genus 0} we compute the corresponding quantum operations explictily when the
rank $N=1$ and the genus $g=0$. In \Cref{subsection: gluing and forgetful
morphisms} we construct gluing morphisms for the stacks of pure Gieseker
bundles, and explain the compatibility of the level line bundle $\cL_{\lev}$
with the gluing morphisms.

\subsection{The level line bundle}

\begin{defn}[Level line bundle] \label{defn: level line bundle} Let $(\varphi:
    \widetilde{C} \to C, \sigma_{\bullet}^{\pm}, E)$ be the universal pure
    Gieseker bundle over $\GBun_{g,n,p}(N)$, and let $\pi: \widetilde{C} \to
    \GBun_{g,n,p}(N)$ denote the structure morphism. We define the line bundle
    $\cL_{\lev}$ on $\GBun_{g,n,p}(N)$ by
    \[ \cL_{\lev}:= \det\left(R\pi_\ast\left(E \oplus E \otimes
    \omega_{\pi}\right)\right)^{\vee}.\]
\end{defn}

It turns out that $\cL_{\lev}$ is the pullback of a line bundle on
$\Pur_{g,n,p}(N)$ via the morphism $\Phi$. We first need to introduce some
notation in order to provide a formula for $\cL_{\lev}$.

\begin{defn} \label{defn: evaluation line bundle on moduli of pure marked
    sheaves} Let $(C, \overline{\sigma}_{\bullet}^{\pm}, F,
    (\sigma_{\bullet}^+)^*(F) \to W_{\bullet})$ be the universal pure marked
    sheaf over $\Pur_{g,n,p}(N)$. Let $K$ denote the kernel of the induced
    morphism $F \twoheadrightarrow \bigoplus_{i=1}^p (\sigma_i^+)_*(W_i)$, which
    is a family of pure sheaves on $C \to \Pur_{g,n,p}(N)$. For each $1 \leq i
    \leq p$, we define the following line bundles on $\Pur_{g,n,p}(N)$:
    \[\cL_{e_i} := \det(W_i), \; \; \; \cL_e := \bigotimes_{i=1}^p \cL_{e_i}, \;
    \;\; \cL_{k_i} = \det((\sigma_i^+)^*(K)), \; \; \text{and} \; \;\cL_k =
    \bigotimes_{i=1}^p \cL_{k_i}.\]
\end{defn}

\begin{lemma} \label{lemma: formula for level line bundle} The line bundle
    $\cL_{\lev}$ from \Cref{defn: level line bundle} is isomorphic to the
    pullback $\Phi^*(\cL_{\deg} \otimes \cL_e^{\vee} \otimes \cL_k^{\vee})$ via
    the morphism $\Phi : \GBun_{g,n,p}(N) \to \Pur_{g,n,p}(N)$.
\end{lemma}
\begin{proof}
    We use the same terminology for the universal family over $\GBun_{g,n,p}(N)$
    as in \Cref{defn: level line bundle}. We denote by $\overline{\pi}: C \to
    \GBun_{g,n,p}(N)$ the structure morphism for the target curve, and set
    $\overline{\sigma}_{\bullet}^{\pm} := \varphi \circ \sigma_{\bullet}^{\pm}$.
    We have that $\cL_{\lev}^{\vee}$ is given by
    \begin{gather*}
        \det\left(R\pi_*\left( E \oplus E \otimes \omega_{\pi}\right)\right)= \det\left(R\overline{\pi}_*\left( R\varphi_*\left(E\right) \oplus R\varphi_*\left(E\otimes \omega_{\pi}\right) \right)\right).
    \end{gather*}

    Let $(C, \overline{\sigma}_{\bullet}^{\pm}, F,
    (\overline{\sigma}_{\bullet}^+)^*(F) \to W_i)$ be the pullback of the
    universal marked sheaf via the morphism $\Phi$. We denote by $K$ the
    corresponding kernel of the morphism $F \to \bigoplus_{i=1}^p
    (\sigma_i^+)_*(W_i)$. Then we have
    \begin{itemize}
        \item $\Phi^*(\cL_{\deg}^{\vee}) = \det\left(R\overline{\pi}_*(K \oplus
        K \otimes \omega_{\overline{\pi}})\right)$;
        \item $\Phi^*(\cL_e) = \det(\bigoplus_{i=1}^p W_i) = \det
        \left(R\overline{\pi}_*\left( \bigoplus_{i=1}^p
        (\overline{\sigma}_i^+)_*(W_i)\right)\right)$; and
        \item $\Phi^*(\cL_k) = \det( \bigoplus_{i=1}^p
        (\overline{\sigma}_i^+)^*(K)) = \det\left(R\overline{\pi}_*\left(
        \bigoplus_{i=1}^p (\overline{\sigma}_i^+)_*
        (\overline{\sigma}_i^+)^*(K))\right)\right)$.
    \end{itemize}
    Hence, in order to prove the desired identification $\cL^{\vee}_{\lev} \cong
    \Phi^*(\cL_{\deg}^{\vee} \otimes \cL_e \otimes \cL_k)$, it suffices to show
    the following equality in $K$-theory of perfect complexes on $C$
    \[
        \left[R\varphi_*\left(E\right)\right] + \left[R\varphi_*\left(E\otimes \omega_{\pi}\right)\right] = [K] + \left[K\otimes \omega_{\overline{\pi}}\right] + \left[\bigoplus_{i=1}^p (\overline{\sigma}_i^+)_*(W_i)\right] + \left[\bigoplus_{i=1}^p(\overline{\sigma}_i^+)_* (\overline{\sigma}_i^+)^*(K))\right].
    \]
    Note that $R^1\varphi_*\left(E\right)=0$ by \Cref{lemma: pushforward lemma},
    and by the definition of the morphism $\Phi$ the pushforward $\varphi_*(E)=
    F$ fits into a short exact sequence
    \[ 0 \to K \to \varphi_*(E) \to \bigoplus_{i=1}^p
    (\overline{\sigma}_i^+)_*(W_i) \to 0. \] Therefore we have $[R\varphi_*(E)]
    = [K] + \left[\bigoplus_{i=1}^p (\overline{\sigma}_i^+)_*(W_i)\right]$, and
    the desired equality in $K$-theory reduces to
    \[  \left[R\varphi_*\left(E\otimes \omega_{\pi}\right)\right] =
    \left[K\otimes \omega_{\overline{\pi}}\right]+
    \left[\bigoplus_{i=1}^p(\overline{\sigma}_i^+)_*
    (\overline{\sigma}_i^+)^*(K))\right].\] Since $\varphi: (\widetilde{C},
    \sigma_{\bullet}^{\pm}) \to (C, \overline{\sigma}_{\bullet}^{\pm})$ is
    log-crepant, we have an isomorphism $\omega_{\pi}(\sum
    \sigma_{\bullet}^{\pm}) \cong \varphi^*(\omega_{\overline{\pi}}(\sum
    \overline{\sigma}_{\bullet}^{\pm}))$. Furthermore, by \Cref{lemma: control
    of bubblings} we know that $\varphi$ is an isomorphism over the negative
    markings $\overline{\sigma}_{\bullet}^-$ of $C$, and hence
    $\omega_{\pi}(-\sum_{i=1}^p \sigma_{i}^+) \cong
    \varphi^*(\omega_{\overline{\pi}}(-\sum_{i=1}^p \overline{\sigma}_{i}^+))$.
    Using the projection formula, we get
    \[ R\varphi_*\left(E\otimes \omega_{\pi}\right) =
    R\varphi_*\left(E\left(-\sum_{i=1}^p \sigma_i^+\right)\right) \otimes
    \omega_{\overline{\pi}}\left(\sum_{i=1}^p \overline{\sigma}_i^+\right).\] By
    \Cref{lemma: pushforward lemma}, we have $R^1\varphi_*\left(E(-\sum_{i=1}^p
    \sigma_i^+)\right) = 0$, and by definition of the morphism $\Phi$ we have
    $\varphi_*(E(-\sum_{i=1}^p \sigma_i^+)) = K$. Hence, we may rewrite the
    equality above as $R\varphi_*\left(E\otimes \omega_{\pi}\right) = K \otimes
    \omega_{\overline{\pi}}(\sum_{i=1}^p\overline{\sigma}_i^+)$. The short exact
    sequence
    \[ 0 \to K\otimes \omega_{\overline{\pi}} \to K \otimes
    \omega_{\overline{\pi}}(\sum_{i=1}^p\overline{\sigma}_i^+) \to
    \bigoplus_{i=1}^p (\overline{\sigma}_i^+)_*\left(
    (\overline{\sigma}_i^+)^*(K) \otimes
    (\overline{\sigma}_i^+)^*(\omega_{\overline{\pi}}(\overline{\sigma}_i^+))\right)
    \to 0\] and the residue isomorphism
    $(\overline{\sigma}_i^+)^*(\omega_{\overline{\pi}}(\overline{\sigma}_i^+))
    \cong \cO_{\GBun_{g,n,p}(N)}$ imply the desired equality
    \[ \left[R\varphi_*\left(E\otimes \omega_{\pi}\right)\right] = [K \otimes
    \omega_{\overline{\pi}}(\sum_{i=1}^p\overline{\sigma}_i^+)] = \left[K\otimes
    \omega_{\overline{\pi}}\right]+
    \left[\bigoplus_{i=1}^p(\overline{\sigma}_i^+)_*
    (\overline{\sigma}_i^+)^*(K))\right].\qedhere\]
\end{proof}

\subsection{Admissibility of tautological classes}

In this subsection, we assume that the base scheme $S$ is a $\mathbb{Q}$-scheme. We fix a degree $d$ and consider the stack $\StMap_{g,n,p}(\BGL_N, d)$ equipped with its corresponding morphism $c\times \ev_+: \StMap_{g,n,p}(\BGL_N, d) \to \StMap_{g, n+p} \times (\BGL_N)^p$.
Let $\tot_\gamma : \cZ_\gamma \to \StMap_{g,n,p}(\BGL_N,d)$ denote the centers of the $\Theta$-stratification from \Cref{coroll: theta stratification for Gieseker bundles over stable curves}, which are indexed by elements $\gamma \in \mathbb{R}[\epsilon]$.
\begin{defn}[{\cite[Sect. 2.1]{halpernleistner2016equivariantverlindeformulamoduli}}] \label{defn: admissible}
    Given a line bundle $\cL$ on $\StMap_{g,n,p}(\BGL_N,d)$, we say that a perfect complex $F \in \Perf(\StMap_{g,n,p}(\BGL_N,d))$ is $\cL$-admissible if for any $m>0$, the $\mathbb{G}_m$-weights appearing in
\[
\tot_\gamma^\ast \left(\cL \otimes F^{\otimes m} \right) \otimes \det((\tot_\gamma^\ast(\bL_{c \times \ev_+}))^+)^{-1}
\]
are negative for all but finitely many indices $\gamma$. Here $(-)^+$ denotes the direct summand that has positive weight with respect to the central $\bG_m$-automorphism on $\cZ_\gamma$.
\end{defn}

\begin{prop} \label{P:local_cohomology_vanishing}
    For any $\cL$-admissible $F \in \Perf(\StMap_{g,n,p}(\BGL_N,d))$ as in \Cref{defn: admissible}, we have that $R(c \times \ev_+)_\ast(\cL \otimes F) \in \Dqc(\StMap_{g,n+p} \times (\BGL_N)^p)$ is perfect, and its formation commutes with flat base change.
\end{prop}
\begin{proof}
    This follows from \cite[Lem.~3.24\&3.25]{halpernleistner2025categoricalperspectivenonabelianlocalization} and the fact that the centers of the stratification are cohomologically proper over $\StMap_{g,n+p} \times (\BGL_N)^p$ by \Cref{coroll: theta stratification for Gieseker bundles over stable curves}.
\end{proof}

Before we investigate admissibility on $\StMap_{g,n,p}(\GL_N,d)$, we will analyze its $\Theta$-stratification. 
\begin{remark} \label{remark: relation between strata in gauged maps and pure sheaves}
    Note that, if a pure Gieseker bundle has underlying unstable pure marked sheaf, then its HN filtration is the unique lift of the HN filtration of the underlying pure marked sheaf \cite[Thm.~5.5.10(3)]{halpernleistner2018structure}. Furthermore, only finitely many strata of $\StMap_{g,n,p}(\GL_N,d)$ lie above each stratum of $\Pur^{\pol}_{g,n,p}(N)$.
\end{remark} 
In view of \Cref{remark: relation between strata in gauged maps and pure sheaves}, it will suffice for our purposes to index the stratification of $\Pur^{\pol}_{g,n,p}(N)$. Since we are working with the base-change $\StMap_{g,n,p}(\GL_N,d)$ as in \Cref{notn: gieseker bundles over stable curves}, we replace $\Pur^{\pol}_{g,n,p}(N)$ with the following fiber product
    \[\begin{tikzcd}
   \overline{\Pur}_{g,n,p}(N,d)   \ar[r] \ar[d] & \Pur^{\pol, N \cdot h, d}(N) \ar[d]  \\
   \StMap_{g,n+p} \ar[r] & \prestable^{\pol}_{g,n,p},
\end{tikzcd} \]
equipped with its corresponding $\Theta$-stratification. Here we set $h:= 2g-2+n+p$ to be the degree of the log-canonical polarization of any stable curve in $\StMap_{g,n+p}$, so that $N \cdot h$ is the rank with respect to the log-canonical bundle of every sheaf that has rank $N$ at every generic point on any such stable curve.

According to \Cref{lemma: technical optimization lemma hn boundednes marked sheaves}, the strata of $\overline{\Pur}_{g,n,p}(N,d)$ are labeled by a collection of integer valued degrees and ranks $(d_0,\ldots,d_q; r_0,\ldots,r_q)$, and an index $0\leq j\leq q$ such that
\begin{enumerate}
    \item $d_0 + \cdots + d_q = d$,\\
    \item $0<r_i$ for all $i$ and $r_0+\cdots+r_q = N \cdot h$,\\
    \item $\frac{d_0}{r_0} < \ldots < \frac{d_{j-1}}{r_{j-1}} < \max(\frac{d_j}{r_j}, \frac{d}{N \cdot h}) < \frac{d_{j+1}}{r_{j+1}} <\cdots < \frac{d_q}{r_q}$
\end{enumerate}
This data indexes the stratum of objects whose HN filtration has $\deg(\overline{G}_i) = d_i$ and $\rk(\overline{G}_i) = r_i$ in the notation of \Cref{lemma: technical optimization lemma hn boundednes marked sheaves}.

\begin{lemma}\label{L:recover_degrees_from_weights}
    The indexing data for the stratum is uniquely determined by the vector
    \[
    \vec{w} = (w_1,\ldots,w_{N\cdot h}) = \left( \frac{d_0}{r_0}, \ldots, \frac{d_0}{r_0}, \frac{d_1}{r_1},\ldots \right),
    \]
    where the $i^{th}$ term of the inequality (3) above is repeated $r_i$ times, and the smallest index $j'$ such that $w_{j'} = \max(\frac{d_j}{r_j},\frac{d}{Nh})$.
\end{lemma}
\begin{proof}
    The value of $q$ and $r_0,\ldots,r_q$ can be uniquely recovered from $\vec{w}$ by counting the number of unique entries and the number of repetitions of each entry. Then $d_0,\ldots,d_{j-1},d_{j+1},\ldots,d_q$ can be recovered by multiplying the corresponding entry of $\vec{w}$ by $r_i$. We have then $d_j = d - \sum_{i \neq j} d_i$, and the index $j$ is determined by $j'$.    
\end{proof}

Therefore, we will index the strata of $\overline{\Pur}_{g,n,p}(N,d)$ by vectors $\vec{w} \in \frac{1}{M} \bZ^{N\cdot h} \subset \bQ^{N\cdot h}$, where $M={\rm{lcm}}(1,\ldots,N \cdot h)$, along with a choice of index $1 \leq j \leq N\cdot h$ such that $w_1 \leq w_2 \leq \cdots \leq w_n$, and $w_j>w_{j-1}$ if $j>1$. The last condition is to avoid duplicate indexing. Such a vector only encodes a possible stratum if the values of $d_i$ determined by $\vec{w}$ as in the proof of \Cref{L:recover_degrees_from_weights} are integers.

\begin{thm}\label{T:admissibility}
    Suppose that $S$ is a $\mathbb{Q}$-scheme. For every $a>0$ and $b,d \in \bZ$, and for every $V \in \Perf((\BGL_N)^{n})$, the perfect complex $\ev_-^\ast(V)$ is $\cL_{\lev}^a \otimes \cL_{\rk}^b$-admissible on $\StMap_{g,n,p}(\BGL_N, d)$. It follows that $\transfer^a_{g,n,p}(N,d) : \Dqc((\BGL_N)^n) \to \Dqc(\StMap_{g,n+p} \times (\BGL_N)^p)$ takes perfect complexes to perfect complexes.
\end{thm}
\begin{proof}
    Consider the center of a stratum $\tot_\gamma : \cZ_\gamma \to \StMap_{g,n,p}(N,d)$, where $\gamma$ is an abstract index. By \Cref{remark: relation between strata in gauged maps and pure sheaves}, only finitely many $\cZ_{\gamma}$ lie over the semistable locus of $\overline{\Pur}_{g,n,p}(N,d)$. For the purposes of admissibility, we may ignore those and only consider centers $\cZ_\gamma$ lying over the unstable locus of $\overline{\Pur}_{g,n,p}(N,d)$. Again by \Cref{remark: relation between strata in gauged maps and pure sheaves}, any such $\cZ_{\gamma}$ lies over the center of a stratum of $\overline{\Pur}_{g,n,p}(N,d)$, which is indexed by a set of rational numbers with bounded denominators $w_1 \leq \cdots \leq w_{N\cdot h}$ and a choice of index $j$ with $w_j > w_{j-1}$. The induced assignment of indexing sets $\gamma \mapsto \vec{w}$ is finite-to-one. The center of the stratum in $\overline{\Pur}_{g,n,p}(N,d)$ classifies graded marked sheaves $(C,\sigma^\pm_\bullet,F = \bigoplus_v F_v,(\sigma^+_\bullet)^\ast (F_0) \twoheadrightarrow W_\bullet)$. We will let $K_w = F_w$ for $w \neq 0$ and $K_0 = \ker(F_0 \to \bigoplus_i (\sigma^+_i)_\ast(W_i))$.
    
    Let us now re-index, so $K = K_0 \oplus \cdots \oplus K_q$ is the decomposition of $K$ given by the nonvanishing summands $K_v$ and the summand $K_i$ is given weight $v_i = M \cdot (w'_i - \frac{d}{N\cdot h})$, where we use $w'_i$ to denote the $i^{th}$ unique entry of $\vec{w}$ for $i=0,\ldots,q$. We will again denote the index of the weight $v$ corresponding to the distinguished entry of $\vec{w}$ by $j$. \Cref{lemma: technical optimization lemma hn boundednes marked sheaves} implies that $w'_i = \deg(K_i) / \rk(K_i)$ for $i \neq j$ and $w'_j = \max \left(\frac{\deg(K_i)}{\rk(K_i)}, \frac{d}{N\cdot h}\right)$. By \Cref{lemma: formula for level line bundle}, we have that $\cL_{\lev}$ is the pullback of the line bundle $\cL_{\deg} \otimes \cL_{e}^{\vee} \otimes \cL_{k}^{\vee}$ on $\overline{\Pur}_{g,n,p}(N,d)$. Similarly, the complexes $\cL_{\rk}$, and $\ev_-^\ast(V)$ are pulled back from corresponding complexes on $\Pur^{\pol,h}_{g,n,p}(N,d)$, which we denote by the same names. To compute the weights of these complexes restricted to the center $\cZ_\gamma$, it suffices to compute their weights on the centers of the stratum of $\overline{\Pur}_{g,n,p}(N,d)$ corresponding to $\vec{w}$.
    \begin{itemize}
        \item We have
        \[
            \wt(\tot_\gamma^\ast(\cL_{\rk})) = \sum_i \rk(K_i) v_i = M (\vec{1} \cdot \vec{w} - d) = M \vec{1} \cdot (\vec{w} - \frac{d}{Nh} \vec{1}),
        \]
        where we are using the standard inner product in $\mathbb{R}^{N \cdot h}$ and we denote by $\vec{1}$ the vector with all entries equal to $1$. Here we have used the fact that the entry $w_i'$ is repeated $\rk(K_i)$ times in $\vec{w}$, and $\sum_i \rk(K_i) = N\cdot h$.\smallskip

        \item Next, we compute
        \begin{align*}\wt(\tot_\gamma^\ast(\cL_{\deg}^{N h} \otimes \cL_{\rk}^{2d})) &= - \sum_i 2 (Nh \deg(K_i) - d \rk(K_i))  v_i \\
        &= - \sum_i 2 M N h\rk(K_i) \left( \frac{\deg(K_i)}{\rk(K_i)} - \frac{d}{Nh} \right) \left(w'_i-\frac{d}{Nh} \right)\\
        &= - 2 M N h \lVert \vec{w} - \frac{d}{Nh} \vec{1} \rVert^2.
        \end{align*}
        The last equality again uses the fact that $w_i'$ is repeated $\rk(K_i)$ times and that $w_i' = \deg(K_i)/\rk(K_i)$ if $i \neq j$. For the $j^{th}$ term in the sum, either $w_j' - d/(Nh)$ vanishes or $w_j' = \deg(K_j)/\rk(K_j)$, so we can replace $\deg(K_j)/\rk(K_j)$ with $w_j'$ without changing the sum.\smallskip

        \item Combining the previous two computations gives
        \[\wt(\tot_\gamma^\ast(\cL_{\deg})) = - 2 M \vec{w} \cdot (\vec{w} - \frac{d}{Nh} \vec{1}).\]

          \item We have $\wt(\tot_{\gamma}(\cL_{e}))=0$, since $\cL_{e}$ is pulled back from $(\BGL_N)^p$ via the morphism $\ev_+$ and the $\Theta$-stratification is relative to the base $\StMap_{g,n+p} \times (\BGL_N)^p$.\smallskip

        \item We may break $V$ into a direct sum of irreducible representations and assume that it is of the form $V = V_1 \otimes \cdots \otimes V_n$, where each $V_i$ is given by a Schur functor of the standard representation of $\GL_N$. For each section $\sigma^-_i$, let $r^-_{i,\ell} = \rk((\sigma_i^-)^\ast(K_\ell))$ for $\ell = 0,\ldots,q$. These numbers depend on $i$, because $K_\ell$ is not locally free and can even vanish on certain components of $C$. However, we have $r_{i,\ell}^- \geq 0$ and $r_{i,0}^-+\cdots r_{i,q}^-=N$ for all $i=1,\ldots,n$. Then the weights appearing in $\ev_-^\ast(V)$ are bounded above by some fixed positive constant multiple of $|\max_i (v_i)| +  |\min_i(v_i)|$. 
        
        \item Finally, for each section $\sigma_i^+$ let us set $r_{i,\ell}^+ = \rk((\sigma_i^+)^\ast(K_\ell))$ for $\ell = 0,\ldots,q$. Then we have $\wt(\tot_{\gamma}^*(\cL_{k})) = \sum_{i=1}^p \sum_{\ell =0}^q r_{i,\ell}^+ v_i$. Again, since we have $r_{i,\ell}^- \geq 0$ and $r_{i,0}^-+\cdots r_{i,q}^-=N$ for all $i=1,\ldots,n$, it follows that $\wt(\tot_{\gamma}^*(\cL_{k}))$ is bounded between $nN\max_i (v_i)$  and $nN\min_i (v_i)$.\smallskip
       
    \end{itemize}

    Because $c \times \ev_+: \StMap_{g,n,p}(\BGL_N,d) \to \StMap_{g,n+p} \times (\BGL_N)^p$ is smooth, the weight of $\det((\tot_\gamma^\ast(\bL_{c \times \ev_+}))^+)$ is automatically positive. Therefore, to verify that a complex $F \in \Perf(\StMap_{g,n,p}(\BGL_N,d))$ is $\cL^a_{\det}\otimes \cL_{\rk}^b$-admissible, it suffices to show that for all $m >0$, the weights appearing in $\tot_\gamma^\ast(\cL^a_{\det} \otimes \cL_{\rk}^b \otimes F^{\otimes m})$ are negative for all but finitely many $\gamma$. The key observation above is that for $a>0$ and $b$ arbitrary, the weights appearing in 
    \[\cL_{\lev}^a \otimes \cL_{\rk}^b \otimes \ev_-^\ast(V) = \Phi^*(\cL_{\deg}^a \otimes \cL_{e}^{-a} \otimes \cL_{k}^{-a} \otimes \cL_{\rk}^b \otimes \ev_-^\ast(V))\]
    are $-2 M a \lVert \vec{w} \rVert^2 + O(\lVert \vec{w} \rVert)$ so only finitely many values of $\vec{w}$ allow for a nonnegative weight. This concludes the proof of admissibility. The fact that $\transfer^a_{g,n,p}(d)$ preserves perfect complexes is now a consequence of \Cref{P:local_cohomology_vanishing}.
\end{proof}

\subsection{The case \texorpdfstring{$N=1$}{N=1} and \texorpdfstring{$g=0$}{g
equals 0}} \label{subsection: explicit rank 1 and genus 0}

In this subsection, we specialize to the case when our base Noetherian scheme
$S$ is of the form $\Spec(k)$ for some field $k$ of characteristic $0$. We
investigate quantum operations in the case when the target is $B\mathbb{G}_m$
and the genus is $0$.

\subsubsection*{The two-point operation}

Let $\StMap_{0,1,1}(B\mathbb{G}_m,d) \subset \GBun_{0,1,1}(1,d)$ denote the
preimage via $c$ of the open substack $B\bG_m \subset \prestable_{0,1,1}$
parameterizing the smooth $(1,1)$-marked curve $(\bP^1_k,x^-,x^+)$ with its
$\bG_m$ automorphism group. Likewise, we let $\Pur_{0,1,1}^{\circ}(1,d) \subset
\Pur_{0,1,1}(1,d)$ and $\Pur^{\circ}_{0,2,0}(1,d) \subset \Pur_{0,2,0}(1,d)$
denote the open substacks of pure marked sheaves whose underlying marked curve
is isomorphic to $(\bP^1,x^-,x^+)$ and the degree of the pure sheaf kernel is
$d$ (as in \Cref{notn: fixed degree pure equirank marked sheaves}).

The stack $\Pur^{\circ}_{0,2,0}(1,d)$ parameterizes triples $(\pi : C \to T,
x^-,x^+, \cL)$, where $x_\pm : T \to C$ are sections of $\pi$, every geometric
fiber of $\pi$ is isomorphic to $\bP^1$, and $\cL$ is an invertible sheaf on $C$
of degree $d$ on every fiber. (We choose to also denote the two marked points
$x^\pm$ to align with the analysis of $\Pur_{0,1,1}^{\circ}(1,d)$.) To such
data, we can assign the sheaf $L := \pi_\ast(\cO_C(x^+ - x^-))$ on $T$, which is
an invertible sheaf because $\cO_C(x^+-x^-)$ is an invertible sheaf of degree
$0$ on every $T$-fiber. There are isomorphisms
\[
    L \cong (x^+)^\ast(\cO_C(x^+)) \cong (x^-)^\ast(\cO_C(-x^-)),
\]
induced by pulling back the morphisms $\cO(x^+-x^-) \to \cO(x^+)$ and $\cO(-x^-)
\to \cO(x^+-x^-)$ under $x^+$ and $x^-$ respectively. In other words, $L$ is the
normal bundle of the divisor $x^+$, and the conormal bundle of the divisor
$x^-$.

\begin{lemma}\label{L:Pur_{0,2,0}} The morphism $(L,\ev_-) :
    \Pur^{\circ}_{0,2,0}(1,d) \to B\bG_m^2$ is an isomorphism. Under this
    identification, the morphism $B\mathbb{G}_m^2 \cong
    \Pur^{\circ}_{0,2,0}(1,d) \xrightarrow{\ev_+}  B\bG_m$ corresponds to the
    character $(d,1): \mathbb{G}_m^2 \to \mathbb{G}_m$.
\end{lemma}
\begin{proof}
    Consider the $k$-point $f: \Spec(k) \to \Pur^{\circ}_{0,2,0}(1,d)$ given by
    $(\mathbb{P}^1, x_{\pm}, \cO(d \cdot x^+))$. The automorphism group scheme
    of $f$ fits into a short exact sequence
    \[ 1 \to \mathbb{G}_m \to \Aut(f) \to \mathbb{G}_m \to 1\] where the
    quotient morphism $\Aut(f) \twoheadrightarrow \mathbb{G}_m \cong
    \Aut((\mathbb{P}^1, x_{\pm}))$ is induced from $\Pur^{\circ}_{0,2,0}(1,d)
    \to B\mathbb{G}_m \subset \prestable_{0,2,0}$. We get an induced
    representable morphism $h:B\!\Aut(f) \to \Pur^{\circ}_{0,2,0}(1,d)$. The
    base-change of $h$ via the faithfully flat morphism $\Spec(k) \to
    B\mathbb{G}_m \subset \prestable_{0,2,0}$ becomes an isomorphism
    $B\mathbb{G}_m \xrightarrow{\sim} \mathcal{P}ic^d(\mathbb{P}^1)$. We
    conclude that $h$ is an isomorphism. The composition $B\!\Aut(f)
    \xrightarrow{\sim} \Pur^{\circ}_{0,2,0}(1,d) \xrightarrow{(L, \ev_-)}
    B\bG_m^2$ is readily seen to induce an isomorphism of automorphism group
    schemes, where the subgroup $\mathbb{G}_m \hookrightarrow \Aut(f)$ from the
    short exact sequence above corresponds to the cocharacter $(0,1):
    \mathbb{G}_m \to \mathbb{G}_m^2$. This implies that
    $\Pur^{\circ}_{0,2,0}(1,d) \xrightarrow{(L, \ev_-)}  B\bG_m^2$ is an
    isomorphism of gerbes over $\Spec(k)$. The last assertion regarding the
    morphism $\ev_+$ follows from the fact that $\Spec(k) \xrightarrow{f}
    \Pur^{\circ}_{0,2,0}(1,d) \xrightarrow{\ev_+} \mathbb{G}_m$ is given by
    $(x^+)^*(\cO(d \cdot x^+))$. \endnote{Since the subgroup $\mathbb{G}_m
    \subset \Aut(f)$ evidently acts with weight $1$ on the $x^+$-fiber $\Spec(k)
    \xrightarrow{f}\Pur^{\circ}_{0,2,0}(1,d) \xrightarrow{\ev_+} B\mathbb{G}_m$
    of the line bundle, it follows that the composition
        \[B\mathbb{G}_m \xrightarrow{(0,1)} B\mathbb{G}_m^2 \xrightarrow{\sim}
        \Pur^{\circ}_{0,2,0}(1,d) \xrightarrow{\ev_+} B\mathbb{G}_m\]
        corresponds to the identity character $1: \mathbb{G}_m \to
        \mathbb{G}_m$. Note also that the composition $\Spec(k) \xrightarrow{f}
        \Pur^{\circ}_{0,2,0}(1,d) \xrightarrow{\ev_+} \mathbb{G}_m$ is given by
        $(x^+)^*(\cO(d \cdot x^+))$. In view of the identification $L\cong
        (x^+)^*(\cO(x^+))$, we conclude that the composition
        \[B\mathbb{G}_m \xrightarrow{(1,0)} B\mathbb{G}_m^2 \xrightarrow{\sim}
        \Pur^{\circ}_{0,2,0}(1,d) \xrightarrow{\ev_+} B\mathbb{G}_m\]
        corresponds to the character $d: \mathbb{G}_m \to \mathbb{G}_m$.}
\end{proof}

Recall that the morphism $\Psi : \Pur_{0,1,1}^{\circ}(1,d) \to
\Pur_{0,2,0}^{\circ}(1,d)$ sends a $T$-point $(\pi : C \to T, x^-,x^+, s : F
\twoheadrightarrow (x^+)_\ast(E))$ of $\Pur_{0,1,1}^{\circ}(1,d)$ to the tuple
$(\pi : C \to T, x_\pm, \ker(s))$.

\begin{lemma} \label{lemma: identification two points pure with A1 mod G_m
    cubed} There is a commutative diagram
    \[
        \begin{tikzcd}
            \Pur_{0,1,1}^{\circ}(1,d) \ar[d, "\Psi"] \ar[r, "\tau"] & \bA^1_k / \bG_m^3 \ar[d] \\
            \Pur_{0,2,0}^{\circ}(1,d) \ar[r, "{(L, \ev_-)}"] & B\mathbb{G}_m^2,
        \end{tikzcd}
    \]
    where $\tau$ is an isomorphism, $\bG_m^3$ acts linearly on $\bA^1_k$ with
    weights $(d+1,1,-1)$, and $\bA^1_k / \bG_m^3 \to B\mathbb{G}_m^2$ is induced
    from the group homomorphism $\bG_m^3 \to \bG_m^2$ given by projection onto
    the first two factors. Under this identification, the morphism
    $\mathbb{A}^1_k/\mathbb{G}_m^3 \cong \Pur_{0,1,1}^{\circ}(1,d)
    \xrightarrow{\ev_+} B\bG_m$ corresponds to the group homomorphism $\bG_m^3
    \to \bG_m$ that projects onto the third factor.
\end{lemma}
\begin{proof}
    The morphism $\Psi \times \ev_+: \Pur_{0,1,1}^{\circ}(1,d) \to
    \Pur_{0,2,0}^{\circ}(1,d) \times B\mathbb{G}_m \cong B\mathbb{G}^3_m$ is
    affine and of finite type by \Cref{proposition: marked sheaves affine over
    nonmarked sheaves}. More precisely, if we denote by $\pi: (C, x_{\pm}) \to
    \Pur_{0,2,0}^{\circ}(1,d) \times B\mathbb{G}_m \cong B\mathbb{G}^3_m$ the
    universal $(1,1)$-marked curve equipped with a universal pure sheaf $K$
    pulled back from the first factor $\Pur_{0,2,0}^{\circ}(1,d)$, and we denote
    by $E$ the universal line bundle pulled-back from the last factor
    $B\mathbb{G}_m$, then we have an isomorphism of $B\mathbb{G}^3_m$-stacks
    $\Pur_{0,1,1}^{\circ}(1,d) \cong \Hom(E, (x^+)^*(K(x^+)))$. In this case
    $\Hom(E, (x^+)^*(K(x^+))) \cong \mathbb{A}^1_k/\mathbb{G}_m^3$ is the total
    space of a line bundle on $B\mathbb{G}_m^3$, corresponding to a linear
    action of $\mathbb{G}_m^3$ on $\mathbb{A}^1_k$. To conclude, we just need to
    check that the weights of $\mathbb{G}_m^3$ on $\mathbb{A}^1_k$ are as
    claimed in the statement of the lemma, which may be checked directly using
    the identification from \Cref{L:Pur_{0,2,0}}. \endnote{First, we note that
    the last coordinate $\mathbb{G}_m \xrightarrow{(0,0,1)} \mathbb{G}_m^3$ acts
    on $\Hom(E, (x^+)^*(K(x^+)))$ via its natural weight 1 action on $E$, and so
    it acts with weight $-1$ on the corresponding $\mathbb{A}^1_k$. On the other
    hand, the first two coordinates of $\mathbb{G}_m^3$ act via their natural
    action on the corresponding line bundle $(x^+)^*(K(x^+)) = (x^+)^*(K)
    \otimes (x^+)^*(\cO(x^+))$ on $\Pur_{0,2,0}^{\circ}(1,d) \cong
    B\mathbb{G}_m^2$. Note that $(x^+)^*(K)$ corresponds to the morphism $\ev_+:
    \Pur_{0,2,0}^{\circ}(1,d) \to B\mathbb{G}_m$, and hence $\mathbb{G}_m^2$
    acts with weight $(d,1)$ on $(x^+)^*(K)$ by \Cref{L:Pur_{0,2,0}}. A similar
    reasoning as in the proof of \Cref{L:Pur_{0,2,0}} also shows that
    $\mathbb{G}_m^2$ acts with weight $(1,0)$ on $(x^+)^*(\cO(x^+))$, and hence
    we get that the first two coordinates $\mathbb{G}_m^2$ act with weight
    $(d+1,1)$ on $(x^+)^*(K(x^+))$.}
\end{proof}

\begin{lemma} \label{lemma: two points Phi is an iso} The proper morphism $\Phi:
    \StMap_{0,1,1}(B\mathbb{G}_m,d) \to \Pur^{\circ}_{0,1,1}(1,d)$ (as in
    \Cref{lemma: Phis is well-defined}) is an isomorphism, and the following
    diagram commutes
    \[
        \begin{tikzcd}
            \StMap_{0,1,1}(B\mathbb{G}_m,d) \ar[d, "\Phi"] \ar[rr, "{(L,\ev_-,\ev_+)}"] & & B\bG_m^3 \ar[d, "="] \\
            \Pur^{\circ}_{0,1,1}(1,d) \ar[r, "\tau"] & \mathbb{A}^1_k/\mathbb{G}_m^3 \ar[r] & B\mathbb{G}_m^3.
        \end{tikzcd}
    \]
\end{lemma}
\begin{proof}
    The commutativity of the diagram follows from \Cref{lemma: identification
    two points pure with A1 mod G_m cubed}. To show that $\Phi$ is an
    isomorphism, we first claim that $\Phi: \StMap_{0,1,1}(B\mathbb{G}_m,d) \to
    \Pur^{\circ}_{0,1,1}(1,d) \cong \mathbb{A}^1_k/\mathbb{G}_m^3$ is an
    isomorphism away from the origin $0: \Spec(k) \to
    \mathbb{A}^1_k/\mathbb{G}_m^3$. To see this, note that $\Phi$ is an
    isomorphism over the locus of geometric points $(\varphi: (\widetilde{C},
    \sigma_{-}, \sigma_+) \to (\mathbb{P}^1_K, x^-, x^+), E)$ in
    $\StMap_{0,1,1}(B\mathbb{G}_m,d)$ valued in an algebraically closed field
    extension $\widetilde{k} \supset k$ where the underlying source curve
    $\widetilde{C}$ is smooth (thus forcing $\varphi$ to be an isomorphism). In
    view of the explicit description of pure Gieseker bundles provided in
    \Cref{lemma: control of bubblings}, if $\widetilde{C}$ is not smooth, this
    forces $\widetilde{C}$ to be the nodal curve with a single node obtained by
    gluing to $(\mathbb{P}^1_{\widetilde{k}}, x^-, x^+)$ another copy of
    $\mathbb{P}^1_{\widetilde{k}}$ at $x^+$. Furthermore, in this case the
    restriction of $E$ to the contracted copy of $\mathbb{P}^1_{\widetilde{k}}$
    is isomorphic to $\cO_{\mathbb{P}^1_{\widetilde{k}}}(1)$, and the image of
    this point under $\Phi$ is the marked sheaf $(\mathbb{P}^1_{\widetilde{k}},
    x^-, x^+), \cO_{\mathbb{P}^1_{\widetilde{k}}}(d) \oplus
    (x^+)_*(\widetilde{k}), \phi: \cO_{\mathbb{P}^1_{\widetilde{k}}}(d) \oplus
    (x^+)_*(\widetilde{k}) \to (x^+)_*(\widetilde{k}))$ with $\phi$ given by the
    projection onto the second component. This is the marked sheaf corresponding
    to the split extension $0: \Spec(\widetilde{k}) \to
    \mathbb{A}^1_k/\mathbb{G}_m^3 \cong \Pur^{\circ}_{0,1,1}(1,d)$ under the
    identification in \Cref{lemma: identification two points pure with A1 mod
    G_m cubed}. This concludes the proof of the claim that $\Phi$ is an
    isomorphism away from the origin $0$.

    By \Cref{L:cBun algebraic} and \Cref{L: connected components of Gbun} the
    source $\StMap_{0,1,1}(B\mathbb{G}_m,d)$ is a smooth irreducible stack, and
    the same holds for the target $\Pur^{\circ}_{0,1,1}(1,d)$ by \Cref{lemma:
    identification two points pure with A1 mod G_m cubed}. To conclude, it
    suffices to show that the birational morphism $\Phi$ is finite. Since $\Phi$
    is proper, this is equivalent to showing that $\Phi$ is representable and
    quasifinite. Furthermore, since $\Phi$ is an isomorphism away from the
    origin, we just need to check that for any algebraically closed field
    extension $\widetilde{k} \supset k$ the fiber of $0: \Spec(\widetilde{k})
    \to \mathbb{A}^1_k/\mathbb{G}_m^3$ consists of a single point with trivial
    stabilizer group scheme. Recall that $0$ corresponds to the marked sheaf
    $(\mathbb{P}^1_{\widetilde{k}}, x^-, x^+), F, \phi: F \to
    (x^+)_*(\widetilde{k}))$ as above, where $F =
    \cO_{\mathbb{P}^1_{\widetilde{k}}}(d) \oplus (x^+)_*(\widetilde{k})$. In
    this case the Quot scheme $Q_F^1 \to \mathbb{P}^1_{\widetilde{k}}$ as in
    \Cref{defn: quot scheme} is the nodal curve with a single node obtained by
    gluing to $(\mathbb{P}^1_{\widetilde{k}}, x^-, x^+)$ another copy of
    $\mathbb{P}^1_{\widetilde{k}}$ at $x^+$. By \Cref{lemma: open immersion of
    GBun into kontsevich stable maps}, there is an injection from the set of
    $\widetilde{k}$-points of the $0$-fiber of $\StMap_{0,1,1}(B\mathbb{G}_m,d)$
    to the set of $\widetilde{k}$-points of the moduli $\cK_{0,1,1}(Q_F^1)$ of
    Kontsevich stable maps as in \Cref{defn: kontsevich stable maps}. In view of
    our discussion in the paragraph above, given a point $(\varphi:
    (\widetilde{C}, \sigma_{-}, \sigma_+) \to (\mathbb{P}^1_{\widetilde{k}},
    x^-, x^+), E)$ of the $0$-fiber, the corresponding Kontsevich stable map
    $\widetilde{C} \to Q_F^1$ is forced to be an isomorphism. This shows that
    there is a single $\widetilde{k}$-point in the $0$-fiber of $\Phi$, which
    moreover has trivial stabilizer group scheme.
\end{proof}

\begin{lemma} \label{lemma: weight of L_lev N=1} Under the identification
    $\StMap_{0,1,1}(B\mathbb{G}_m,d) \cong \bA^1_k / \bG_m^3$ constructed in
    \Cref{lemma: identification two points pure with A1 mod G_m cubed} and
    \Cref{lemma: two points Phi is an iso}, the line bundle $\cL_{\lev}$ from
    \Cref{defn: level line bundle} corresponds to the trivial line bundle
    tensored with the character of weight $((d-1)^2,-2d-1,-1)$.
\end{lemma}
\begin{proof}
    By \Cref{lemma: formula for level line bundle}, we have $\cL_{\lev} =
    \Phi^*(\cL_{\deg} \otimes \cL_{e}^{\vee} \otimes \cL_{k}^{\vee})$.The line
    bundle $\Phi^*(\cL_{e})$ corresponds to the morphism
    $\ev_+:\StMap_{0,1,1}(\B\mathbb{G}_m,d) \to B\mathbb{G}_m$, which
    corresponds to the morphism $\mathbb{A}^1_k/\mathbb{G}_m^3 \to
    B\mathbb{G}_m^3 \to B\mathbb{G}_m$ induced by the projection $\mathbb{G}_m^3
    \to \mathbb{G}_m$ onto the last coordinate. Hence $\Phi^*(\cL_{e})$
    corresponds to the trivial line bundle tensored by the character with weight
    $(0,0,-1)$. On the other hand, the line bundle $\Phi^*(\cL_k)$ corresponds
    to the composition $\StMap_{0,1,1}(B\mathbb{G}_m,d) \xrightarrow{\Psi \circ
    \Phi} \Pur^{\circ}_{0,2,0}(1,d) \xrightarrow{\ev_+} B\mathbb{G}_m$. In view
    of \Cref{L:Pur_{0,2,0}} and the identifications above, we conclude that
    $\Phi^*(\cL_k)^{\vee}$ corresponds to the trivial line bundle tensored by
    the character $(-d,-1,0)$.

    To conclude the proof, we need to show that
    the stabilizer $\mathbb{G}_m^3$ of the origin $0$ acts on the fiber
    $\Phi^*(\cL_{\deg})_0$ with character $(d^2-d+1, -2d,0)$. Under the
    identifications

    constructed above, the morphism $\Psi \circ \Phi$ is identified with the
    morphism $\mathbb{A}^1_k/\mathbb{G}_m^3 \to B\mathbb{G}_m^3 \to
    B\mathbb{G}_m^2$ induced by the homomorphism $\mathbb{G}_m^3 \to
    \mathbb{G}_m^2$ given by the projection to the first two factors. Hence, we
    are reduced to showing that the automorphism group $\mathbb{G}_m^2$ of the
    unique point $f: \Spec(k) \to \Pur^{\circ}_{0,2,0}(1,d)$ given by the triple
    $(\mathbb{P}^1_k, x^-, x^+, \cO_{\mathbb{P}^1_k}(d))$ acts on the pullback
    \[ f^*(\cL_{\deg}) = \det\left(R\Gamma(\mathbb{P}^1_k,
    \cO_{\mathbb{P}^1_k}(d))\right)^{\vee} \otimes
    \det\left(R\Gamma(\mathbb{P}^1_k, \cO_{\mathbb{P}^1_k}(d) \otimes
    \omega_{\mathbb{P}^1_k})\right)^{\vee}\] with weight $(d^2-d+1,-2d,0)$. This
    follows from a direct computation.
\end{proof}

For the following proposition, we use the notation $k\langle n \rangle$ to
denote the invertible sheaf on $B\mathbb{G}_m$ where $\mathbb{G}_m$ acts with
weight $-n$. We set $\transfer^a_{0,1,1}(1,d) := R(c\times \ev_+)_*\left(
\cL^a_{\lev} \otimes (\ev^-)^\ast(-) \right)$ for the stack
$\StMap_{0,1,1}(B\mathbb{G}_m,d)$.
\begin{prop}\label{P: N=1 transfer functor} For any $a>0$, the functor
    $\transfer^a_{0,1,1}(1,d) : \Dqc(B\bG_m) \to \Dqc(B\bG_m \times B\bG_m)$
    maps
    \[
        k\langle n \rangle \mapsto \left\{ \begin{array}{ll} k\langle -a(d-1)^2+ (d+1)m, a+m\rangle, & \text{if } m \leq 0 \\ 0, & \text{otherwise} \end{array} \right.,
    \]
    where $m:= (2d+1)a +n$.
\end{prop}
\begin{proof}
    This is a now direct computation, using our previous results. \endnote{We
    use the identification $\StMap_{0,1,1}(B\mathbb{G}_m,d) \cong
    \mathbb{A}^1_k/\mathbb{G}_m^3$ from \Cref{lemma: identification two points
    pure with A1 mod G_m cubed} and \Cref{lemma: two points Phi is an iso}. We
    have that the functor $\transfer^a_{0,1,1}(d)$ is given by
    $R(\ev_+)_*\left(\cL_{\lev} \otimes \ev_-^*(-)\right)$, where the morphism
    \[c \times \ev_+: \mathbb{A}^1_k/\mathbb{G}_m^3 \cong
    \StMap_{0,1,1}(B\mathbb{G}_m,d) \to B\mathbb{G}_m \times B\mathbb{G}_m
    \hookrightarrow \prestable_{0,1,1} \times B\mathbb{G}_m\] is isomorphic to
    the composition $f: \mathbb{A}^1_k/\mathbb{G}_m^3 \to B\mathbb{G}_m^3 \to
    B\mathbb{G}_m^2$ induced by the character $\left( (1,0,0),
    (0,0,1)\right):\mathbb{G}_m^3 \to \mathbb{G}_m^2$. Since $\ev_-:
    \mathbb{A}^1_k/\mathbb{G}_m^3 \to B\mathbb{G}_m^3 \to \mathbb{G}_m$ is
    induced by the character $(0,1,0): \mathbb{G}_m^3 \to \mathbb{G}_m$, we may
    use \Cref{lemma: weight of L_lev N=1} to write the image of $k\langle n
    \rangle$ under $\transfer^a_{0,1,1}(d)$ as
    \[R(c \times \ev_+)_*\left(\cL_{\deg}^{\otimes a}\langle 0, n,
    0\rangle\right) = Rf_*\left(\cO_{\mathbb{A}^1_k/\mathbb{G}_m^3}\langle
    -a(d-1)^2, n + (2d+1)a, a\rangle\right).\] This may be identified with the
    bi-graded vector space of invariants
    \[ \left(k[t]\langle -a(d-1)^2, n+ (2d+1)a, a\rangle\right)^{1 \times
    \mathbb{G}_m \times 1}\] under the subgroup $1 \times \mathbb{G}_m \times 1
    \subset \mathbb{G}_m^3$, where $\mathbb{G}_m^3$ acts on $t$ with weight
    $(-d-1,-1,1)$. We may rewrite this as
    \[ \bigoplus_{i \geq 0} \left(k\langle -a(d-1)^2 + (d+1)i, n + (2d+1)a +i,
    a-i\rangle\right)^{1 \times \mathbb{G}_m \times 1}.\] If $m:= (2d+1)a+n >0$,
    then the subspace of invariants above is $0$. Otherwise, the only summand
    that has nonzero invariants corresponds to $i= -m = -(2d+1)a -n $, and the
    bigraded vector space of invariants is $k\langle - a(d-1)^2 +(d+1)m,
    a+m\rangle$.}
\end{proof}

Now let $u$ be the generator of $\Rep(B\bG_m)$ for the ``target" $B\bG_m$ in our
theory corresponding to a character of weight $-1$.

\begin{lemma} \label{lemma: two point generating function for N=1} Consider the
    generating function defined by $\transfer^a_{0,1,1}(d)$ at the level of
    $K$-theory
    \[I^a_{0,1,1} := \sum_{d \in \bZ} z^d \cdot I^a_{0,1,1}(d) : \bZ[u^{\pm 1}]
    \to \bZ[q^{\pm 1}, u^{\pm 1}](\!(z^{-1})\!).\] If one forgets the action of
    the first factor of $B\bG_m$ in \Cref{P: N=1 transfer functor}
    (corresponding to setting $q=1$), then $I^a_{0,1,1}$ is given by
    \[
        u^n \mapsto \frac{z^{-\lceil \frac{n+a}{2a}\rceil} u^{n+a -2a \lceil \frac{n+a}{2a}\rceil}}{1-z^{-1}u^{-2a}}.
    \]
\end{lemma}
\begin{proof}
    This is a direct consequence of \Cref{P: N=1 transfer functor}.
\end{proof}

We note that the homomorphism $I^{a}_{0,1,1} : \bZ[u^{\pm 1}] \to \bZ[u^{\pm
1}](\!(z^{-1})\!)$ is uniquely determined by the difference equation $I(u^{2a}
f(u)) = z^{-1} I(f(u))$ and the identity
\[I(u^n) = \frac{u^n}{u^a - z^{-1}u^{-a}} \quad \text{for} \quad -a<n\leq a.\]
If we adjoin a root $\zeta = z^{1/2a}$ and let $\tilde{I}(u^n):=\zeta^{n-a}
I(u^n)$
, then
$\tilde{I}(\sum_i z^i f(u)) := \sum z^{i} \tilde{I}(f(u))$ is a
$\bZ(\!(\zeta^{-1})\!)$-linear extension $\tilde{I} : \bZ[u^{\pm
1}](\!(\zeta^{-1})\!) \to \bZ[u^{\pm 1}](\!(\zeta^{-1})\!)$ that satisfies
\[\tilde{I}(u^{2a} f(\zeta,u)) = \tilde{I}(f(\zeta,u)) \quad \text{and} \quad
\tilde{I}(u^n) = \frac{(\zeta u)^n}{(\zeta u)^a - (\zeta u)^{-a}}, \text{ for
}-a<n\leq a.\]

\subsubsection*{Three-point operations}

Just as in the case of two-point operations, we denote by
$\StMap_{0,2,1}(B\mathbb{G}_m,d) \subset \GBun_{0,2,1}(1,d)$ the preimage of the
open substack $\Spec(k) \hookrightarrow \prestable_{0,2,1}$ given by the smooth
$(2,1)$-marked curve $(\mathbb{P}^1_k, x_1^-, x_2^-, x^+)$, and similarly for
$\Pur^{\circ}_{0,2,1}(1,d)$.

\begin{lemma} \label{lemma: N=1 two vs three points} There is a well-defined
    forgetful morphism of stacks $\theta: \StMap_{0,2,1}(B\mathbb{G}_m,d) \to
    \StMap_{0,1,1}(B\mathbb{G}_m,d)$ given by
    \[(\varphi:\widetilde{C} \to C, \sigma_1^-, \sigma_2^-, \sigma^+, E) \mapsto
    (\varphi: \widetilde{C} \to C, \sigma_1^-, \sigma^+, E)\] fitting into a
    Cartesian diagram
    \[
        \begin{tikzcd}
            \StMap_{0,2,1}(B\mathbb{G}_m,d) \ar[d] \ar[r, "\theta"] & \StMap_{0,1,1}(B\mathbb{G}_m,d) \ar[d, "c"] \\
            \Spec(k) \ar[r] & B\mathbb{G}_m,
        \end{tikzcd}
    \]
    where the morphism $c:\StMap_{0,1,1}(B\mathbb{G}_m,d) \to B\mathbb{G}_m
    \hookrightarrow \prestable_{0,1,1}$ is given by $(\varphi: \widetilde{C} \to
    C, \sigma^-, \sigma^+, E)\mapsto (C, \varphi \circ \sigma^-, \varphi \circ
    \sigma^+)$.
\end{lemma}
\begin{proof}
    In order to see that $\theta$ is well-defined, we just need to check that
    for any given $(\varphi:\widetilde{C} \to C, \sigma_1^-, \sigma_2^-,
    \sigma^+, E)$ the associated morphism after forgetting one of the markings
    $\varphi: ( \widetilde{C}, \sigma_1^-, \sigma^+) \to (C, \varphi \circ
    \sigma_1^-, \varphi \circ \sigma^+)$ is a log-crepant contraction of
    $(1,1)$-marked curves. This is automatic since $\varphi$ is an isomorphism
    at $\sigma_2^+$ by \Cref{lemma: control of bubblings}. This fact about
    $\varphi$ can also be used to see that the diagram is Cartesian. \endnote{To
    see that the diagram is Cartesian, consider the fiber product $\mathcal{X}
    := \StMap_{0,1,1}(B\mathbb{G}_m,d) \times_{B\mathbb{G}_m} \Spec(k)$. For any
    $k$-scheme $T$, the $T$-points of $\mathcal{X}$ consists of tuples
    $(\varphi:\widetilde{C} \to C, \sigma^-, \sigma^+, E, \psi)$, where
    $(\varphi:\widetilde{C} \to C, \sigma^-, \sigma^+, E)$ is a $T$-point of
    $\GBun_{0,1,1}(1,d)$ and $\psi$ is an isomorphism $\psi: (C, \varphi \circ
    \sigma^-, \varphi \circ \sigma^+) \xrightarrow{\sim} (\mathbb{P}^1_T, x_1^1
    \times T, x^+ \times T)$ of $(1,1)$-marked curves. Since $\psi \circ
    \varphi: \widetilde{C} \to \mathbb{P}^1_T$ is an isomorphism over
    $\mathbb{P}^1_T \setminus x^+\times T$ (\Cref{lemma: control of bubblings}),
    the inverse image $(\psi \circ \varphi)^{-1}(x_2^-\times T)$ is well-defined
    and lies in the smooth locus of $\widetilde{C}$. We may then define a
    morphism of stacks $\mathcal{X} \to \StMap_{0,1,1}(B\mathbb{G}_m,d)$ given
    by
        \[ (\varphi:\widetilde{C} \to C, \sigma^-, \sigma^+, E, \psi) \mapsto
        (\varphi:\widetilde{C} \to C, \sigma^-, (\psi \circ
        \varphi)^{-1}(x_2^-\times T),  \sigma^+, E)\] which is an inverse of
        $\theta$.}
\end{proof}

Combining \Cref{lemma: N=1 two vs three points} with \Cref{lemma: identification
two points pure with A1 mod G_m cubed}, we get a commutative diagram
\[
    \begin{tikzcd}
        \StMap_{0,2,1}(B\mathbb{G}_m,d) \ar[d, "\sim"] \ar[r, "\theta"] & \StMap_{0,1,1}(B\mathbb{G}_m,d) \ar[d, "\sim"] \\
        \mathbb{A}^1_k/\mathbb{G}_m^2 \ar[r] & \mathbb{A}^1_k/\mathbb{G}_m^3,
    \end{tikzcd}
\]
where the morphism $\mathbb{A}^1_k/\mathbb{G}_m^2 \to
\mathbb{A}^1_k/\mathbb{G}_m^3$ is induced by the homomorphism $\mathbb{G}_m^2
\hookrightarrow \mathbb{G}_m^3$ given by the inclusion into the last two
coordinates of $\mathbb{G}_m^3$. We note that the evaluation morphism for the
negative markings $\ev_-: \StMap_{0,2,1}(B\mathbb{G}_m,d) \to B\mathbb{G}_m^2$
corresponds under the identification above to the morphism
$\mathbb{A}^1_k/\mathbb{G}_m^2 \to B\mathbb{G}_m^2 \to B\mathbb{G}_m^2$ induced
by the homomorphism $\mathbb{G}_m^2 \to \mathbb{G}_m^2$ given by $(t,s) \mapsto
(t,t)$.

\begin{prop}
    For any $a>0$, the generating function defined by $\transfer^a_{0,2,1}(1,d)$
    at the level of $K$-theory
    \[I^a_{0,2,1} := \sum_{d \in \bZ} z^d \cdot I^a_{0,2,1}(d) : \bZ[u^{\pm 1}]
    \otimes \bZ[u^{\pm 1}] \to \bZ[u^{\pm 1}](\!(z^{-1})\!)\] is given by
    \[
        u^m \otimes u^n \mapsto \frac{z^{-\lceil \frac{m+n+a}{2a}\rceil} u^{m+n+a-2a \lceil \frac{m+n+a}{2a}\rceil}}{1-z^{-1}u^{-2a}}.
    \]
\end{prop}
\begin{proof}
    In view of the discussion above, the image of $u^m \otimes u^n$ under
    $I_{0,2,1}^a$ agrees with the image of $u^{m+n}$ under $I_{0,1,1}^a$ after
    setting $q=1$ as in \Cref{lemma: two point generating function for N=1}.
    Hence, this proposition is an immediate consequence of the formula in
    \Cref{lemma: two point generating function for N=1}.
\end{proof}

\subsection{Gluing and orientation changing maps} \label{subsection: gluing and
forgetful morphisms}

In order to have a cleaner combinatorial description of the natural operations
for the moduli problem $\StMap_{g,n,p}(\GL_N, d)$, we use an oriented version of
modular graphs, inspired by the notion of modular graphs from
\cite{behrend-manin} and \cite{lee_quantum_k_theory_i}.

\begin{defn}[Directed modular graphs]
    A directed modular graph $\tau$ consists of the following:
    \begin{enumerate}
        \item Finite sets of vertices $V_{\tau}$ and edges $E_{\tau}$, along
        with two maps $s,t : E_{\tau} \to V_{\tau}$ called the source and target
        morphisms;
        \item Finite sets $\mathbf{P}_{\tau}$ and $\mathbf{N}_{\tau}$ along with
        ``marking functions" $m_{\tau, \rm out} : \mathbf{P}_{\tau} \to
        V_{\tau}$ and $m_{\tau, \rm in} : \mathbf{N}_{\tau} \to V_{\tau}$. We
        will typically identify $\mathbf{P}_{\tau} = \{1,\ldots,p\}$ and
        $\mathbf{N}_{\tau} = \{1,\ldots,n\}$ for some $p,n \geq 0$;
        \item A genus function $g_{\tau} : V_{\tau} \to \bZ_{\geq 0}$ and a
        degree function $d_{\tau}: V_{\tau} \to \bZ$.
    \end{enumerate}
\end{defn}

Given a vertex $v \in V$, we think of the preimage $m_{\tau, \rm in}^{-1}(v)$ as
a set of ``incoming'' tails, and we think of $m_{\tau, \rm out}^{-1}(v)$ as a
set of ``outgoing'' tails.

\begin{defn}[Incoming and outgoing valence]
    Let $\tau$ be a modular graph. For each vertex $v \in V$, let $\#_{\rm
    in}(v) = |t^{-1}(v)| + |m_{\tau, \rm in}^{-1}(v)|$ and $\#_{\rm out}(v) =
    |s^{-1}(v)| + |m_{\tau, \rm out}^{-1}(v)|$.
\end{defn}

\begin{defn}[Cutting an edge] \label{defn: cutting an edge} Let $\tau$ be a
    directed modular graph, and let $e \in E_{\tau}$ be an edge. The directed
    modular graph obtained by cutting $e$ is the unique directed modular graph
    $\sigma$ given as follows:
    \begin{itemize}
        \item $V_{\sigma} = V_{\tau}$, and the genus and degree functions of
        $\sigma$ and $\tau$ agree.
        \item $E_{\sigma} = E\setminus \{e\}$, and the source and target
        morphisms are obtained by restricting the corresponding morphisms for
        $\tau$.
        \item $\mathbf{P}_{\sigma} = \mathbf{P}_{\tau} \sqcup \{p+1\}$ and
        $\mathbf{N}_{\sigma} = \mathbf{N}_{\tau} \sqcup \{n+1\}$. The morphisms
        $m_{\sigma, \mathrm{in}}$ and $m_{\sigma, \mathrm{out}}$ agree with
        $m_{\tau, \mathrm{in}}$ and $m_{\tau, \mathrm{out}}$ on
        $\mathbf{N}_{\tau} \subset \mathbf{N}_{\sigma}$ and $\mathbf{P}_{\tau}
        \subset \mathbf{P}_{\sigma}$, and we set $m_{\sigma, \mathrm{in}}(n+1) =
        s(e)$ and $m_{\sigma, \mathrm{in}}(p+1) = t(e)$
    \end{itemize}
\end{defn}

\begin{defn}[Stack of stable maps with fixed directed modular graph]
\label{defn: stacks of stable maps with fixed modular graph}. Let $\tau \mapsto
\StMap_{\tau}(\GL_N)$ be the unique assignment from directed modular graphs to
smooth algebraic stacks over $S$ satisfying the following properties:
    \begin{itemize}
        \item If $\tau$ satisfies $E_{\tau} = \emptyset$, then we set
        $\StMap_{\tau}(\GL_N) = \prod_{v \in V_{\tau}} \StMap_{g(v),\#_{\rm
        in}(v),\#_{\rm out}(v)}(\GL_N, d_{\tau}(v))$.
        \item If $\sigma$ is obtained by cutting an edge of $\tau$ as in
        \Cref{defn: cutting an edge}, then we have a fiber product diagram
              \[
                  \begin{tikzcd}
                      \StMap_{\tau}(\GL_N) \ar[r] \ar[d] &  \StMap_{\sigma}(\GL_N)\ar[d, "\ev_{n+1}^- \times \ev_{p+1}^+"]  \\
                      \BGL_N \ar[r, "\Delta"] & \BGL_N \times \BGL_N,
                  \end{tikzcd}
              \]
              where $\Delta$ is the diagonal morphims, and the evaluation
              morphisms $\ev_{n+1}^-$ and $\ev_{p+1}^+$ are obtained by
              restricting the relevant universal bundles at the markings
              corresponding to $\{p+1\} \in \mathbf{P}_{\sigma}$ and $\{n+1\}
              \in \mathbf{N}_{\sigma}$.
    \end{itemize}
\end{defn}

\begin{defn}[Contracting an edge] \label{defn: contracting an edge} Let $\tau$
    be a directed modular graph, and let $e \in E_{\tau}$ be an edge. The
    directed modular graph obtained by contracting $e$ is the unique directed
    modular graph $\sigma$ given as follows:
    \begin{itemize}
        \item $V_{\sigma} = (V_{\tau} \setminus \{s(e), t(e)\}) \sqcup \{*\}$
        and $E_{\sigma} = E\setminus \{e\}$. The source and target morphisms are
        obtained by restricting the corresponding morphisms for $\tau$.
        \item $\mathbf{P}_{\sigma} = \mathbf{P}_{\tau}$ and $\mathbf{N}_{\sigma}
        = \mathbf{N}_{\tau}$. The morphisms $m_{\sigma, \mathrm{in}}$ and
        $m_{\sigma, \mathrm{out}}$ are given by compositions $\mathbf{N}_{\tau}
        \xrightarrow{m_{\tau, \mathrm{in}}} V_{\tau} \xrightarrow{\pi}
        V_{\sigma}$ and $\mathbf{P}_{\tau} \xrightarrow{m_{\tau, \mathrm{out}}}
        V_{\tau} \xrightarrow{\pi} V_{\sigma}$, where $\pi: V_{\tau}
        \xrightarrow{\pi} V_{\sigma}$ is the morphism that restricts to the
        identity on $V_{\tau} \setminus \{s(e), t(e)\}$ and satisfies $\pi(s(e))
        = \pi(t(e)) = *$.
        \item The genus and degree functions $g_{\sigma}, d_{\sigma}$ restrict
        to the genus and degree functions for $\tau$ on the subset $V_{\tau}
        \setminus \{s(e), t(e)\}$.
        \item We set $g_{\sigma}(*) = g_{\tau}(s(e)) + g_{\tau}(t(e))$ if $s(e)
        \neq t(e)$, and otherwise we set $g_{\sigma}(*) = g_{\tau}(t(e)) +1$.
        \item We set $d_{\sigma}(*) = d_{\tau}(s(e)) + d_{\tau}(t(e))$ if $s(e)
        \neq t(e)$, and otherwise we set $d_{\sigma}(*) = d_{\tau}(t(e))$.
    \end{itemize}
\end{defn}

\begin{prop}[Gluing morphism] \label{prop: gluing morphisms} Let $\tau$ be a
    directed modular graph, and choose an edge $e \in E_{\tau}$. Let $\sigma$ be
    the directed modular graph obtained by contracting $e$, as in \Cref{defn:
    contracting an edge}. Then, there exists a morphism $\Gamma_{\tau \to
    \sigma}: \StMap_{\tau}(\BGL_N) \to \StMap_{\sigma}(\BGL_N)$ compatible with
    the standard gluing morphism of underlying stable target curves.
\end{prop}
\begin{proof}
    Similarly as in the case of stable maps to a projective target, our
    construction of the gluing morphism is local around the markings
    corresponding to the source $s(e)$ and target $t(e)$ of the contracted edge
    $e$. It is then sufficient to give a description of the gluing morphism in
    the case where $\tau$ has a single edge $E_{\tau}=\{e\}$ and
    $\mathbf{P}_{\tau} = \mathbf{N}_{\tau}=\emptyset$, which generalizes in the
    evident way to more general directed modular graphs. We distinguish two
    possible cases.

    \noindent \textbf{Case 1: $\tau$ has two vertices.} Then
    $\StMap_{\tau}(\GL_N)$ fits into a Cartesian diagram
    \[
        \begin{tikzcd}
            \StMap_{\tau}(\GL_N) \ar[r] \ar[d] &  \StMap_{g_1,0,1}(\GL_N, d_1) \times \StMap_{g_2, 1, 0}(\BGL_N, d_2) \ar[d, "\ev_+ \times \ev_-"]  \\
            \BGL_N \ar[r, "\Delta"] & \BGL_N \times \BGL_N.
        \end{tikzcd}
    \]
    For any scheme $T$, a $T$-point of $\StMap_{\tau}(\GL_N)$ corresponds to
    pairs of points $(\varphi_1: \widetilde{C}_1 \to C_1, \sigma^+, E_1) \in
    \StMap_{g_1,0,1}(\GL_N, d_1)(T)$ and $(\varphi_2: \widetilde{C}_2 \to C_2,
    \sigma^-, E_2) \in \StMap_{g_2, 1, 0}(\BGL_N, d_2)(T)$ along with an
    isomorphism $\psi: (\sigma^+)^*(E_1) \xrightarrow{\sim} (\sigma^-)^*(E_2)$.
    We define $\Gamma_{\tau \to \sigma}$ to be the assignment that sends such
    data to:
    \begin{itemize}
        \item The log-crepant contraction of $(0,0)$-marked curves $\varphi:
        \widetilde{C} \to C$, where $\widetilde{C}$ is obtained by gluing
        $\widetilde{C}_1$ and $\widetilde{C}_2$ along the markings $\sigma^+$
        and $\sigma^-$, and $C$ is obtained by gluing $C_1$ and $C_2$ along the
        markings $\varphi_1 \circ \sigma^+$ and $\varphi_2 \circ \sigma^-$.
        \item The vector bundle $E$ on $\widetilde{C}$ obtained by gluing $E_1$
        to $E_2$ via the identification $\psi$.
    \end{itemize}
    It can be checked that $(\varphi: \widetilde{C} \to C, E)$ is $T$-point of
    $\StMap_{g_1+g_2, 0,0}(\BGL_N, d_1+d_2) \cong \StMap_{\sigma}(\BGL_N)$ (one
    way to see this is using the concrete description of pure Gieseker bundles
    in the introduction).

    \noindent \textbf{Case 2: $\tau$ has one vertex.} We have that
    $\StMap_{\tau}(\GL_N)$ fits into a Cartesian diagram
    \[
        \begin{tikzcd}
            \StMap_{\tau}(\GL_N) \ar[r] \ar[d] &  \StMap_{g,1,1}(\GL_N, d) \ar[d, "\ev_+ \times \ev_-"]  \\
            \BGL_N \ar[r, "\Delta"] & \BGL_N \times \BGL_N.
        \end{tikzcd}
    \]
    A $T$-point of $\StMap_{\tau}(\GL_N)$ corresponds a point $(\varphi:
    \widetilde{C} \to C, \sigma^+, E) \in \StMap_{g,1,1}(\GL_N, d)(T)$ along
    with an isomorphism $\psi: (\sigma^+)^*(E) \xrightarrow{\sim}
    (\sigma^-)^*(E)$. We define $\Gamma_{\tau \to \sigma}$ to be the assignment
    that sends any such pair to:
    \begin{itemize}
        \item The log-crepant contraction of $(0,0)$-marked curves $\varphi':
        \widetilde{C}' \to C'$, where $\widetilde{C}'$ is obtained by gluing
        $\widetilde{C}$ along the markings $\sigma^+$ and $\sigma^-$, and $C'$
        is obtained by gluing $C$ along the markings $\varphi \circ \sigma^+$
        and $\varphi \circ \sigma^-$.
        \item The vector bundle $E'$ on $\widetilde{C}'$ obtained by gluing $E$
        via the identification $\psi$.
    \end{itemize}
    It follows as in the previous case that this induces a well-defined morphism
    $\Gamma_{\tau \to \sigma}:  \StMap_{\tau}(\GL_N) \to
    \StMap_{g+1,0,0}(\BGL_N, d) \cong \StMap_{\sigma}(\BGL_N)$.
\end{proof}

\begin{notn}[Level line bundle] \label{notn: level line bundle} By \Cref{defn:
    stacks of stable maps with fixed modular graph}, for every directed graph
    $\tau$ we have a morphism
    \[ \StMap_{\tau}(\BGL_N) \to \StMap_{\overline{\tau}}(\BGL_N) \cong \prod_{v
    \in V_{\tau}} \StMap_{g(v),\#_{\rm in}(v),\#_{\rm out}(v)}(\GL_N, d(v)), \]
    where $\overline{\tau}$ is the modular graph obtained by cutting all the
    edges of $\tau$. We denote by $\cL_{\lev}$ the line bundle on
    $\StMap_{\tau}(\BGL_N)$ obtained by pulling back $\boxtimes_{v \in V_{\tau}}
    \cL_{\lev}$ under this morphism.
\end{notn}

\begin{prop}[Compatibility of level line bundle with gluing] \label{prop:
    compatibility of level line bundle with gluing} Let $\tau$ be a directed
    modular graph, and choose an edge $e \in E_{\tau}$. Let $\sigma$ be the
    directed modular graph obtained by contracting $e$. Then, we have
    $\cL_{\lev} = \Gamma_{\tau \to \sigma}^*(\cL_{\lev})$.
\end{prop}
\begin{proof}
    As in the proof of \Cref{prop: gluing morphisms}, we may reduce to the case
    when $E_{\tau} = \{e\}$ and $\mathbf{N}_{\tau} = \mathbf{P}_{\tau} =
    \emptyset$. We spell out the proof in the case when $\tau$ has two vertices;
    the other case follows from a very similar computation. In this case the
    gluing morphism is $\Gamma_{\tau \to \sigma}$ may be written as
    \[ \Gamma_{\tau \to \sigma}: \StMap_{g_1,0,1}(\GL_N,
    d_1)\times_{\BGL_N}\StMap_{g_2,1,0}(\GL_N, d_2) \to \StMap_{g_1+g_2,
    0,0}(\BGL_N, d_1+d_2)\] as described in the proof of \Cref{prop: gluing
    morphisms}. Our goal is to show that there is an isomorphism $\Gamma_{\tau
    \to \sigma}^*(\cL_{\lev}) \cong (\mathrm{pr}_1)^*(\cL) \otimes
    (\mathrm{pr}_2)^*(\cL)$, where $\mathrm{pr}_1$ and $\mathrm{pr}_2$ denote
    the first and second projections of the fiber product
    $\StMap_{g_1,0,1}(\GL_N, d_1)\times_{\BGL_N}\StMap_{g_2,1,0}(\GL_N, d_2)$.

    Let $(\varphi: \widetilde{C}_1 \to C_1, \sigma^+, E_1)$ and $(\varphi:
    \widetilde{C}_2 \to C_2, \sigma^-, E_2)$ denote the pullbacks of the
    universal pure Gieseker bundles via $\mathrm{pr}_1$ and $\mathrm{pr}_2$, and
    let $(\varphi: \widetilde{C} \to C, E)$ denote the pullback of the universal
    pure Gieseker bundle via $\Gamma_{\tau \to \sigma}$. Set
    $E_{\ev}=(\sigma^+)^*(E_1) \cong (\sigma^-)^*(E_2)$. The relative curve
    $\pi: \widetilde{C} \to \StMap_{\tau}(\BGL_N)$ is obtained by gluing $\pi_1:
    \widetilde{C}_1 \to \StMap_{\tau}(\BGL_N)$ and $\pi_2: \widetilde{C}_2 \to
    \StMap_{\tau}(\BGL_N)$ along $\sigma^+$ and $\sigma^-$. Let $j_1:
    \widetilde{C}_1 \hookrightarrow \widetilde{C}$ and $j_2: \widetilde{C}_2
    \hookrightarrow \widetilde{C}$ denote the corresponding closed immersions,
    and set $\sigma = j_1 \circ \sigma^+ = j_2 \circ \sigma^-$. By the
    definition of $\Gamma_{\tau \to \sigma}$, we have a short exact sequence
    \[ 0 \to E \to (j_1)_*(E_1) \oplus (j_2)_*(E_2) \to \sigma_*(E_{\ev}) \to
    0,\] which yields an equality $[E] = [(j_1)_*(E_1)] +  [(j_2)_*(E_2)] -
    [\sigma_*(E_{\ev})]$ in $K$-theory of $\widetilde{C}$. Similarly, we have a
    short exact sequence
    \[ 0 \to E \otimes \omega_{\pi} \to (j_1)_*(E_1 \otimes
    \omega_{\pi_1}(\sigma^+)) \oplus (j_2)_*(E_2 \otimes
    \omega_{\pi_2}(\sigma^-)) \to \sigma_*(E_{\ev} \otimes
    \sigma^*(\omega_{\pi})) \to 0.\] Using the residue isomorphism
    $\sigma^*(\omega_{\pi})) \cong (\sigma^+)^*(\omega_{\pi_1}(\sigma^+)) \cong
    \cO_{\StMap_{\tau}(\BGL_N)}$, this yields:
    \[[ E \otimes \omega_{\pi}] =  [ (j_1)_*(E_1 \otimes
    \omega_{\pi_1}(\sigma^+))] + [(j_2)_*(E_2 \otimes \omega_{\pi_2}(\sigma^-))]
    - [\sigma_*(E_{\ev})].\] Now using the short exact sequence
    \[ 0 \to E_1 \otimes \omega_{\pi_1} \to E_1 \otimes \omega_{\pi_1}(\sigma^+)
    \to (\sigma^+)_*\left( E_{\ev}
    \otimes(\sigma^+)^*(\omega_{\pi_1}(\sigma^+))\right) \to 0 \] and the
    residue isomorphism $(\sigma^+)^*(\omega_{\pi_1}(\sigma^+)) \cong
    \cO_{\StMap_{\tau}(\BGL_N)}$ we get $[(j_1)_*(E_1 \otimes
    \omega_{\pi_1}(\sigma^+))] = [(j_1)_*(E_1 \otimes \omega_{\pi_1})] +
    [\sigma_*(E_{\ev})]$. Similarly we have $[(j_2)_*(E_2 \otimes
    \omega_{\pi_2}(\sigma^-))] = [(j_2)_*(E_2 \otimes \omega_{\pi_2})] +
    [\sigma_*(E_{\ev})]$. Putting everything together, we conclude that
    \[ [E] + [E \otimes \omega_{\pi}] = [(j_1)_*(E_1)] + [(j_1)_*(E_1 \otimes
    \omega_{\pi_1})] + [(j_2)_*(E_2)]+ [(j_2)_*(E_2 \otimes \omega_{\pi_2})].\]
    Pushing forward to $\StMap_{g_1,0,1}(\GL_N,
    d_1)\times_{\BGL_N}\StMap_{g_2,1,0}(\BGL_N, d_2)$ and taking determinant we
    see that $(\mathrm{pr}_1)^*(\cL_{\lev}) \otimes
    (\mathrm{pr}_2)^*(\cL_{\lev}) = \Gamma_{\tau \to \sigma}^*(\cL_{\lev})$, as
    desired.
\end{proof}

\begin{defn}[Changing the orientation of a negative marking] \label{defn:
    changing the orientation of the last negative marking} Let $\tau$ be a
    directed modular graph such that $\mathbf{N}_{\tau}$ is nonempty. Let $v_n:=
    m_{\mathrm{in}}(n) \in V_{\tau}$. The directed modular graph obtained by
    changing the orientation of the last negative marking is the unique directed
    modular graph $\sigma$ given as follows:
    \begin{itemize}
        \item $V_{\tau}= V_{\sigma}$, $E_{\tau} = E_{\sigma}$ and the functions
        $s,t,g_{\sigma}$ agree with those for $\tau$.
        \item The degree function for $d_\sigma$ agrees with $d_\tau$ on
        $V_{\tau} \setminus \{v_n\}$ and sends the vertex $v_n$ to
        $d_{\tau}(v_n) -1$.
        \item $\mathbf{N}_{\sigma} = \mathbf{N}_{\tau} \setminus \{n\}$ and
        $m_{\sigma, \mathrm{in}}$ is obtained by restricting $m_{\tau,
        \mathrm{in}}$ to $\mathbf{N}_{\sigma} \subset \mathbf{N}_{\tau}$.
        \item $\mathbf{P}_{\sigma} = \mathbf{P}_{\tau} \sqcup \{p+1\}$, and the
        function $m_{\mathrm{out}}: \mathbf{P}_{\sigma} \to V_{\sigma}$ agrees
        with $m_{\tau, \mathrm{out}}$ on $\mathbf{P}_{\tau} \subset
        \mathbf{P}_{\sigma}$ and sends $p+1$ to $v_n$.
    \end{itemize}
\end{defn}

\begin{prop}[Negative orientation changing morphisms]
    Let $\tau$ be a modular graph with $\mathbf{N}_{\tau}$ nonempty. Let
    $\sigma$ be the directed modular graph obtained by changing the orientation
    of the last negative marking, as in \Cref{defn: changing the orientation of
    the last negative marking}. Then, there is a morphism $\Gamma_{- \to +}:
    \StMap_{\tau}(\BGL_N) \to \StMap_{\sigma}(\BGL_N)$ satisfying $\Gamma_{- \to
    +}^*(\cL_{\lev}) \cong \cL_{\lev}$.
\end{prop}
\begin{proof}
    The construction is local around the negative marking corresponding to
    $\{n\} \in \mathbf{N}_{\tau}$. We may reduce to the case when $\tau$ has a
    single vertex, and satisfies $\mathbf{P}_{\tau} = E_{\tau} = \emptyset$ and
    $\mathbf{N}_{\tau} = \{1\}$. In this case we have $\StMap_{\tau}(\BGL_N)
    \cong \StMap_{g, 1, 0}(\BGL_N, d)$ and $\StMap_{\sigma}(\BGL_N) \cong
    \StMap_{g, 0,1}(\BGL_N, d-1)$. We define the assignment $\Gamma_{- \to +}$
    at the level of $T$-points by sending a family $(\varphi: \widetilde{C} \to
    C, \sigma^-,E)$ in $\StMap_{g, 1, 0}(\BGL_N, d)$ to the same tuple
    $(\varphi: \widetilde{C} \to C, \sigma^-, E)$, which we view as a family in
    $\StMap_{g, 0,1}(\BGL_N, d-1)$ by regarding $\sigma^-$ as a positive
    marking. The compatibility of this morphism with the level line bundle
    $\cL_{\lev}$ is evident.
\end{proof}

\begin{defn}[Changing the orientation of a positive marking] \label{defn:
    changing the orientation of the last positive marking} Let $\tau$ be a
    directed modular graph such that $\mathbf{P}_{\tau}$ is nonempty. We say
    that $\sigma$ is the oriented modular graph obtained by changing the
    orientation of the positive marking if $\tau$ satisfies the properties of
    the oriented graph obtained by changing the last negative marking of
    $\sigma$, with the exception that we require $d_{\sigma} = d_{\tau}$.
\end{defn}

\begin{prop}[Positive orientation changing morphisms]
    Let $\tau$ be a modular graph with $\mathbf{P}_{\tau}$ nonempty. Let
    $\sigma$ be the directed modular graph obtained by changing the orientation
    of the last positive marking, as in \Cref{defn: changing the orientation of
    the last positive marking}. Then, we have a nontrivial morphism $\Gamma_{+
    \to -}: \StMap_{\tau}(\BGL_N) \to \StMap_{\sigma}(\BGL_N)$.
\end{prop}
\begin{proof}
    The construction is local around the positive marking corresponding to
    $\{p\} \in \mathbf{P}_{\tau}$. We may reduce to the case when $\tau$ has a
    single vertex, and satisfies $\mathbf{N}_{\tau} = E_{\tau} = \emptyset$ and
    $\mathbf{P}_{\tau} = \{1\}$. In this case we have $\StMap_{\tau}(\BGL_N)
    \cong \StMap_{g, 0, 1}(\BGL_N, d)$ and $\StMap_{\sigma}(\BGL_N) \cong
    \StMap_{g,1,0}(\BGL_N, d)$. Let $T$ be a scheme, and choose a $T$-point
    $(\varphi: \widetilde{C} \to C, \sigma^+,E)$ of $\StMap_{g, 0,1}(\BGL_N,
    d)$. Recall that by \Cref{lemma: control of bubblings} the morphism
    $\varphi$ is an isomorphism over the open complement $U \subset C^{sm}$ of
    $\varphi \circ \sigma^+$ inside the $T$-smooth locus $C^{sm} \subset C$. We
    denote by $\widetilde{C}'$ the curve obtained by gluing $C^{sm}$ and
    $\widetilde{C} \setminus \sigma^+$ along their common open subset $U$. Then
    we have a factorization $\varphi: \widetilde{C} \xrightarrow{f}
    \widetilde{C}' \xrightarrow{\varphi'} C$, and if we set $\sigma^- := f \circ
    \sigma^+$, then $(\varphi': \widetilde{C}' \to C, \sigma^-)$ is a
    log-crepant contraction of $(1,0)$-marked curves. It follows from the
    definition of pure Gieseker bundles that the pushforward $E' :=
    f_*(E(-\sigma^+))$ is a vector bundle on $\widetilde{C}'$, and that the
    tuple $(\varphi': \widetilde{C}' \to C, \sigma^-, E')$ is a $T$-family of
    pure Gieseker bundles in $\StMap_{g,1,0}(\BGL_N, d)$. We define $\Gamma_{+
    \to -}$ by the assignment $(\varphi: \widetilde{C} \to C, \sigma^+,E)
    \mapsto (\varphi': \widetilde{C}' \to C, \sigma^-, E')$.
\end{proof}

\begin{remark}
    The level line bundle $\cL_{\lev}$ does not seem to be preserved by pulling
    back under the positive orientation changing morphism $\Gamma_{+ \to -}$.
    This illustrates the lack of symmetry between the positive and negative
    markings in our moduli problem.
\end{remark}

\begin{remark}[Forgetful morphisms]
    Unlike in the case of the stack $\overline{\cM}_{g,n}(X)$ of stable maps to
    a projective scheme $X$, it is unclear how to define forgetful morphisms for
    $\StMap_{\tau}(\BGL_N)$ that remove a marking.
\end{remark}

\section{Projective targets} \label{section: projective targets} In this
section, we explain how to adapt the techniques developed in \Cref{section:
moduli of gieseker bundles} to study the moduli theory of gauged maps from
marked prestable curves into a projective $\GL_N$-scheme. From now on, we fix
the following.
\begin{context} \label{context: moduli section} Let $Z$ be a projective
    $S$-scheme equipped with an action of $\GL_N$. We fix a choice $\cO(1)$ of
    $\GL_N$-equivariant very ample line bundle on $Z \to S$.
\end{context}

\subsection{The moduli stack}

\begin{defn}[Stack of Gieseker maps] \label{defn: stack of gieseker stable maps}
    We define $\GMap_{g,n,p}(Z/\GL_N)$ to be the stack that sends an $S$-scheme
    $T$ into the groupoid of tuples $(\varphi: \widetilde{C} \to C,
    \sigma_{\bullet}^{\pm}, u)$, where $(\varphi: \widetilde{C} \to C,
    \sigma_{\bullet}^{\pm}, u)$ is a $T$-family of $(n,p)$-marked log-crepant
    contractions in $\precont_{g,n,p}$ and $u: \widetilde{C} \to Z/\GL_N$ is a
    morphism such that the following are satisfied:
    \begin{enumerate}[(1)]
        \item  Let $E$ be the rank $N$ vector bundle on $\widetilde{C}
        \xrightarrow{u} Z/\GL_N \to \BGL_N$. Then the counit $\varphi^*
        \varphi_*(E) \to E$ is surjective and the pushforward $\cF:=
        \varphi_*(E(\sum_{i=1}^p \sigma_i^+))$ is a $T$-family of pure sheaves
        of dimension $1$ on $C$. This induces a morphism $\widetilde{C} \to
        Q^N_{\cF}$ as in \Cref{defn: quot scheme}. We denote by $E_{univ}$ the
        universal rank $N$ vector bundle quotient on $ Q^N_{\cF}$, thought of as
        a $\GL_N$-bundle.
        \item Let $s$ be the section of the associated fiber bundle $E(Z) := E
        \times^{\GL_N} Z \to C$ corresponding to $u: \widetilde{C} \to Z/\GL_N$.
        Then the induced morphism $\widetilde{C} \to E_{univ} \times^{\GL_N} Z$
        is a $T$-family of $(n+p)$-pointed Kontsevich stable maps to the
        $T$-projective scheme $E_{univ} \times^{\GL_N} Z$.
    \end{enumerate}
\end{defn}

We may think of a point $(\varphi: \widetilde{C} \to C, \sigma_{\bullet}^{\pm},
u)$ alternatively as the data of $(\varphi: \widetilde{C} \to C,
\sigma_{\bullet}^{\pm}, E,s)$. We will freely go back and forth between these
two points of view.

\begin{remark}
    The condition in \Cref{defn: stack of gieseker stable maps}(1) concerning
    the purity of the $T$-fibers of $\varphi_*(E(\sum_{i=1}^p \sigma_i^+))$ is
    meaningful because its formation commutes with base-change on $T$
    (\Cref{lemma: pushforward lemma}).
\end{remark}

\begin{remark} \label{remark: alternative description kontsevich stability} With
    notation as in \Cref{defn: stack of gieseker stable maps}, there is a line
    bundle $\cO_{E(Z)}(1)$ on $E(Z)$ induced by the $\GL_N$-equivariant line
    bundle on $Z$. Kontsevich stability is equivalent to the condition that for
    all $t \in T$ the line bundle \[\omega_{\widetilde{C}_{t}/t}(\sum_{i=1}^p
    \sigma_{i,t}^+ + \sum_{i=1}^n \sigma_{i,t}^-) \otimes
    \text{det}(E|_{\widetilde{C}_t})^{\otimes m} \otimes
    s^*(\cO_{E(Z)}(1))|_{\widetilde{C}_t}\] is $\varphi_t$-ample for all $m\gg
    0$.
\end{remark}

We consider $\GMap_{g,n,p}(Z/\GL_N)$ as a stack over $\prestable_{g,n,p}$ via
the morphism $(\varphi: \widetilde{C} \to C, \sigma_{\bullet}^{\pm}, u) \to (C,
\varphi \circ \sigma_{\bullet}^{\pm})$.
\begin{prop}
    The stack $\GMap_{g,n,p}(Z/\GL_N)$ is an algebraic stack locally of finite
    type and with affine relative diagonal over $\prestable_{g,n,p}$.
\end{prop}
\begin{proof}
    Let $\cC \to \precont_{g,n,p}$ be the universal source curve.
    $\GMap_{g,n,p}(Z/\GL_N)$ is an open substack of the mapping stack
    $\Map_{\precont_{g,n,p}}(\cC, (Z/\GL_N)\times\precont_{g,n,p})$, and hence
    it is an algebraic stack locally of finite type over $\precont_{g,n,p}$ by
    \cite[Thm. 1.2]{hall-rydh-tannakahom}. In particular, it is also locally of
    finite type over $\prestable_{g,n,p}$. The proof that
    $\GMap_{g,n,p}(Z/\GL_N) \to \prestable_{g,n,p}$ has affine relative diagram
    is similar to the argument in \Cref{L:cBun algebraic}, replacing the usage
    of $\det(E)$ with the $\varphi$-ample line bundle
    $\omega_{\widetilde{C}}(\sum_i \sigma_{i,t}^{\pm}) \otimes
    \text{det}(E)^{\otimes m} \otimes s^*(\cO_{E(Z)}(1))$ with $m$ sufficiently
    large as in \Cref{remark: alternative description kontsevich stability}
\end{proof}

\begin{defn}[{\cite[Sect. 2.2.]{gauged_theta_stratifications}}] \label{defn:
    numerical neron-severi} Let $N^1_S(Z/\GL_N)$ denote the group of isomorphism
    classes of line bundles on $Z/\GL_N$ modulo the subgroup of those line
    bundles $\cL \in \Pic(Z/\GL_N)$ such that for all proper reduced curves $C$
    over a field $k$ and morphisms $C \to Z/\GL_N$ such that the composition $C
    \to S$ factors through the structure morphism $C \to \Spec(k)$, we have
    $\text{deg}(\cL|_C) =0$. By definition $N^1_S(Z/\GL_N)$ is a torsion-free
    abelian group, we denote by $H_2(Z/\GL_N):=
    \Hom_{\mathbb{Z}}(N^1_S(Z/\GL_N), \mathbb{Z})$ its linear dual.
\end{defn}

Let $T$ be a scheme, and choose $T \to \GMap_{g,n,p}(Z/\GL_N)$ corresponding to
a $T$-family $(\varphi: \widetilde{C} \to C, \sigma_{\bullet}^{\pm}, u:
\widetilde{C} \to Z/\GL_N)$. There is an induced locally constant function
$\delta: |T| \to H_2(Z/\GL_N)$, which sends a point $t \in T$ to the
homomorphism $\delta_t: N^1_S(Z/\GL_N) \to \mathbb{Z}$ given on equivalence
classes $[\cL]$ of $\cL \in \Pic(Z/\GL_N)$ by $ \delta_t([\cL]) :=
\text{deg}(u^*(\cL)|_{X_t})$.

\begin{defn}[Degree of Gieseker maps] \label{defn: degree of Gieseker maps} We
    say that the family $T \to \GMap_{g,n,p}(Z/\GL_N)$ has degree $d \in
    H_2(Z/\GL_N)$ if $\delta = d$ identically on $|T|$. We denote by
    $\GMap_{g,n,p}(Z/\GL_N,d) \subset \GMap_{g,n,p}(Z/\GL_N)$ the corresponding
    open and closed substack parameterizing points of degree $d$.
\end{defn}

\subsection{Properness over the stack of pure marked sheaves}

\begin{lemma} \label{lemma: PhiZ is well-defined} The assignment
    \[(\varphi : \tilde{C} \to C, \sigma^\pm_\bullet, E,s) \mapsto (C, \varphi
    \circ \sigma^\pm_\bullet, \varphi_*(E), (\varphi \circ
    \sigma_{\bullet}^+)^*(\varphi_*(E)) \to E|_{\sigma_{\bullet}^+})\] defines a
    morphism of algebraic stacks $\Phi_Z : \GMap_{g,n,p}(Z/\GL_N) \to
    \Pur_{g,n,p}(N)$.
\end{lemma}
\begin{proof}
    The same proof as in \Cref{lemma: Phis is well-defined} applies here.
\end{proof}

As usual, we denote by $\GMap_{g,n,p}^{\pol}(Z/\GL_N)$ the stack of points of
$(\varphi : \tilde{C} \to C, \sigma^\pm_\bullet, E,s)$ in
$\GMap_{g,n,p}^{\pol}(Z/\GL_N)$ along with the data of a relative polarization
$\cH$ on the family of curves $C$. Our next goal is to state the analog of
\Cref{T:projective_morphism}. The corresponding relatively ample family of line
bundles will be the following.
\begin{defn} \label{defn: LCor for projective targets} Let $(\varphi:
    \widetilde{C} \to C, \sigma_{\bullet}^{\pm}, E, s, \cH) \to
    \GMap_{g,n,p}^{\pol}(Z/\GL_N)$ denote the universal family, and let $\pi:
    \widetilde{C} \to \GMap_{g,n,p}(Z/\GL_N)$ denote the structure morphism. For
    any $m \in \mathbb{Z}$, we denote by $\cL_{\Cor,m}$ the line bundle on the
    stack $\GMap_{g,n,p}(Z/\GL_N)$ defined by
    \begin{gather*}
        \mathcal{L}_{\Cor, m} := \text{det}\left(R\pi_*\left(\left(\left[\omega_{\pi}(\sum_i \sigma_i^{\pm})\otimes \left(s^*(\cO_{E(Z)}(1)) \otimes \text{det}(E)^{\otimes m} \otimes \varphi^*(\cH)^{\otimes m^2}\right)^{\otimes3} \right]- [\mathcal{O}_{C}]\right)^{\otimes 2}\right)\right),
    \end{gather*}
    where we are taking the pushforward and determinant at the level of
    $K$-theory.
\end{defn}

\begin{remark}
    Just as in the case of the corresponding sequence of line bundles in
    $\GBun_{g,n,p}(N)$, the bilinearity of the Deligne pairing (see \cite[Notn.
    2.8]{gauged_theta_stratifications}) implies that the assigment $m \mapsto
    \mathcal{L}_{\Cor, m}$ defines a polynomial on the variable $m$ with values
    on the Picard group of $\GMap_{g,n,p}(Z/\GL_N)$.
\end{remark}

\begin{prop} \label{prop: projective morphism projective target} Fix $d \in
    H_2(Z/\GL_N)$. The morphism $\Phi_Z : \GMap_{g,n,p}^{\pol}(Z/\GL_N,d) \to
    \Pur^{\pol}_{g,n,p}(N)$ is relatively DM and locally projective. For any
    given quasicompact open substack $\cW \subset \Pur^{\pol}_{g,n,p}(N)$, there
    exists some $m\gg 0$ such that the restriction of the line bundle
    $\cL_{\Cor,m}$ to $\Phi_Z^{-1}(\cW)$ is relatively ample on all
    $\cW$-fibers.
\end{prop}

The proof of \Cref{prop: projective morphism projective target} is very similar
to that of \Cref{T:projective_morphism}; we just need to introduce a slightly
different auxiliary stack of Kontsevich stable maps. Let $\cY$ be a locally
Noetherian stack, and let $(C, \overline{\sigma}_{\bullet}^{\pm})$ be a
$\cY$-family of $(n,p)$-marked prestable curves of genus $g$. For any coherent
sheaf $\cF$ on $C$, we recall that there is an associated locally projective
Quot scheme $Q^N_{\cF} \to C$ (\Cref{defn: quot scheme}), which comes equipped
with a universal rank $N$ vector bundle quotient $E_{univ}$. We denote by
$E_{univ}\times^{\GL_N} Z \to Q^N_{\cF}$ the associated $Z$-fiber bundle.

\begin{defn} \label{defn: kontsevich stable maps projective target} The stack
    $\cK_{g,n,p}(E_{univ}(Z)/\cY)$ over $\cY$ is the pseudofunctor that sends a
    $\cY$-scheme $T$ to the groupoid of $T$-families of $(n+p)$-pointed genus
    $g$ Kontsevich stable morphisms $(\tilde{C}, \sigma_{\bullet}^{\pm}) \to
    E_{univ}(Z)\times_{\cY}T$ such that we have an equality of tuples $(\varphi
    \circ \sigma_{\bullet}^{\pm}) = (\overline{\sigma}_{\bullet}^{\pm})$ and the
    composition $\tilde{C} \to C_T$ has degree $1$ on every $T$-fiber.
\end{defn}

\begin{defn}
    Suppose that $C \to \cY$ is equipped with a relative polarization $\cH$. Let
    $(\varphi: \widetilde{C}, C_{\cK_{g,n,p}(E_{univ}(Z)/\cY)}, s: \widetilde{C}
    \to E_{univ}(Z))$ denote the universal family over
    $\cK_{g,n,p}(E_{univ}(Z)/\cY)$, and let $\pi: \widetilde{C} \to
    cK_{g,n,p}(E_{univ}(Z)/\cY)$ denote the structure morphism. For any $m \in
    \mathbb{Z}$, we denote by $\cL_{\Cor,m}$ the line bundle on the stack
    $cK_{g,n,p}(E_{un}(Z)/\cY)$ defined by
    \begin{gather*}
        \mathcal{L}_{\Cor, m} := \text{det}\left(R\pi_*\left(\left(\left[\omega_{\pi}(\sum_i \sigma_i^{\pm})\otimes \left(s^*(\cO_{E_{univ}(Z)}(1)) \otimes \text{det}(E_{univ}|_{\widetilde{C}})^{\otimes m} \otimes \varphi^*(\cH)^{\otimes m^2}\right)^{\otimes3} \right]- [\mathcal{O}_{C}]\right)^{\otimes 2}\right)\right).
    \end{gather*}
\end{defn}

\begin{lemma} \label{lemma: stack kontsevich stable maps into Z is union of
    proper DM stacks} The stack $\cK_{g,n,p}(E_{univ}(Z)/\cY)$ is smooth-locally
    on $\cY$ a disjoint union of $\cY$-proper DM stacks. Suppose furthermore
    that we are given a $\cY$-ample line bundle $\cH$ on $C$, and that $\cY$ is
    quasicompact. Then for all sufficiently large $m$ the line bundle
    $\mathcal{L}_{\Cor, \cH, m}$ restricts to a $\cY$-ample line bundle on any
    proper substack of $\cK_{g,n,p}(E_{univ}(Z)/\cY)$.
\end{lemma}
\begin{proof}
    This follows from a similar argument as in \Cref{lemma: stack kontsevich
    stable maps is union of projective DM stacks}, using the fact that
    $\cO_{E_{univ}(Z)}(1) \otimes \text{det}(E_{univ})^{\otimes m} \otimes
    \cH^{\otimes m^2}$ is relatively ample on $E_{univ}(Z) \to \cY$ for $m\gg0$.
\end{proof}

Let $\cal{F}_{univ}$ be the universal family of pure sheaves on the universal
polarized curve $(C, \cH) \to \Pur^{\pol}_{g,n,p}(N)$. Set $\cK_{g,n,p}(N,Z):=
\cK_{g,n,p}(E_{univ}(Z)/\Pur^{\pol}_{g,n,p}(N))$ as in \Cref{defn: kontsevich
stable maps projective target}.

\begin{lemma} \label{lemma: open immersion of GBun into kontsevich stable maps
    projective target} There is a locally closed immersion
    $\GMap_{g,n,p}^{\pol}(Z/\GL_N) \hookrightarrow \cK_{g,n,p}(N,Z)$ such that
    the composition $\GMap_{g,n,p}^{\pol}(Z/\GL_N)  \hookrightarrow
    \cK_{g,n,p}(N,Z) \to \Pur^{\pol}_{g,n,p}(N)$ is isomorphic to $\Phi_Z$.
\end{lemma}
\begin{proof}
    This follows very similarly as in \Cref{lemma: open immersion of GBun into
    kontsevich stable maps}.
\end{proof}

\begin{prop} \label{prop: valuative criterion over torsion-free moduli for maps
    into Z} The morphism $\Phi_Z : \GMap_{g,n,p}(Z/\GL_N) \to \Pur_{g,n,p}(N)$
    satisfies the valuative criterion for properness for discrete valuation
    rings.
\end{prop}
\begin{proof}
    The proof of this statement is similar to the one in \Cref{L: valuative
    criterion over torsion-free moduli}. We may equivalently consider the
    morphism $\Phi_Z : \GMap_{g,n,p}^{\pol}(Z/\GL_N) \to
    \Pur_{g,n,p}^{\pol}(N)$. We may use the inclusion
    $\GMap_{g,n,p}^{\pol}(Z/\GL_N) \hookrightarrow \cK_{g,n,p}(N,Z)$ to find a
    limit as a morphism to $\cK_{g,n,p}(N,Z)$, since $\cK_{g,n,p}(N,Z) \to
    \Pur^{\pol}_{g,n,p}(N)$ satisfies the valuative criterion for properness
    (\Cref{lemma: stack kontsevich stable maps into Z is union of proper DM
    stacks}). Then one need to check that the resulting limit actually lies in
    the locally closed substack $\GMap_{g,n,p}^{\pol}(Z/\GL_N) \hookrightarrow
    \cK_{g,n,p}(N,Z)$. Again, it suffices to check that there is an induced
    isomorphism of short exact sequences as in diagram \eqref{diagram 1} in the
    proof of \Cref{L: valuative criterion over torsion-free moduli}. This only
    depends on the composition with $\cK_{g,n,p}(N,Z) \to \cK_{g,n,p}(N)$ (the
    additional data of the equivariant morphism to $Z$ is irrelevant), and so
    the same proof as in \Cref{L: valuative criterion over torsion-free moduli}
    applies here.
\end{proof}

\begin{proof}[Proof of \Cref{prop: projective morphism projective target}] As in
the proof of \Cref{T:projective_morphism}, the valuative criterion from
\Cref{prop: valuative criterion over torsion-free moduli for maps into Z}
implies that the morphism $\GMap_{g,n,p}^{\pol}(Z/\GL_N) \hookrightarrow
\cK_{g,n,p}(N,Z)$ defined in \Cref{lemma: open immersion of GBun into kontsevich
stable maps projective target} is a closed immersion. We are only left to show
that $\Phi_Z: \GMap_{g,n,p}^{\pol}(Z/\GL_N,d) \to \Pur^{\pol}_{g,n,p}(N)$ is
quasicompact; the rest of the desired statements will then follow from
\Cref{lemma: stack kontsevich stable maps into Z is union of proper DM stacks}.

    The quasicompactness of $\Phi_Z$ can be checked locally on $S$. For every
    open subscheme $U \subset S$ there is a morphism induced by restriction
    $r_U: H_2(Z_U/\GL_N) \to H_2(Z/\GL_N)$. It suffices to show that the
    morphism $\bigsqcup_{d' \in r_U^{-1}(d)}
    \GMap_{g,n,p}^{\pol}(Z_U/\GL_N)_{d'} \to \Pur^{\pol}_{g,n,p}(N)_U$ is
    quasicompact for any open affine subscheme $U \subset S$. Since $U \times
    \BGL_N$ satisfies the resolution property \cite[Thm. A]{gross-resolution},
    we can use the equivariant polarization $\cO(1)$ to obtain a closed
    immersion $i: Z_U/\GL_N\hookrightarrow \mathbb{P}(V)/\GL_N$ for some
    representation $\GL_N \to \GL(V)$ on some vector bundle $V$ over $U$. The
    morphism $i$ induces a closed immersion $\GMap_{g,n,p}^{\pol}(Z_U/\GL_N)
    \hookrightarrow \GMap_{g,n,p}^{\pol}(\mathbb{P}(V)/\GL_N)$ and a
    homomorphism $\psi: H_2(Z_U/\GL_N) \to H_2(\mathbb{P}(V)/\GL_N)$. It
    suffices to show that $\bigsqcup_{d' \in \psi(r_U^{-1}(d))}
    \GMap_{g,n,p}^{\pol}(\mathbb{P}(V)/\GL_N, d') \to
    (\Pur^{\pol}_{g,n,p}(N))_U$ is quasicompact. Note that
    $N^1_U(\mathbb{P}(V)/\GL_N)$ is generated by classes of line bundles in
    $\BGL_N$ and the tautological equivariant line bundle
    $\cO_{\mathbb{P}(V)}(1)$. Since $\cO(1)$ on $Z/\GL_N$ restricts to the
    pullback of $\cO_{\mathbb{P}(V)}(1)$ on $Z_U/\GL_N$, it follows that
    $\psi(r_U^{-1}(d))$ consists of a single element. Therefore, we can assume
    without loss of generality that $Z=\mathbb{P}(V)$ for some representation
    $\GL_N \to \GL(V)$.

    Let $T$ be a Noetherian scheme, and choose a morphism $T \to
    \Pur^{\pol}_{g,n,p}(N)$ corresponding to a family $(C,
    \overline{\sigma}_{\bullet}^{\pm}, \cH, F, (\sigma^+_\bullet)^\ast(F)
    \twoheadrightarrow W_\bullet)$. Let $E_{univ}$ denote the universal quotient
    vector bundle on the Quot scheme $Q^N_{\cF} \to C \to T$, and let
    $E_{univ}(Z) = \mathbb{P}(E_{univ}(V)) \to Q^N_{\cF}$ denote the associated
    projective bundle. It suffices to show that the image of the open and closed
    immersion $\GMap_{g,n,p}^{\pol}(Z/\GL_N,d) \times_{\Pur^{\pol}_{g,n,p}(N)} T
    \hookrightarrow \cK_{g,n,p}(\mathbb{P}(E_{univ}(V))/T)$ is quasicompact. We
    can define similarly as in \Cref{defn: numerical neron-severi} the numerical
    groups $N^1_T(\mathbb{P}(E_{univ}(V)))$ and $H_2(\mathbb{P}(E_{univ}(V)))$,
    which induce a decomposition
    \[\cK_{g,n,p}(\mathbb{P}(E_{univ}(V))/T) = \bigsqcup_{h \in
    H_2(\mathbb{P}(E_{univ}(V)))}\cK_{g,n,p}(\mathbb{P}(E_{univ}(V))/T, h).\] By
    working locally on $T$ and using boundedness of Kontsevich stable maps with
    fixed degree into projective space, we are reduced to showing that the image
    of
    \[\GMap_{g,n,p}^{\pol}(Z/\GL_N,d) \times_{\Pur^{\pol}_{g,n,p}(N)} T
    \hookrightarrow \cK_{g,n,p}(\mathbb{P}(E_{univ}(V))/T)\] is contained in
    finitely many components $\cK_{g,n,p}(\mathbb{P}(E_{univ}(V))/T, h)$. Let
    $x$ be a geometric point in $\GMap_{g,n,p}^{\pol}(Z/\GL_N,d)
    \times_{\Pur^{\pol}_{g,n,p}(N)} T$ whose first coordinate corresponds to a
    family $(\varphi: \widetilde{C} \to C_x, \sigma_{\bullet}^{\pm}, E, s)$. We
    denote by $h_x \in H_2(\mathbb{P}(E_{univ}(V)))$ the degree of the
    corresponding point in $\cK_{g,n,p}(\mathbb{P}(E_{univ}(V))/T)$. We need to
    show that only finitely many $h_x$ arise. $N^1_T(\mathbb{P}(E_{univ}(V)))$
    is generated by the tautological class $[\cO(1)]$ on the projective bundle
    $\mathbb{P}(E_{univ}(V)) \to Q^N_{F}$ and equivalence classes of line
    bundles that are pulled back from $Q^N_{F}$. The integer $h_x([\cO(1)])$
    agrees with $d([\cO(1)])$, and so it is determined by the fixed degree $d$.
    Therefore, it suffices to show that only finitely many restrictions
    $h_x|_{H_2(Q^N_{F})}$ under the morphism $H_2(Q^N_{F}) \to
    H_2(\mathbb{P}(E_{univ}(V)))$ arise. Let $t$ denote the geometric point
    image of $x$ in the base $T$. Consider the composition $f:(\widetilde{C},
    \sigma_{\bullet}^{\pm}) \to \mathbb{P}(E_{univ}(V))_t \to (Q^N_{\cF})_t$. By
    contracting the $\mathbb{P}^1$-components contracted by $f$, we obtain a
    partial stabilization $(X, \sigma_{\bullet}^{\pm}) \to (Y,
    \tau_{\bullet}^{\pm}) \to (C, \overline{\sigma}_{\bullet}^{\pm})$ and a
    Kontsevich stable morphism $u:(Y, \tau_{\bullet}^{\pm}) \to (Q^N_{F})_t$
    fitting into a commutative diagram
    \begin{equation} \label{equation: diagram 1}
        \begin{tikzcd}[ampersand replacement=\&]
            (X, \sigma_{\bullet}^{\pm}) \ar[r] \ar[d] \& \mathbb{P}(E_{univ}(V))_t  \ar[d] \\  (Y, \tau_{\bullet}^{\pm}) \ar[r] \& (Q^N_{F})_t
        \end{tikzcd}
    \end{equation}

    The morphism $u: (Y, \tau_{\bullet}^{\pm}) \to (Q^N_{F})_t$ yields a
    geometric point $y$ of the substack $\GBun_{g,n,p}^{\pol}(N)
    \times_{\Pur^{\pol}_{g,n,p}(N)} T \hookrightarrow \cK_{g,n,p}(Q^N_{F}/T)$.
    By \Cref{T:projective_morphism}, $\GBun_{g,n,p}^{\pol}(N)
    \times_{\Pur^{\pol}_{g,n,p}(N)} T$ is quasicompact, and so there are
    finitely many possible degrees $h_y \in H_2(Q^N_{F})$ arising from such
    geometric points. By diagram \eqref{equation: diagram 1}, we have
    $h_x|_{H_2(Q^N_{F})}= h_y$, and hence it follows that only finitely many
    restrictions $h_x|_{H_2(Q^N_{F})}$ arise.
\end{proof}

\subsection{The stratification and moduli spaces}
There is an evaluation morphism $\ev^+: \GMap_{g,n}(Z/\GL_N) \to (Z/\GL_N)^p$
which sends a family $(\varphi: \widetilde{C} \to C, \sigma_{\bullet}^{\pm}, u)$
to the $p$-tuple of compositions $\left(T \xrightarrow{\sigma_i} X
\xrightarrow{u} Z/\GL_N\right)_{i=1}^n$. In this subsection, we view
$\GMap_{g,n,p}^{\pol}(Z/\GL_N)$ as a stack over $\prestable^{\pol}_{g,n,p}
\times (Z/\GL_N)^p$ via the natural morphism to $\prestable^{\pol}_{g,n,p}$ that
keeps the target curve, and the evaluation morphism $\ev^+$.

\begin{notn}
    We denote by $\cL_0$ the pullback of the line bundle from \Cref{defn:
    rational line bundle L mu} via the relatively DM morphism $\Phi_Z:
    \GMap_{g,n,p}^{\pol}(Z/\GL_N) \to \Pur_{g,n,p}^{\pol}(N)$. Similarly, we
    denote by $b$ the nondegenerate rational norm on graded points on
    $\GMap_{g,n,p}(Z/\GL_N)$ relative to $\prestable^{\pol}_{g,n,p} \times
    (Z/\GL_N)^p$ obtained by pulling back the corresponding norm on
    $\Pur_{g,n,p}^{\pol}(N)$ (see \Cref{defn: rational norm on Pur}).
\end{notn}

For the following definition, recall the polynomial sequence of line bundles
$\cL_{\Cor,m}$ on $ \GMap_{g,n,p}^{\pol}(Z/\GL_N)$ from \Cref{defn: LCor for
projective targets}.
\begin{defn} \label{defn: numerical invariant projective target} The
    formal polynomial numerical invariant $\nu_{\epsilon}$ on
    $\GMap_{g,n,p}^{\pol}(Z/\GL_N)$ relative to $\prestable^{\pol}_{g,n,p}
    \times (Z/\GL_N)^p$ is defined to be $\left(\wt(\cL_{0}) - \epsilon^5 \cdot
    \wt(\cL_{\Cor,\epsilon^{-1}})\right) / \sqrt{b}$.
\end{defn}

\begin{thm} \label{thm: theta stratification for projective target} Fix $d \in
    H_2(Z/\GL_N)$. The formal numerical invariant $\nu_{\epsilon}$ defines a
    weak $\Theta$-stratification of $\GMap_{g,n,p}^{\pol}(Z/\GL_N,d)$ relative
    to $\prestable^{\pol}_{g,n,p} \times (Z/\GL_N)^p$ satisfying the following:
    \begin{enumerate}
        \item Every open stratum $\GMap_{g,n,p}^{\pol}(Z/\GL_N,d)_{\nu_\epsilon
        \leq \gamma}$ is of finite type over $\prestable^{\pol}_{g,n,p} \times
        (Z/\GL_N)^p$. In particular, the relative weak $\Theta$-stratification
        is well-ordered.
        \item If the base scheme $S$ is a $\mathbb{Q}$-scheme, then
              \begin{enumerate}
                  \item The weak $\Theta$-stratification is a
                  $\Theta$-stratification;
                  \item The semistable locus
                  $\GMap_{g,n,p}^{\pol}(Z/\GL_N,d)^{\nu_{\epsilon}\dash ss}$
                  admits a proper relative good moduli space over
                  $\prestable^{\pol}_{g,n,p} \times (Z/\GL_N)^p$;
                  \item The center of every $\nu_\epsilon$-unstable stratum
                  admits another $\Theta$-stratification whose centers have
                  proper good moduli spaces relative to
                  $\prestable^{\pol}_{g,n,p} \times (Z/\GL_N)^p$.
              \end{enumerate}
    \end{enumerate}
\end{thm}
\begin{proof}
    The same formal proof as for \Cref{thm: theta stratification for Gieseker
    bundles} applies here, we only need to replace the use of
    \Cref{T:projective_morphism} with \Cref{prop: projective morphism projective
    target}.
\end{proof}

\printendnotes

\footnotesize{\bibliography{compactification_via_bubblings.bib}}

@Misc{stacks-project,
  author       = {The {Stacks Project Authors}},
  howpublished = {\url{https://stacks.math.columbia.edu}},
  title        = {\textit{Stacks Project}},
  year         = {2024},
  shorthand    = {Stacks},
}

@misc{halpernleistner2025categoricalperspectivenonabelianlocalization,
      title={A categorical perspective on non-abelian localization}, 
      author={Daniel Halpern-Leistner},
      year={2025},
      eprint={2509.24009},
      archivePrefix={arXiv},
      primaryClass={math.AG},
      url={https://arxiv.org/abs/2509.24009}, 
}

@article {alper2019existence,
    AUTHOR = {Alper, Jarod and Halpern-Leistner, Daniel and Heinloth,
              Jochen},
     TITLE = {Existence of moduli spaces for algebraic stacks},
   JOURNAL = {Invent. Math.},
  FJOURNAL = {Inventiones Mathematicae},
    VOLUME = {234},
      YEAR = {2023},
    NUMBER = {3},
     PAGES = {949--1038},
      ISSN = {0020-9910,1432-1297},
   MRCLASS = {14D23 (14A20)},
  MRNUMBER = {4665776},
       DOI = {10.1007/s00222-023-01214-4},
       URL = {https://doi.org/10.1007/s00222-023-01214-4},
}

@Misc{halpernleistner2018structure,
  author        = {Daniel Halpern-Leistner},
  howpublished  = {\url{https://arxiv.org/abs/1411.0627}},
  title         = {On the structure of instability in moduli theory},
  year          = {2021},
  note = {arXiv preprint},
  archiveprefix = {arXiv},
  eprint        = {1411.0627},
  primaryclass  = {math.AG},
}

@article {egaiv,
    AUTHOR = {Grothendieck, A.},
     TITLE = {\'{E}l\'{e}ments de g\'{e}om\'{e}trie alg\'{e}brique. {IV}. \'{E}tude locale des
              sch\'{e}mas et des morphismes de sch\'{e}mas. {III}},
   JOURNAL = {Inst. Hautes \'{E}tudes Sci. Publ. Math.},
  FJOURNAL = {Institut des Hautes \'{E}tudes Scientifiques. Publications
              Math\'{e}matiques},
    VOLUME = {28},
      YEAR = {1966},
     PAGES = {255},
      ISSN = {0073-8301},
   MRCLASS = {14.55},
  MRNUMBER = {217086},
MRREVIEWER = {J. P. Murre},
       URL = {http://www.numdam.org.proxy.library.cornell.edu/item?id=PMIHES_1966__28__255_0},
}

@article {alper-good-moduli,
    AUTHOR = {Alper, Jarod},
     TITLE = {Good moduli spaces for {A}rtin stacks},
   JOURNAL = {Ann. Inst. Fourier (Grenoble)},
  FJOURNAL = {Universit\'{e} de Grenoble. Annales de l'Institut Fourier},
    VOLUME = {63},
      YEAR = {2013},
    NUMBER = {6},
     PAGES = {2349--2402},
      ISSN = {0373-0956},
   MRCLASS = {14D23 (14L24 14L30)},
  MRNUMBER = {3237451},
MRREVIEWER = {Arvid Siqveland},
       URL = {http://aif.cedram.org/item?id=AIF_2013__63_6_2349_0},
}

@article {torsion-freepaper,
    AUTHOR = {Halpern-Leistner, Daniel and Fernandez Herrero, Andres and
              Jones, Trevor},
     TITLE = {Moduli spaces of sheaves via affine {G}rassmannians},
   JOURNAL = {J. Reine Angew. Math.},
  FJOURNAL = {Journal f\"ur die Reine und Angewandte Mathematik. [Crelle's
              Journal]},
    VOLUME = {809},
      YEAR = {2024},
     PAGES = {159--215},
      ISSN = {0075-4102,1435-5345},
   MRCLASS = {14D23 (14L24 14M15)},
  MRNUMBER = {4726568},
       DOI = {10.1515/crelle-2023-0099},
       URL = {https://doi.org/10.1515/crelle-2023-0099},
}

@article {rho-sheaves-paper,
    AUTHOR = {G\'omez, Tom\'as L. and Fernandez Herrero, Andres and Zamora,
              Alfonso},
     TITLE = {The moduli stack of principal {$\rho$}-sheaves and
              {G}ieseker-{H}arder-{N}arasimhan filtrations},
   JOURNAL = {Math. Z.},
  FJOURNAL = {Mathematische Zeitschrift},
    VOLUME = {307},
      YEAR = {2024},
    NUMBER = {3},
     PAGES = {Paper No. 51, 67},
      ISSN = {0025-5874,1432-1823},
   MRCLASS = {14D20 (14D23 14F06 14J60 14L24)},
  MRNUMBER = {4756614},
       DOI = {10.1007/s00209-024-03497-6},
       URL = {https://doi.org/10.1007/s00209-024-03497-6},
}

@misc{gauged_theta_stratifications,
 title= {The structure of the moduli of gauged maps from a smooth curve},
      author={Daniel Halpern-Leistner and Andres Fernandez Herrero},
      year={2023},
      howpublished = {\url{https://arxiv.org/abs/2305.09632}},
}

@article {deligne_mumford,
    AUTHOR = {Deligne, P. and Mumford, D.},
     TITLE = {The irreducibility of the space of curves of given genus},
   JOURNAL = {Inst. Hautes \'{E}tudes Sci. Publ. Math.},
  FJOURNAL = {Institut des Hautes \'{E}tudes Scientifiques. Publications
              Math\'{e}matiques},
    VOLUME = {36},
      YEAR = {1969},
     PAGES = {75--109},
      ISSN = {0073-8301},
   MRCLASS = {14.20},
  MRNUMBER = {262240},
MRREVIEWER = {Manfred Herrmann},
       URL = {http://www.numdam.org.proxy.library.cornell.edu/item?id=PMIHES_1969__36__75_0},
}

@article {behrend-manin,
    AUTHOR = {Behrend, K. and Manin, Yu.},
     TITLE = {Stacks of stable maps and {G}romov-{W}itten invariants},
   JOURNAL = {Duke Math. J.},
  FJOURNAL = {Duke Mathematical Journal},
    VOLUME = {85},
      YEAR = {1996},
    NUMBER = {1},
     PAGES = {1--60},
      ISSN = {0012-7094},
   MRCLASS = {14D20 (14C25 14D22)},
  MRNUMBER = {1412436},
MRREVIEWER = {Barbara Fantechi},
       DOI = {10.1215/S0012-7094-96-08501-4},
       URL = {https://doi-org.proxy.library.cornell.edu/10.1215/S0012-7094-96-08501-4},
}

@article {hall-rydh-tannakahom,
    AUTHOR = {Hall, Jack and Rydh, David},
     TITLE = {Coherent {T}annaka duality and algebraicity of {H}om-stacks},
   JOURNAL = {Algebra Number Theory},
  FJOURNAL = {Algebra \& Number Theory},
    VOLUME = {13},
      YEAR = {2019},
    NUMBER = {7},
     PAGES = {1633--1675},
      ISSN = {1937-0652},
   MRCLASS = {14A20 (14D23 18M05)},
  MRNUMBER = {4009673},
MRREVIEWER = {Frank Gounelas},
       DOI = {10.2140/ant.2019.13.1633},
       URL = {https://doi.org/10.2140/ant.2019.13.1633},
}

@article {halpernleistner2019mapping,
    AUTHOR = {Halpern-Leistner, Daniel and Preygel, Anatoly},
     TITLE = {Mapping stacks and categorical notions of properness},
   JOURNAL = {Compos. Math.},
  FJOURNAL = {Compositio Mathematica},
    VOLUME = {159},
      YEAR = {2023},
    NUMBER = {3},
     PAGES = {530--589},
      ISSN = {0010-437X,1570-5846},
   MRCLASS = {14A20 (14A30)},
  MRNUMBER = {4560539},
MRREVIEWER = {Jon\ Eivind\ Vatne},
       DOI = {10.1112/S0010437X22007667},
       URL = {https://doi.org/10.1112/S0010437X22007667},
}

@article {nagaraj-seshadri-2,
    AUTHOR = {Nagaraj, D. S. and Seshadri, C. S.},
     TITLE = {Degenerations of the moduli spaces of vector bundles on
              curves. {II}. {G}eneralized {G}ieseker moduli spaces},
   JOURNAL = {Proc. Indian Acad. Sci. Math. Sci.},
  FJOURNAL = {Indian Academy of Sciences. Proceedings. Mathematical
              Sciences},
    VOLUME = {109},
      YEAR = {1999},
    NUMBER = {2},
     PAGES = {165--201},
      ISSN = {0253-4142},
   MRCLASS = {14H60 (14D20 14D22)},
  MRNUMBER = {1687729},
MRREVIEWER = {Francesco Bottacin},
       DOI = {10.1007/BF02841533},
       URL = {https://doi-org.proxy.library.cornell.edu/10.1007/BF02841533},
}

@article {teleman-woodward,
    AUTHOR = {Teleman, Constantin and Woodward, Christopher T.},
     TITLE = {The index formula for the moduli of {$G$}-bundles on a curve},
   JOURNAL = {Ann. of Math. (2)},
  FJOURNAL = {Annals of Mathematics. Second Series},
    VOLUME = {170},
      YEAR = {2009},
    NUMBER = {2},
     PAGES = {495--527},
      ISSN = {0003-486X,1939-8980},
   MRCLASS = {58J22 (14D23 19L50 57R56)},
  MRNUMBER = {2552100},
MRREVIEWER = {Daniel\ S.\ Freed},
       DOI = {10.4007/annals.2009.170.495},
       URL = {https://doi-org.proxy.library.upenn.edu/10.4007/annals.2009.170.495},
}

@article {thomas_casson_invariants,
    AUTHOR = {Thomas, R. P.},
     TITLE = {A holomorphic {C}asson invariant for {C}alabi-{Y}au 3-folds,
              and bundles on {$K3$} fibrations},
   JOURNAL = {J. Differential Geom.},
  FJOURNAL = {Journal of Differential Geometry},
    VOLUME = {54},
      YEAR = {2000},
    NUMBER = {2},
     PAGES = {367--438},
      ISSN = {0022-040X,1945-743X},
   MRCLASS = {14J32 (14J60 32J17 32Q25)},
  MRNUMBER = {1818182},
MRREVIEWER = {Bal\'azs\ Szendr\H oi},
       URL = {http://projecteuclid.org.proxy.library.upenn.edu/euclid.jdg/1214341649},
}

@article {kontsevich-manin,
    AUTHOR = {Kontsevich, M. and Manin, Yu.},
     TITLE = {Gromov-{W}itten classes, quantum cohomology, and enumerative
              geometry},
   JOURNAL = {Comm. Math. Phys.},
  FJOURNAL = {Communications in Mathematical Physics},
    VOLUME = {164},
      YEAR = {1994},
    NUMBER = {3},
     PAGES = {525--562},
      ISSN = {0010-3616,1432-0916},
   MRCLASS = {14N10 (53C15 58D10 58F05)},
  MRNUMBER = {1291244},
MRREVIEWER = {Dietmar\ A.\ Salamon},
       URL = {http://projecteuclid.org.proxy.library.upenn.edu/euclid.cmp/1104270948},
}

@article {dolce_deligne_pairing,
    AUTHOR = {Dolce, Paolo},
     TITLE = {Explicit {D}eligne pairing},
   JOURNAL = {Eur. J. Math.},
  FJOURNAL = {European Journal of Mathematics},
    VOLUME = {8},
      YEAR = {2022},
     PAGES = {S101--S129},
      ISSN = {2199-675X,2199-6768},
   MRCLASS = {14C17},
  MRNUMBER = {4452839},
MRREVIEWER = {P.\ E.\ Newstead},
       DOI = {10.1007/s40879-021-00482-9},
       URL = {https://doi-org.proxy.library.upenn.edu/10.1007/s40879-021-00482-9},
}

@incollection {donaldson_thomas_gauged,
    AUTHOR = {Donaldson, S. K. and Thomas, R. P.},
     TITLE = {Gauge theory in higher dimensions},
 BOOKTITLE = {The geometric universe ({O}xford, 1996)},
     PAGES = {31--47},
 PUBLISHER = {Oxford Univ. Press, Oxford},
      YEAR = {1998},
      ISBN = {0-19-850059-9},
   MRCLASS = {57R57 (14J32 32J18 53C07 57R58 58D27)},
  MRNUMBER = {1634503},
MRREVIEWER = {Krzysztof\ Galicki},
}

@article {frenkel-teleman-tolland-ggw,
    AUTHOR = {Frenkel, Edward and Teleman, Constantin and Tolland, A. J.},
     TITLE = {Gromov-{W}itten gauge theory},
   JOURNAL = {Adv. Math.},
  FJOURNAL = {Advances in Mathematics},
    VOLUME = {288},
      YEAR = {2016},
     PAGES = {201--239},
      ISSN = {0001-8708},
   MRCLASS = {14N35 (14A20 53D45)},
  MRNUMBER = {3436385},
MRREVIEWER = {Sergiy Koshkin},
       DOI = {10.1016/j.aim.2015.10.008},
       URL = {https://doi-org.proxy.library.cornell.edu/10.1016/j.aim.2015.10.008},
}

@article {schmitt-hilbert-compactification,
    AUTHOR = {Schmitt, Alexander},
     TITLE = {The {H}ilbert compactification of the universal moduli space
              of semistable vector bundles over smooth curves},
   JOURNAL = {J. Differential Geom.},
  FJOURNAL = {Journal of Differential Geometry},
    VOLUME = {66},
      YEAR = {2004},
    NUMBER = {2},
     PAGES = {169--209},
      ISSN = {0022-040X},
   MRCLASS = {14H60 (14D20 14H10)},
  MRNUMBER = {2106123},
MRREVIEWER = {Arvid Siqveland},
       URL = {http://projecteuclid.org.proxy.library.cornell.edu/euclid.jdg/1102538609},
}

@article {behrend-gromov-witten-theory,
    AUTHOR = {Behrend, K.},
     TITLE = {Gromov-{W}itten invariants in algebraic geometry},
   JOURNAL = {Invent. Math.},
  FJOURNAL = {Inventiones Mathematicae},
    VOLUME = {127},
      YEAR = {1997},
    NUMBER = {3},
     PAGES = {601--617},
      ISSN = {0020-9910},
   MRCLASS = {14D20 (14C25 14D22)},
  MRNUMBER = {1431140},
MRREVIEWER = {Barbara Fantechi},
       DOI = {10.1007/s002220050132},
       URL = {https://doi.org/10.1007/s002220050132},
}

@article {starr-dejong-amost-proper,
    AUTHOR = {Starr, Jason and de Jong, Johan},
     TITLE = {Almost proper {GIT}-stacks and discriminant avoidance},
   JOURNAL = {Doc. Math.},
  FJOURNAL = {Documenta Mathematica},
    VOLUME = {15},
      YEAR = {2010},
     PAGES = {957--972},
      ISSN = {1431-0635},
   MRCLASS = {14J10 (14L15)},
  MRNUMBER = {2745688},
}

@incollection {hoffmann-ks-rationality,
    AUTHOR = {Hoffmann, Norbert},
     TITLE = {Moduli stacks of vector bundles on curves and the
              {K}ing-{S}chofield rationality proof},
 BOOKTITLE = {Cohomological and geometric approaches to rationality
              problems},
    SERIES = {Progr. Math.},
    VOLUME = {282},
     PAGES = {133--148},
 PUBLISHER = {Birkh\"{a}user Boston, Boston, MA},
      YEAR = {2010},
   MRCLASS = {14H60 (14D23 14E08)},
  MRNUMBER = {2605167},
MRREVIEWER = {P. E. Newstead},
       DOI = {10.1007/978-0-8176-4934-0\_5},
       URL = {https://doi.org/10.1007/978-0-8176-4934-0_5},
}

@article {langer-boundedness,
    AUTHOR = {Langer, Adrian},
     TITLE = {Semistable sheaves in positive characteristic},
   JOURNAL = {Ann. of Math. (2)},
  FJOURNAL = {Annals of Mathematics. Second Series},
    VOLUME = {159},
      YEAR = {2004},
    NUMBER = {1},
     PAGES = {251--276},
      ISSN = {0003-486X},
   MRCLASS = {14F05 (14D20 14J60)},
  MRNUMBER = {2051393},
MRREVIEWER = {Vikram B. Mehta},
       DOI = {10.4007/annals.2004.159.251},
       URL = {https://doi.org/10.4007/annals.2004.159.251},
}

@article {gross-resolution,
    AUTHOR = {Gross, Philipp},
     TITLE = {Tensor generators on schemes and stacks},
   JOURNAL = {Algebr. Geom.},
  FJOURNAL = {Algebraic Geometry},
    VOLUME = {4},
      YEAR = {2017},
    NUMBER = {4},
     PAGES = {501--522},
      ISSN = {2313-1691},
   MRCLASS = {14A20 (14F05 14L30)},
  MRNUMBER = {3683505},
MRREVIEWER = {XiaoWen Hu},
       DOI = {10.14231/2017-026},
       URL = {https://doi.org/10.14231/2017-026},
}

@article {romagny-composantes,
    AUTHOR = {Romagny, Matthieu},
     TITLE = {Composantes connexes et irr\'{e}ductibles en familles},
   JOURNAL = {Manuscripta Math.},
  FJOURNAL = {Manuscripta Mathematica},
    VOLUME = {136},
      YEAR = {2011},
    NUMBER = {1-2},
     PAGES = {1--32},
      ISSN = {0025-2611,1432-1785},
   MRCLASS = {14D23 (14D06 14H10)},
  MRNUMBER = {2820394},
MRREVIEWER = {Orsola\ Tommasi},
       DOI = {10.1007/s00229-010-0424-7},
       URL = {https://doi.org/10.1007/s00229-010-0424-7},
}

@misc{martens_thaddeus_variations_grothendieck,
      title={Variations on a theme of {G}rothendieck}, 
      author={Johan Martens and Michael Thaddeus},
      year={2015},
      eprint={1210.8161},
      archivePrefix={arXiv},
      primaryClass={math.AG},
      url={https://arxiv.org/abs/1210.8161}, 
howpublished = {\url{https://arxiv.org/abs/1210.8161}},
}

@article {lee_quantum_k_theory_i,
    AUTHOR = {Lee, Y.-P.},
     TITLE = {Quantum {$K$}-theory. {I}. {F}oundations},
   JOURNAL = {Duke Math. J.},
  FJOURNAL = {Duke Mathematical Journal},
    VOLUME = {121},
      YEAR = {2004},
    NUMBER = {3},
     PAGES = {389--424},
      ISSN = {0012-7094,1547-7398},
   MRCLASS = {14N35 (19E08 53D45 55N15)},
  MRNUMBER = {2040281},
MRREVIEWER = {Andrew\ Kresch},
       DOI = {10.1215/S0012-7094-04-12131-1},
       URL = {https://doi.org/10.1215/S0012-7094-04-12131-1},
}

@incollection {kausz_stable_maps,
    AUTHOR = {Kausz, Ivan},
     TITLE = {Stable maps into the classifying space of the general linear
              group},
 BOOKTITLE = {Teichm\"uller theory and moduli problem},
    SERIES = {Ramanujan Math. Soc. Lect. Notes Ser.},
    VOLUME = {10},
     PAGES = {437--449},
 PUBLISHER = {Ramanujan Math. Soc., Mysore},
      YEAR = {2010},
      ISBN = {978-93-80416-00-7},
   MRCLASS = {14D23},
  MRNUMBER = {2667565},
MRREVIEWER = {Cristina\ Mart\'inez},
}

@article {caporaso_compactification,
    AUTHOR = {Caporaso, Lucia},
     TITLE = {A compactification of the universal {P}icard variety over the
              moduli space of stable curves},
   JOURNAL = {J. Amer. Math. Soc.},
  FJOURNAL = {Journal of the American Mathematical Society},
    VOLUME = {7},
      YEAR = {1994},
    NUMBER = {3},
     PAGES = {589--660},
      ISSN = {0894-0347,1088-6834},
   MRCLASS = {14D20 (14D25 14H10 14H40)},
  MRNUMBER = {1254134},
MRREVIEWER = {P.\ E.\ Newstead},
       DOI = {10.2307/2152786},
       URL = {https://doi.org/10.2307/2152786},
}

@article {gieseker-degeneration,
    AUTHOR = {Gieseker, D.},
     TITLE = {A degeneration of the moduli space of stable bundles},
   JOURNAL = {J. Differential Geom.},
  FJOURNAL = {Journal of Differential Geometry},
    VOLUME = {19},
      YEAR = {1984},
    NUMBER = {1},
     PAGES = {173--206},
      ISSN = {0022-040X,1945-743X},
   MRCLASS = {14D20 (14D22 14F05)},
  MRNUMBER = {739786},
MRREVIEWER = {P.\ E.\ Newstead},
       URL = {http://projecteuclid.org/euclid.jdg/1214438427},
}

@article {pandharipande-compactification,
    AUTHOR = {Pandharipande, Rahul},
     TITLE = {A compactification over {$\overline {M}_g$} of the universal
              moduli space of slope-semistable vector bundles},
   JOURNAL = {J. Amer. Math. Soc.},
  FJOURNAL = {Journal of the American Mathematical Society},
    VOLUME = {9},
      YEAR = {1996},
    NUMBER = {2},
     PAGES = {425--471},
      ISSN = {0894-0347,1088-6834},
   MRCLASS = {14D20 (14D25 14H10 14H60)},
  MRNUMBER = {1308406},
MRREVIEWER = {P.\ E.\ Newstead},
       DOI = {10.1090/S0894-0347-96-00173-7},
       URL = {https://doi.org/10.1090/S0894-0347-96-00173-7},
}

@article {joyce_song,
    AUTHOR = {Joyce, Dominic and Song, Yinan},
     TITLE = {A theory of generalized {D}onaldson-{T}homas invariants},
   JOURNAL = {Mem. Amer. Math. Soc.},
  FJOURNAL = {Memoirs of the American Mathematical Society},
    VOLUME = {217},
      YEAR = {2012},
    NUMBER = {1020},
     PAGES = {iv+199},
      ISSN = {0065-9266,1947-6221},
      ISBN = {978-0-8218-5279-8},
   MRCLASS = {14N35 (14D23 14F05 14J32)},
  MRNUMBER = {2951762},
MRREVIEWER = {Amin\ Gholampour},
       DOI = {10.1090/S0065-9266-2011-00630-1},
       URL = {https://doi.org/10.1090/S0065-9266-2011-00630-1},
}

@incollection {kontsevich_soibelman,
    AUTHOR = {Kontsevich, Maxim and Soibelman, Yan},
     TITLE = {Wall-crossing structures in {D}onaldson-{T}homas invariants,
              integrable systems and mirror symmetry},
 BOOKTITLE = {Homological mirror symmetry and tropical geometry},
    SERIES = {Lect. Notes Unione Mat. Ital.},
    VOLUME = {15},
     PAGES = {197--308},
 PUBLISHER = {Springer, Cham},
      YEAR = {2014},
      ISBN = {978-3-319-06513-7; 978-3-319-06514-4},
   MRCLASS = {14N35 (14J33 53D37)},
  MRNUMBER = {3330788},
MRREVIEWER = {Victor\ Przyjalkowski},
       DOI = {10.1007/978-3-319-06514-4\_6},
       URL = {https://doi.org/10.1007/978-3-319-06514-4_6},
}

@article {kontsevich_soibelman_2,
    AUTHOR = {Kontsevich, Maxim and Soibelman, Yan},
     TITLE = {Cohomological {H}all algebra, exponential {H}odge structures
              and motivic {D}onaldson-{T}homas invariants},
   JOURNAL = {Commun. Number Theory Phys.},
  FJOURNAL = {Communications in Number Theory and Physics},
    VOLUME = {5},
      YEAR = {2011},
    NUMBER = {2},
     PAGES = {231--352},
      ISSN = {1931-4523,1931-4531},
   MRCLASS = {14N35 (14F43 16G20)},
  MRNUMBER = {2851153},
MRREVIEWER = {Mark\ Gross},
       DOI = {10.4310/CNTP.2011.v5.n2.a1},
       URL = {https://doi.org/10.4310/CNTP.2011.v5.n2.a1},
}

@misc {hl_ibanez,
author = {Halpern-Leistner, Daniel and Ib{\'a}{\~n}ez N{\'u}{\~n}ez, Andr{\'e}s},
title = {Spaces of {$\Theta$}-stratifications},
year = {2026},
howpublished = {In preparation},

}

@incollection {okounkovlectures,
    AUTHOR = {Okounkov, Andrei},
     TITLE = {Lectures on {K}-theoretic computations in enumerative
              geometry},
 BOOKTITLE = {Geometry of moduli spaces and representation theory},
    SERIES = {IAS/Park City Math. Ser.},
    VOLUME = {24},
     PAGES = {251--380},
 PUBLISHER = {Amer. Math. Soc., Providence, RI},
      YEAR = {2017},
      ISBN = {978-1-4704-3574-5},
   MRCLASS = {14N35 (14C35 19M05)},
  MRNUMBER = {3752463},
MRREVIEWER = {Emily\ Clader},
       DOI = {10.1090/pcms/024},
       URL = {https://doi.org/10.1090/pcms/024},
}

@incollection {gonzalez-solis-woodward-stable-gauged,
    AUTHOR = {Gonz\'{a}lez, Eduardo and Solis, Pablo and Woodward, Chris T.},
     TITLE = {Stable gauged maps},
 BOOKTITLE = {Algebraic geometry: {S}alt {L}ake {C}ity 2015},
    SERIES = {Proc. Sympos. Pure Math.},
    VOLUME = {97},
     PAGES = {243--275},
 PUBLISHER = {Amer. Math. Soc., Providence, RI},
      YEAR = {2018},
   MRCLASS = {14N35 (14D23)},
  MRNUMBER = {3821152},
MRREVIEWER = {Hsian-Hua Tseng},
       DOI = {10.1016/j.jalgebra.2017.06.015},
       URL = {https://doi.org/10.1016/j.jalgebra.2017.06.015},
}

@incollection {gonzalez-woodward-quantum,
    AUTHOR = {Gonz\'{a}lez, E. and Woodward, C.},
     TITLE = {Quantum {K}irwan for quantum {K}-theory},
 BOOKTITLE = {Facets of algebraic geometry. {V}ol. {I}},
    SERIES = {London Math. Soc. Lecture Note Ser.},
    VOLUME = {472},
     PAGES = {265--332},
 PUBLISHER = {Cambridge Univ. Press, Cambridge},
      YEAR = {2022},
   MRCLASS = {14C35 (14L24 19M99)},
  MRNUMBER = {4381905},
MRREVIEWER = {Cheolgyu Lee},
}

@article {gonzalez-solis-woodward-properness,
    AUTHOR = {Gonz\'{a}lez, Eduardo and Solis, Pablo and Woodward, Chris T.},
     TITLE = {Properness for scaled gauged maps},
   JOURNAL = {J. Algebra},
  FJOURNAL = {Journal of Algebra},
    VOLUME = {490},
      YEAR = {2017},
     PAGES = {104--157},
      ISSN = {0021-8693},
   MRCLASS = {14N35 (14D23 14L24 53D45)},
  MRNUMBER = {3690329},
MRREVIEWER = {Xiaobin Li},
       DOI = {10.1016/j.jalgebra.2017.06.015},
       URL = {https://doi.org/10.1016/j.jalgebra.2017.06.015},
}

@article {woodward-quantum-quotients,
    AUTHOR = {Woodward, Chris T.},
     TITLE = {Quantum {K}irwan morphism and {G}romov-{W}itten invariants of
              quotients {I}},
   JOURNAL = {Transform. Groups},
  FJOURNAL = {Transformation Groups},
    VOLUME = {20},
      YEAR = {2015},
    NUMBER = {2},
     PAGES = {507--556},
      ISSN = {1083-4362},
   MRCLASS = {14N35 (55N32)},
  MRNUMBER = {3348566},
MRREVIEWER = {Xiaobin Li},
       DOI = {10.1007/s00031-015-9313-1},
       URL = {https://doi.org/10.1007/s00031-015-9313-1},
}

@misc{halpernleistner2016equivariantverlindeformulamoduli,
      title={The equivariant Verlinde formula on the moduli of Higgs bundles}, 
      author={Daniel Halpern-Leistner},
      year={2016},
      eprint={1608.01754},
      archivePrefix={arXiv},
      primaryClass={math.AG},
      url={https://arxiv.org/abs/1608.01754}, 
    howpublished = {\url{https://arxiv.org/abs/1608.01754}},
}

@Misc{cooper_compactified_jac,
  author        = {George Cooper},
  howpublished  = {\url{https://arxiv.org/abs/2210.11457}},
  title         = {GIT Constructions of Compactified Universal Jacobians over Stacks of Stable Maps},
  year          = {2022},
  note = {arXiv preprint},
}

@misc{donagi2024meromorphichitchinfibrationstable,
      title={The meromorphic {H}itchin fibration over stable pointed curves: moduli spaces}, 
      author={Ron Donagi and Andres Fernandez Herrero},
    howpublished = {\url{https://arxiv.org/abs/2411.16912}},
      year={2024},
      eprint={2411.16912},
      archivePrefix={arXiv},
      primaryClass={math.AG},
      url={https://arxiv.org/abs/2411.16912}, 
}
\bibliographystyle{alpha}
\end{document}